\documentclass[opre,nonblindrev]{informs3}

\usepackage{pgfplots}
\usepackage[english]{babel}
\usepackage[autostyle, english = american]{csquotes}
\MakeOuterQuote{"}
\usepackage{graphicx}
\usepackage{subcaption}
\usepackage{bbm,xspace,multirow,multicol}
\usepackage{caption}
\usepackage{amsmath}
\usepackage{bm}
\usepackage{array}
\newcolumntype{P}[1]{>{\centering\arraybackslash}p{#1}}
\usepackage[normalem]{ulem}
\usepackage{pdflscape}
\usepackage{cancel}
\usepackage{mathrsfs}
\usepackage{enumerate}
\usepackage[hypertexnames=false]{hyperref}
\usepackage{cleveref}
\crefname{subsection}{subsection}{subsections}

\newcommand{\eps}{\varepsilon}

\newcommand{\bR}{\mathbb{R}}

\newcommand{\FLU}{\mathsf{FLU}}
\newcommand{\pr}{\mathbb{P}}

\NewDocumentEnvironment{myproof}{o}
{\IfNoValueTF{#1}{\paragraph{{Proof.} }} {\paragraph{{#1.} }} }
{\hfill$\Halmos$}

\usepackage{natbib}
 \bibpunct[, ]{(}{)}{,}{a}{}{,}%
 \def\bibfont{\small}%
\TheoremsNumberedThrough     
\ECRepeatTheorems

\EquationsNumberedThrough    

\MANUSCRIPTNO{} 

\begin{document}



\RUNAUTHOR{Balseiro, Ma, and Zhang}

\RUNTITLE{Dynamic Pricing for Reusable Resources}

\TITLE{Dynamic Pricing for Reusable Resources:\\The Power of Two Prices}

\ARTICLEAUTHORS{%
\AUTHOR{Santiago R. Balseiro}
\AFF{Graduate School of Business, Columbia University, \EMAIL{srb2155@columbia.edu}}

\AUTHOR{Will Ma}
\AFF{Graduate School of Business, Columbia University, \EMAIL{wm2428@columbia.edu}}

\AUTHOR{Wenxin Zhang}
\AFF{Graduate School of Business, Columbia University, \EMAIL{wz2574@columbia.edu}}
}
\ABSTRACT{  
\indent 

Motivated by real-world applications such as rental and cloud computing services, we investigate pricing for reusable resources. We consider a system where a single resource with a fixed number of identical copies serves customers with heterogeneous willingness-to-pay (WTP), and the usage duration distribution is general. Optimal dynamic policies are computationally intractable when usage durations are not memoryless, so existing literature has focused on static pricing, which incurs a steady-state performance loss of $\mathcal{O}(\sqrt{c})$ compared to optimality when supply and demand scale with $c$. We propose a class of dynamic "stock-dependent" policies that 1) are computationally tractable and 2) can attain a steady-state performance loss of $o(\sqrt{c})$. We give parametric bounds based on the local shape of the reward function at the optimal fluid admission probability and show that the performance loss of stock-dependent policies can be as low as $\mathcal{O}((\log{c})^2)$. We characterize the tight performance loss for stock-dependent policies and show that they can in fact be achieved by a simple two-price policy that sets a higher price when the stock is below some threshold and a lower price otherwise. We extend our results to settings with multiple resources and multiple customer classes. Finally, we demonstrate this "minimally dynamic" class of two-price policies performs well numerically, even in non-asymptotic settings, suggesting that a little dynamicity can go a long way.

}

\KEYWORDS{reusable resources; dynamic pricing; online algorithms; loss network; asymptotic analysis}

\maketitle

\section{Introduction}\label{sec: intro}

In many real-world applications, companies manage resources that are reusable, i.e., the same resource can be repeatedly used to serve customers after the previous service finishes. Prominent examples are car rental, hotel booking, parking lots, and cloud computing services.\footnote{According to \href{https://www.gartner.com/en/newsroom/press-releases/2022-10-31-gartner-forecasts-worldwide-public-cloud-end-user-spending-to-reach-nearly-600-billion-in-2023}{Gartner}, the 2023 market size of the worldwide public cloud services is estimated to be close to 600 billion dollars, and this market is expected to grow at an annual rate of 20\%.} In these markets, customers request resources (rooms in the case of hotel booking or CPU and storage in the case of cloud computing) and the usage duration can be deterministic or random. An important operational question is how to manage reusable resources efficiently. 
The challenge with reusable resources is that the number of system states can grow exponentially with the initial stock, i.e., the total number of units in the system. When the usage duration is memoryless, tracking the number of available units of stock is sufficient, but more is needed when studying general usage durations. In particular, a Markovian state descriptor must also include information regarding how long each unavailable unit has been in use, and the time that has elapsed since the last arrival. This is true even if the units are all identical, i.e., corresponding to a single resource.
Therefore, computing optimal policies using dynamic programming quickly becomes intractable, and one simple heuristic that has been studied is static policies that post a fixed price regardless of the system state. \cite{besbes2019static} construct a static policy that garners 78.9\% of the reward of the optimal policy in steady state, in reusable resource settings with exponential usage durations. \cite{levi2010provably} also proposed a static policy that is asymptotically optimal when the inventory is large, i.e., the performance loss is sublinear in inventory.

In this paper, we focus on the following research question:

\emph{Can we develop a tractable dynamic policy for pricing reusable resources with general usage durations that offers improved performance compared to static policies?}

To answer this question, we introduce a class of \textit{stock-dependent} policies (see \Cref{def: sd policy}), where decisions are made only based on how many units of stock are available, i.e., current stock level, and are agnostic to how long each unavailable unit of stock has been in use. Therefore, the number of system states required by such a policy grows only linearly with the initial stock. Our policy class shares the spirit of congestion-dependent pricing policies in queueing models \citep{paschalidis2000congestion}---both policies adjust decisions based only on the number of available units.
Note that when the usage duration is memoryless, the optimal dynamic policy lies in the class of stock-dependent policies, which is the benchmark used by \cite{besbes2019static}. \\

\noindent \textbf{Model.\ }
There is one reusable resource with initial stock $c\in \mathbb{N}$ that is used to serve customers. A decision-maker dynamically posts prices to maximize revenue. (Our results also hold for a generic reward function that captures welfare maximization, as we discuss in \Cref{sec: model detailed}.) Customers arrive according to a Poisson process {with rate $\lambda >0$} and decide whether to pay the posted price based on whether their private WTP exceeds it. Once service starts, one of the resource units will be unavailable for a stochastic amount of time, and we denote its expectation by $d$. The WTPs and usage durations are assumed to be i.i.d.~drawn across customers, and also independent of each other. Our results can also be extended to multiple resources and multiple customer classes with different WTPs and usage duration distributions under personalized pricing (see \Cref{sec: extension}).  

\noindent \textbf{Performance Metric.\ }
We consider an infinite time horizon and measure the performance of policies by their long-run average revenue.
This is perhaps the most natural model of reusable resources originally considered in \citet{levi2010provably}, even though there has been much work recently for finite horizons (see \Cref{sec:litRev}).
We benchmark policies by comparing their performance with a fluid relaxation ($\mathsf{FLU}$)---the problem obtained when the inventory constraint only needs to be satisfied in expectation (see \Cref{subsec: flu}, where we show $\mathsf{FLU}$ upper bounds the performance of any online policy). We analyze the \textit{performance loss} of policies, i.e., the additive gap between their long-run average revenue per period and the value of $\mathsf{FLU}$.
We study a regime in which resources are scarce, i.e., {$c/(\lambda d)<1$}, and adopt the large-market asymptotic framework, i.e., scale $c$ and {$\lambda$} proportionally.
This {scaling} regime is motivated by modern large-scale service systems, {where we would like to understand how the performance loss of the service provider changes as both the supply and demand grow.}

\noindent \textbf{Overview of Analysis.\ }
Our analysis is built on deriving the steady-state distribution of the stock level under stock-dependent policies. As a stock-dependent policy prescribes prices, i.e., ex-ante selling probabilities for each stock level, our system can be seen as a loss network with state-dependent arrival rates. Therefore, the steady-state distribution of the stock level under a given policy satisfies a simple set of detailed balance equations akin to those of a birth-and-death process. This "insensitivity result" is striking as the same equations hold for every distribution of usage durations with the same mean. We can further write down an optimization problem to derive the optimal stock-dependent policy, which we show is convex and can be solved efficiently using first-order methods (in the case of continuous WTP distributions) or linear programming (in the case of discrete WTP distributions). These results are presented in \Cref{sec: balance equation}.

\noindent \textbf{performance loss Results.\ }
We first show that static policies cannot achieve a performance loss better than the $\mathcal{O}(\sqrt{c})$ in non-trivial cases.
Our main result is then to provide a sharp characterization of the performance loss of stock-dependent policies via matching upper and lower bounds up to logarithmic factors, which can strictly improve upon the performance loss of static policies.
Specifically, optimal stock-dependent policies achieve a performance loss of $\tilde{\Theta}(c^{1/(1+\alpha)})$, where $\alpha \in [1,\infty]$ is a parameter that depends on the local shape of the revenue curve at $c/(\lambda d)$.
Generally, for continuous WTP distributions we can have $\alpha=2$ and an improved performance loss of $\tilde{\Theta}(c^{1/3})$; 
for discrete WTP distributions we can have $\alpha=\infty$ and a further improved performance loss scaling of ${\mathcal{O}}((\log{c})^2)$ with a lower bound $\Omega(\log{c})$.
These results hold under certain regularity conditions such as differentiability (for continuous WTP) and non-degeneracy (for discrete WTP), which we elaborate on when we define the parameter $\alpha$ in \Cref{sec: main results}.
{We comprehensively characterize performance loss in the different cases (see \Cref{tab: compare policy performance}).}

\noindent\textbf{The Power of Two Prices.\ }
We show that the optimal performance loss of stock-dependent policies can in fact be achieved by a simple two-price policy characterized by 3 parameters instead of $c$ parameters. The two-price policy is a special stock-dependent policy that
sets a high price when the stock level is below some threshold, and
sets a low price otherwise (see \Cref{{def: 2-price policy}}).

The intuition for why two-price policies achieve the optimal performance loss is that they are sufficient to induce a mean reversion process in which the stock level concentrates around the switching threshold and balances between loss due to stock-outs against the gains of setting higher prices (see \Cref{sec: two-price} for further intuition). We note that this "power of two prices" has been previously observed by \cite{kim2018value} in queueing systems with price- and delay-sensitive customers and exponential service times, where they show that a simple policy using only two prices can achieve most of the benefits of dynamic pricing. We establish a related insight for a different, revenue management setting where customers do not enter the queue; also, we are able to explicitly characterize a spectrum of optimal performance loss under different customer WTPs.

\noindent\textbf{Numerical Experiments for Stock-dependent and Two-price Policies.\ }
In \Cref{sec: numerics} we provide some numerical examples to compare the performance of the optimal stock-dependent policy, our two-price policy, a static policy with optimal admission probabilities, and the static policy based on the fluid relaxation. We first test large values of $c$ and display numerical results that complement our theoretical findings with additional insights. We then test their performance under small values of $c$. Results demonstrate that stock-dependent policies provide significant gains over static policies even when $c$ is small. Moreover, we show that a two-price policy performs very closely to the optimal stock-dependent policy.

\noindent\textbf{Organization of the Paper.\ } 
The paper is organized as follows. A formal description of our model is presented in \Cref{sec: model}. The definition of stock-dependent policies and the convex optimization program for deriving the optimal policies within this class are presented in \Cref{sec: balance equation}. The main results on performance guarantees for optimal stock-dependent policies and two-price policies, based on the convex optimization program, are developed in \Cref{sec: main results}. {Extensions to incorporate multiple customer classes and resource types are discussed in \Cref{sec: extension}.} Numerical experiments for both asymptotic and non-asymptotic settings are described in \Cref{sec: numerics}. Concluding remarks and future directions are summarized in \Cref{sec: conclusion}.

\section{Literature Review} \label{sec:litRev}

\noindent\textbf{Reusable Resource Allocation in Infinite-horizon/Steady-state Settings.}
Early works here include \cite{gans2007pricing}, who consider a rental firm facing contract customers and walk-in customers and the firm needs to decide when to serve the former and how much to charge the latter, and \cite{papier2010capacity}, who also investigate capacity rationing policies for rental firms, with part of the customers booking the service in advance. Both works focus on the trade-off between serving two types of customers with different reservation behaviors, while our work focuses more on the trade-off between serving current and future customers. Along these lines, \cite{levi2010provably} and \cite{chen2017revenue} study admission control with general usage durations, and develop near-optimal static policies based on a deterministic LP. \cite{chen2017revenue} assume that usage duration is known upon arrival, whereas our model does not assume this. \cite{besbes2019static} study pricing in reusable resource settings with memoryless usage durations and show that a static pricing policy guarantees 78.9\% times the reward from the optimal policy under regular valuations. Building on this work, \cite{elmachtoub2023power} improve the performance guarantees of static pricing for monotone hazard rate distributions and provide new guarantees for multi-class systems. Our model goes beyond memoryless usage durations and focuses on asymptotic guarantees. In fact, the optimal dynamic policy studied in \cite{besbes2019static} lies in our class of stock-dependent policies. 
\cite{xie2022dynamic} study network revenue management problems with reusable resources where finite types of customers have different value and usage duration distributions. Their policy, akin to our two-price policy, is essentially a trunk reservation policy that rejects customers with lower values if the number of available resource units falls below the trunk reservation level. Such a policy has been demonstrated to be optimal in simpler models such as $M/M/c/c$ queues \citep{key1990optimal, miller1969queueing, morrison2010optimal}, and our results reaffirm the benefits of dynamicity in more complex scenarios involving many customer types and general usage durations. \cite{xie2022dynamic} requires non-degeneracy of the LP, whereas we show that trunk reservation policies have worse performance under degeneracy (in our formulation, degeneracy means that $g$ is non-differentiable at ${c/(\lambda d)}$, such as when ${c/(\lambda d)} = 0.5$ in \Cref{ex: discrete}).

\noindent{\textbf{Reusable Resource Allocation in Finite-horizon Settings.}}
\cite{xu2013dynamic} study dynamic pricing for cloud computing services and characterize some structural results of the optimal dynamic policy. \cite{owen2018price} consider joint pricing and assortment problems, where they demonstrate that decision-makers can take advantage of demand flexibility to balance resource utilization. 
\cite{lei2020real} consider a pricing model for network revenue management with deterministic usage durations, and propose an LP-based resolving heuristic that is asymptotically optimal with large demand and resource capacity. \cite{rusmevichientong_dynamic_2020} consider stochastic arrival sequences and usage durations, where they propose an algorithm based on approximate dynamic programming. \cite{rusmevichientong2023revenue} further investigate how to incorporate reservations into their models. \cite{feng_near-optimal_2020} investigate the same setting and get a near-optimal competitive ratio (for expected reward) with a simulation-based algorithm that samples according to the solution to the expected LP benchmark. \cite{baek_bifurcating_2022} consider network revenue management with general stochastic usage durations and also demonstrate how reusable resources interact when matroid feasibility structures. \cite{doan2020pricing} consider uncertainty in both usage duration and demand; they formulate robust deterministic approximation models to construct efficient fixed-price heuristics. 
\cite{zhang2022online} study a general reusable resource allocation problem that incorporates assignment decisions such as matching and pricing, and multiple objectives such as revenue and market share. They propose a near-optimal policy (for high probability max-min objective) that trades off between rewards earned and resources occupied, where they are allowed to be arbitrarily correlated. By contrast, our work aims to understand how to design simple dynamic policies that perform better than static ones.

The aforementioned works focus on stochastic arrivals, as in our paper. We note that reusable resource allocation has also been extensively studied under adversarial arrivals
\citep{goyal_asymptotically_2021,gong2021online,delong2022online,huo2022online,baek2023leveraging}.

{While these finite‑horizon formulations yield valuable insights, their applicability diminishes when resources recycle continuously and there is not a clear ``start'' and ``end''. We believe the reason for finite horizons is that much of this literature is built upon online resource allocation without reusability, which is naturally a finite-horizon problem.}

\noindent \textbf{Pricing Controls in Queueing Networks.\ } There is a large literature on pricing controls in queueing networks \citep{waserhole2016pricing,balseiro2019dynamic,kanoria_blind_2022}; here we review several works that are most closely related to our paper. \cite{paschalidis2000congestion} study a class of congestion-dependent pricing for network services, which is similar to our stock-dependent policy. They show that optimal dynamic prices increase with congestion and that static pricing is asymptotically optimal. We provide a tight characterization of the performance loss of both static and stock-dependent policies and show that the convergence rate to optimality of static policies is worse.
\cite{banerjee2015pricing} study the value of dynamic pricing in ride-sharing platforms where prices can dynamically react to the imbalance of demand and supply. 
We focus on the centralized allocation of a fixed number of resources with no strategic behavior and consider general usage durations.  
\cite{banerjee2022pricing} study general control policies in a shared vehicle system with fixed numbers of vehicles and stations. They consider state-independent policies under which the resulting Markov chain has a product-form stationary measure. Our formulation is more closely related to \cite{benjaafar2022pricing}, in which the authors focus on static policies and their competitive ratio in networks involving many locations. 
Apart from steady-state analysis, \cite{chen2023real} focus on transient control (pricing and relocation) of a ride-hailing system with non-stationary demand and deterministic travel time between regions. They also consider both static and dynamic policies, providing insights on the value of dynamic pricing. The major difference between our work and these aforementioned papers is that our reusable resources do not travel in space, and our focus is on the inter-dependency of stock level over time.

\section{Problem Formulation}\label{sec: model}
In this section, we formulate a general admission control model for a single reusable resource and describe the performance objective and benchmark.
We describe how it captures the pricing problem from the Introduction as well as other applications.

\subsection{Admission Control Model for a Reusable Resource} \label{sec: model detailed}

We study the admission control problem for a single reusable resource with $c$ units of stock. Customers arrive according to a Poisson process with an arrival rate {$\lambda>0$}. Over continuous time $t$, the decision-maker decides an admission probability $x(t) \in [0,1]$. 
When the admission probability is set to $x$ and there is stock available, customers will enter service following a Poisson process with rate $\lambda x$. The instantaneous reward rate is $\lambda g(x)$, where $g$ is a non-decreasing concave function to be specified later. We assume that when there is no stock available, customers will immediately leave and not enter service. This assumption is motivated by applications such as parking lots and charging stations, where customers are unlikely to wait if no service units are immediately available.
Therefore, without loss of generality and for simplicity of model description, the admission probability must be set to zero when stock-out happens.

The usage duration for each customer is drawn i.i.d. from a distribution with CDF $G$ with expectation $d<\infty$ and $G(0)=0$, during which time one unit of stock would be unavailable. We assume the service is non-preemptive: once the service starts, the decision-maker cannot interrupt the service. Also, while service is ongoing, the decision-maker knows only how long each unit has been in use, but {the realized usage duration is not known until the service finishes.}
The decision-maker needs to design an online policy that controls the admission probability $x(t)$ based on how much stock is available, and how long each unavailable unit has been in use. The objective is to maximize the long-run average reward.

This abstract admission control model captures many general reusable resource allocation problems. Our analysis holds for any reward function $g$ that is concave and non-decreasing with $g(0)=0$, which gives the decision maker the flexibility to choose 1) the objective, 2) the payment collection scheme, such as a one-time fee payment or a rate payment (price per time in use), and 3) the cost structure, such as incorporating operating costs and maintenance fees. For more discussion on the reward function, see \Cref{apx: extension on g}. In particular, our model captures the revenue management problem of dynamic pricing for a single resource. Suppose the distribution of customer WTP $F$ is known, then implementing an admission probability of $x$ corresponds to setting a price such that the probability that the customer has a valuation higher than the price is exactly $x$. In other words, one unit of stock would be consumed w.p. $x$ when there is an arrival, leading to an effective arrival rate of $\lambda x$ by the Poisson splitting theorem. For discrete WTP distributions, the price may have to be randomized to achieve a particular $x$.

In general, for revenue maximization, $g(x)$ is defined to be the increasing concave envelope of the function $x F^{-1}(1-x)$, 
where $F^{-1}(q) = \inf\{v\in \bR: q\leq F(v)\}$ is the inverse CDF for the customer WTP distribution that maps a quantile $q\in [0,1]$ to a valuation for $F$.\footnote{The quantile $q$ is that $q = F(v)$ when the measure corresponding to $F$ has no mass on the realized value of $v$, or uniformly drawn from $(F(v)-\epsilon, F(v)]$ if it has mass $\epsilon>0$ on $v$.}
This increasing concave envelope construction can be implemented randomization or blocking service units.
We provide the following examples to aid the reader and defer further details to \citet[Ch.~3]{hartline2013mechanism}.
\begin{example}\label{ex: discrete}
    If the WTP is either $1$ or $2$ w.p. $1/2$ each, then $x F^{-1}(1-x) = 2x$ when $x \in [0,0.5)$ and $x F^{-1}(1-x) = x$ otherwise. Therefore, by posting a randomized price that is accepted with probability $x$, the revenue obtained is the increasing concave envelope of $x F^{-1}(1-x)$; that is, $g(x) =  \min\{2x, 1\}$.
\end{example}

\begin{example}\label{ex: continuous}
    If the WTP is uniformly distributed in $[0,1]$, then the revenue obtained by posting a price that is accepted with probability $x$ is $x(1-x)$. The reward function is $g(x) = x(1-x)$ for $x\in [0,0.5]$ and $g(x) = 0.25$ for $x\in (0.5, 1]$. The decision-maker posts a price of $1-x$ when $x \leq 0.5$. When $x>0.5$, the decision-maker posts a price of $0.5$, and if an arriving customer cannot afford the price, with probability $2(x-0.5)$, the decision-maker "blocks" a unit by making it unavailable for a random service time drawn from the same distribution as actual service times. This ensures an overall acceptance probability of $0.5+0.5\cdot 2(x-0.5)=x$, while achieving a revenue rate of $g(0.5)$. 
\end{example}

If the payment collection scheme charges customers a price proportional to the time in use, like nightly rates at a hotel, then the formulas for expected revenue should be scaled by the mean usage duration $d$, and in this case, $F$ should be interpreted as the customer WTP distribution for staying one night.
We can similarly handle even more complex payment schemes such as a one-time flat fee plus a pay-per-use charge.

On another note, if the decision-maker aims to maximize welfare instead of revenue, then the formula for the reward function is $g(x)= \int_{1-x}^1 F^{-1}(v)\, dv$.  The reward functions for the WTP distributions from \Cref{ex: discrete,ex: continuous} would then be $g(x) = \min\{2x, x+1/2\}$ and $g(x) = 2x - x^2/2$, respectively.

We note that in the base model, we implicitly assume that the usage duration distribution conditional on the admission probability $x$ is always identical, with mean $d$. This assumption is prevalent throughout the literature \citep{gans2007pricing,papier2010capacity,besbes2019static, rusmevichientong_dynamic_2020, rusmevichientong2023revenue}, and is justified whenever the pricing does not affect the average length of stay of a customer or the firm sells a standardized service that deterministically takes duration $d$ \citep{lei2020real,delong2022online}. However, it's important to note that the assumption of usage duration being independent of WTP can be relaxed. In \Cref{sec: extension}, we extend the model to incorporate multiple customer classes, tackling a {higher level} of complexity. 

\subsection{Performance Metrics}\label{subsec: flu}

Let $\sigma$ denote any feasible adaptive and non-anticipative online policy for our admission control problem, and $x_t^\sigma$ denote the admission probability prescribed by policy $\sigma$ at time $t$ (which is forced to be zero if the system is out of stock).
Motivated by the reusability of resources, we focus on the long-run average reward defined as
\[\liminf_{T\to\infty}\frac{1}{T} \mathbb{E}_{\sigma}\left[\int_{t=0}^T {\lambda} g(x^\sigma_t)\,dt\right],\]
where $\int_{t=0}^T {\lambda} g(x^\sigma_t)\,dt$ is the reward that the policy $\sigma$ collects over the interval $[0, T]$, and the expectation is taken with respect to the probability measure induced by the policy $\sigma$.

To analyze the performance of a policy $\sigma$, we benchmark against the fluid relaxation. This is the problem obtained when the inventory constraint only needs to be satisfied in expectation: the expected number of customers in service cannot exceed the number of total units of stock.
It is analogous to the deterministic LP for
online decision-making problems that have appeared widely in the literature \citep{gallego1994optimal, jasin_re-solving_2012, bumpensanti_re-solving_2018}.
We let $\FLU$ denote the optimal objective value of the fluid relaxation, and show that it serves as an upper bound for the optimal performance of an online policy, in the limit as $T\to\infty$.

Formally, the fluid relaxation is defined as follows:
\begin{align*}
    \mathsf{FLU}=\max_{x \in [0,1]}\ & {\lambda} g(x) \notag \\
    \text{s.t. }& {\lambda} x  d \leq c.
\end{align*}
Note that the constraint is motivated by Little's Law, which says that the expected number of customers in service is ${\lambda} x d $, the effective arrival rate times the expected length of stay, which must be upper-bounded by the total number of available units $c$. As stated earlier, we assume the resource is scarce in that $c/(\lambda d)<1$. In fact, we have {$x^*:= c/(\lambda d)$} in an optimal solution because $g$ is non-decreasing, and hence {$\FLU = \lambda g(x^*) = \lambda g(x^*)$}.
The optimal fluid solution $x^*$ can be interpreted as the expected probability of selling one unit of stock, or equivalently, the proportion of customers that get served in the fluid relaxation.
We will later see that the performance of the policies we are interested in depends solely on the local shape of $g$ around {$x^*$}.

\begin{lemma}\label{lem: FLU benchmark}
For any feasible non-anticipative policy $\sigma$, 
\begin{align*}
\FLU \geq  \liminf_{T\to\infty}\frac{1}{T} \mathbb{E}_{\sigma}\left[\int_{t=0}^T {\lambda} g(x^\sigma_t)\,dt\right].
\end{align*}
\end{lemma}

The proof of \Cref{lem: FLU benchmark} can be found in \Cref{apx: FLU benchmark}.
We refer to the difference between $\FLU$ and $\liminf_{T\to\infty}\frac{1}{T} \mathbb{E}_{\sigma}[\int_{t=0}^T {\lambda} g(x^\sigma_t)\,dt]$ as the \textit{long-run average performance loss} of policy $\sigma$.
For brevity, we will refer to long-run average performance loss as just "performance loss" in the remainder of this paper.

\subsection{Asymptotic Regime}\label{sec: asymptotic regime}
In this paper, we focus on the long-run average reward and, moreover, on an asymptotic setting with a large number of units of stock $c$. This is precisely the setting where dynamic programming is intractable, and the policies we develop will also be shown to empirically perform well for non-asymptotic values of $c$. More specifically, we consider a sequence of increasing problems where the capacity $c$ (supply) and the customer arrival rate $\lambda$ (demand) scales proportionally, with the usage duration distribution and its mean $d$ fixed. Because the fluid benchmark scales with $c$, any policy with a $o(c)$ performance loss would be asymptotically optimal.

\section{Stock-Dependent Policies and Reduction to a Convex Program}\label{sec: balance equation}
In this section, we introduce the class of stock-dependent policies and investigate the steady-state distribution of stock levels induced by such policies. In \Cref{pro: balance equations}, we cite an insensitivity result that allows us to deal with general usage durations analytically. We then formulate an optimization problem to derive the optimal stock-dependent policy and show that the formulation is convex. Finally, we discuss several important properties of this optimization problem.

\begin{definition}\label{def: sd policy}
A \textit{stock-dependent} policy (for $c$ units of initial stock) is defined by parameters $x_1,\ldots,x_c$, where $x_j\in[0,1]$ prescribes the admission probability when $j$ units of stock are available, regardless of how long the unavailable units have been in use.
\end{definition}

Let $S_t$ denote the number of units of stock available at time $t$, then we can define the limiting average probability that exactly $j$ units of stock are available as $\pi_j = \lim_{T\to\infty} \frac{1}{T} \int_0^T \mathbb{P}\{S_t = j \}\,dt$. Then we have the following characterization of the steady state using existing results for the loss network.

\begin{proposition}[\cite{brumelle1978generalization} Theorem 3]\label{pro: balance equations}
If the decision-maker adopts admission probabilities $\bm{x}:=\{x_j\}_{j=1}^c$, then the steady-state probabilities $\bm{\pi}:=\{\pi_j\}_{j=0}^c$ of our system is the unique solution to the equations:
\begin{align}
    \pi_j {\lambda} x_j= \pi_{j-1}\frac{c-j+1}{d}, \forall j \in [c], \label{balance equations}
\end{align}
that also satisfies $\sum_{j=0}^c \pi_j = 1$. 
\end{proposition}

\Cref{pro: balance equations} is an insensitivity result---the steady-state distribution of the number of customers in the system is insensitive to the form of the usage duration distribution, and requires only the mean usage durations. The insensitivity property of the Erlang loss model with Poisson input is well-established \citep{tijms2003first}, and the intuition is that arriving customers never queue. \Cref{pro: balance equations} extends this result to state-dependent arrival rates.

To get some intuition for \Cref{pro: balance equations}, suppose the usage duration is exponentially distributed with mean $d$. Then,~\eqref{balance equations} would be the detailed balance equations. More specifically, when there are $j$ units of stock available, the transition rate to one unit less is $\pi_j{\lambda}x_j$ as the admission probability is ${\lambda}x_j$. On the other hand, when there are $j-1$ units of stock available, {there are $c-(j-1)$ units in service, and hence} the transition rate to one unit more is $\pi_{j-1}\cdot(c-j+1)/d$, as each unit has a service rate of $1/d$ implying a total service rate of $(c-j+1)/d$.

Since $\bm{\pi}$ is uniquely defined from $\bm{x}$, we can formulate the optimization problem over stock-dependent policies as follows:
\begin{align}
    \max_{\bm{\pi}\geq \bm{0}, \bm{x} \geq \bm{0}} & \sum_{j=1}^c\pi_j {\lambda} g(x_j) \label{OPTSD} \\
   \text{s.t. }  &\pi_j{\lambda}x_j = \pi_{j-1}\frac{c-j+1}{d}, \forall j \in [c], \notag \\
   & x_j \leq 1, \forall j \in [c], \notag \\
   & \sum_{j=0}^c\pi_j = 1. \notag
\end{align}

Note that the admission probabilities $x_j$ can be set to arbitrary values in [0,1], and the objective function multiplies each reward rate ${\lambda} g(x_j)$ by the steady-state probability of having exactly $j$ units in stock.

We now introduce some properties of the stock-dependent policy optimization problem stated in \eqref{OPTSD}. 

\begin{proposition}\label{pro: optsd convex}
The optimization problem \eqref{OPTSD} can be reformulated as a convex program in $\bm{\pi}$.
\end{proposition}
The proof of \Cref{pro: optsd convex} (see \Cref{apx: optsd convex}) relies on eliminating the control variables $\bm{x}$ to restate the problem in terms of $\bm{\pi}$. Then, the objective function can be written as a perspective function of $g$, which preserves convexity in $\bm{\pi}$. For piecewise linear $g$, e.g., when there are finite customer types,
we can further simplify \eqref{OPTSD} to an LP.

\begin{theorem}\label{thm: optsd to lpsd}
Suppose there exists $k\in [m]$, $a_k \geq 0$, and $b_k\in \bR$ such that $g(x)=\min_{k\in [m]}\{a_k + b_k x\}$, i.e., $g$ is the minimum of non-decreasing linear functions. Then, \eqref{OPTSD} can be reformulated into the following LP:
    \begin{align}
    \max_{\bm{\pi}\geq \bm{0}} & \sum_{j=1}^c {\lambda} g_j \label{LPSD} \\
   \text{s.t. }  & g_j \leq a_k\pi_j+b_k\frac{c-j+1}{{\lambda}d}\pi_{j-1}, \forall k \in [m], \forall j \in [c], \notag \\
   & \pi_j \geq \frac{c-j+1}{\lambda  d}\pi_{j-1}, \forall j \in [c], \notag \\
   & \sum_{j=0}^c\pi_j = 1. \notag
   \end{align}
\end{theorem}
The proof of \Cref{thm: optsd to lpsd} (see \Cref{apx: optsd to lpsd}) is based on showing the equivalence between $g(x_j) \leq a +b x_j$ and $\pi_jg(x_j) \leq a \pi_{j} + b \pi_{j-1} (c-j+1)/({\lambda}d)$, where $\pi_jx_j$ is replaced by $\pi_{j-1}(c-j+1)/({\lambda}d)$ following \Cref{pro: balance equations}.  Letting $g_j := \pi_jg(x_j)$ be our new decision variables then completes the proof. \Cref{thm: optsd to lpsd} plays a key role in our theoretical analysis as it allows us to use linear programming duality to provide upper bounds on the optimal reward rate of stock-dependent policies after discretizing the reward function $g$. 

\section{Main Results} \label{sec: main results}

In \Cref{sec: balance equation} we reduced the computation of the long-run average reward of any stock-dependent policy to a convex program.
In this section, we focus on analyzing the objective value of this program and its performance loss compared to the fluid benchmark, i.e., $\FLU = {\lambda}g(x^*)$, as the problem scale grows large.

\subsection{Performance Loss of Static Policies}\label{subsec: static policies}

We start with analyzing the performance loss of the fluid policy, which sets a static admission probability of $x^*$ whenever there is stock available. This policy, also known as the certainty-equivalent heuristic, is widely used in various decision-making problems under uncertainty, such as non-reusable resource allocation and queueing network controls.

\begin{proposition}\label{pro: failure of fluid pricing}
The performance loss of the fluid policy is $\Theta(\sqrt{c})$.
\end{proposition}

Specifically, for any optimal fluid solution value $x^*\in(0,1)$, i.e., for any parameter settings of $\lambda, c$ and $d$ such that $c/(\lambda d) \in (0,1)$, and for any reward function $g$ that is concave and non-decreasing with $g(0)=0$, the performance loss of the fluid policy is $\Theta(\sqrt{c})$.
We analyze the stock-out probability to establish performance bounds (see \Cref{apx: failure of fluid pricing}). Such an approach was also used in \citet{levi2010provably} to provide the same upper bound on stock-out probability, and here we complement their result by establishing a matching lower bound.

Note that for a finite $c$, the fluid policy is not necessarily optimal among all static policies that set a fixed admission probability whenever there is stock available. Specifically, the fluid policy can be trivially improved when: 1) it is rejecting the highest valuation; 2) it is accepting customers with zero valuations. 
The following two special cases occur when the reward function is locally linear around $x^*$.
Let $\partial g(x)$ denote the set of supergradients of $g$ at $x\in [0,1]$. The proofs for the following two propositions can be found in \Cref{apx: 2special cases}. 

\begin{proposition}[Extremely Scarce Supply Helps]\label{pro: scarce supply helps}
If $\partial g(x) = \{r\}$ for all $x\in [0,x^*+\varepsilon]$ for some constants $r>0$ and $\varepsilon>0$, the performance loss of the optimal static policy is $\mathcal{O}(1)$. 
\end{proposition}

Recall that $x^*$ = $c/(\lambda d)$ is the load factor that quantifies the expected proportion of customers the decision-maker can (at most) serve in the fluid relaxation. The assumption that $\partial g(x) = \{r\}$ for all $x\in [0,x^*+\varepsilon]$ means that reward increases linearly in admission probability for $x\in [0, x^*+\varepsilon]$.
This implies that the proportion of the highest-valuation customers exceeds maximum feasible average admission rate $x^*$.
Therefore, the optimal static policy admits all customers from this group and rejects all others. In particular, the performance loss of the optimal static policy is $o(\sqrt{c})$ because the decision-maker can reliably fill capacity using only the high-valuation group and avoid any trade-off or risk of allocating capacity to lower types.

\begin{proposition}[Excess Supply Helps]\label{pro: excess supply helps}
If $\partial g(x^*) = \{0\}$, the performance loss of the optimal static policy is $\mathcal{O}(e^{-c})$. 
\end{proposition}

The assumption that $\partial g(x^*) = \{0\}$ means that the reward function is locally flat at $x^*$. This implies that, when posting an admission probability $x^*$, some admitted customers generate zero reward. Thus, with excess supply, it is optimal to set a lower admission probability that excludes these customers, since this would reduce stock-outs without hurting the reward. The performance loss of the static policy is $o(\sqrt{c})$ as setting a lower admission probability yields an exponentially small probability of stocking out.

Beyond the above two special cases, any static policy would still incur $\Omega(\sqrt{c})$ performance loss, as formalized by the following proposition. The proof can be found in \Cref{apx: lower bound statis pricing}.
\begin{proposition}\label{pro: failure of static pricing}
Excluding the cases stated in \Cref{pro: scarce supply helps} and \Cref{pro: excess supply helps}, the performance loss of any static policy is $\Omega(\sqrt{c})$.
\end{proposition}

In the remainder of the paper, we provide conditions under which a stock-dependent policy, {in particular a two-price policy}, achieves $o(\sqrt{c})$ performance loss. \Cref{tab: compare policy performance} delineates the cases analyzed here and previews results in the following \Cref{sec: upper bound,sec: lower bound} (the definition of $\alpha$ can be found in \Cref{assmp: condition for upper bound} and \Cref{assmp: condition for lower bound}).

\begin{table}[h]
\centering
\renewcommand{\arraystretch}{1.5}
\begin{tabular}{|>{\centering\arraybackslash}m{4cm}|>{\centering\arraybackslash}m{1.5cm}|>{\centering\arraybackslash}m{2cm}|>{\centering\arraybackslash}m{2.5cm}|>{\centering\arraybackslash}m{3cm}|>{\centering\arraybackslash}m{2cm}|}
\hline
\multirow{2}{*}{} &  \multirow{2}{*}{$\alpha = 1$} & \multirow{2}{*}{$\alpha \in (1, \infty)$} & \multicolumn{3}{c|}{$\alpha = \infty$} \\
\cline{4-6}
& & & $\partial g(x^*) = \{0\}$ & $\partial g(x) = \{r\}, \forall x \in [0, x^* + \varepsilon], r> 0, \varepsilon > 0$ & Otherwise \\
\hline
Fluid policy & \multicolumn{5}{c|}{$\Theta(\sqrt{c})$ Prop \ref{pro: failure of fluid pricing} } \\
\hline
Optimal static policy & 
\multirow{2}{*}{\parbox{1.5cm}{\centering\vspace{0.5em} $\Theta(\sqrt{c})$ \\ Prop \ref{pro: failure of fluid pricing} \\  Prop \ref{pro: failure of static pricing}\vspace{0.5em}}} 
& $\Theta(\sqrt{c})$  Prop \ref{pro: failure of fluid pricing} Prop \ref{pro: failure of static pricing} & 
\multirow{2}{*}{\parbox{2.5cm}{\centering \vspace{0.5em}$\mathcal{O}(e^{-c})$ \\ Prop \ref{pro: excess supply helps}}} 
& \multirow{2}{*}{\parbox{3cm}{\centering \vspace{0.5em} $\mathcal{O}(1)$ \\ Prop \ref{pro: scarce supply helps}}} 
& $\Theta(\sqrt{c})$ Prop \ref{pro: failure of fluid pricing} Prop \ref{pro: failure of static pricing} \\
\cline{1-1} \cline{3-3} \cline{6-6}
Two-price policy &  & $\tilde{\mathcal{O}}(c^{1/(1+\alpha)})$ Thm \ref{thm: upper bound}& & &${\mathcal{O}}((\log{c})^2)$ Thm \ref{thm: upper bound} \\
\hline
Stock-dependent policy {(lower bounds)} & \multicolumn{2}{c|}{${\Omega}(c^{1/(1+\alpha)})$ Thm \ref{thm: lower bound}} & $\Omega(e^{-c})$ Thm \ref{thm: lower bound} & {$\Omega(1)$} Thm \ref{thm: lower bound} & $\Omega(\log{c})$ Thm \ref{thm: lower bound} \\
\hline
\end{tabular}
\caption{Performance loss of different policies under different assumptions on the reward function $g$. {All results are new, except the $O(\sqrt{c})$ upper bound for the fluid policy.} The parameter $\alpha$ captures the local curvature of $g$ around the optimal fluid admission probability $x^*$. Here $\partial g(x)$ denote the set of supergradients of $g$ at $x\in [0,1]$.}
\label{tab: compare policy performance}
\end{table}

\subsection{Performance Loss Upper Bound for Two-Price Stock-dependent Policies}\label{sec: upper bound}

We proceed by deriving an upper bound on the performance loss of a simple two-price stock-dependent policy, which also provides the upper bound on the performance loss of the optimal stock-dependent policy. We outline our result in this subsection and provide the formal proof in the next subsection.

\begin{definition}\label{def: 2-price policy}
A \textit{two-price} (stock-dependent) policy is parameterized by $x_L \in [0,1], x_H \in [0,1], \tau \in [c]$, where $x_j=x_L$ for $1\leq j \leq \tau$ and $x_j = x_H$ otherwise.
\end{definition}
That is, a two-price policy dynamically chooses between two admission probabilities $x_L$ (low) and $x_H$ (high) depending on whether the stock level is above a threshold $\tau$. This is not to be confused with randomization between two prices to achieve a specific admission probability, which even the fluid policy may perform when the WTP distribution is discrete.

The shape of the reward function $g$ near the fluid admission probability of $x^*$ is the key driver of performance. The following assumption captures the notion of the shape of interest; the same assumption with $\alpha=2$ has been previously used by \cite{balseiro_survey_2021} to derive performance loss results for other pricing problems. Let $\partial g(x)$ be the set of supergradients of $g$ at $x\in [0,1]$, and an element of $\partial g(x)$ is denoted by $g'(x)$.

\begin{assumption}\label{assmp: condition for upper bound}
There exists $0< \varepsilon < \min\{x^*, 1-x^*\} , \alpha\in[1,\infty], k_1 \ge 0$, and a supergradient $g'(x^*) \geq 0$ such that $g(x)\geq g(x^*)+g'(x^*)(x-x^*)-k_1|x-x^*|^\alpha$ for all $x\in[x^*-\varepsilon,x^*+\varepsilon]$.
\end{assumption}

\begin{figure}[!htbp]
\captionsetup{justification=centering}
\begin{subfigure}[b]{0.3\textwidth}
    \centering
    \begin{tikzpicture}
      \begin{axis}[
        axis lines = middle,
        xlabel = $x$,
        ylabel = $g(x)$,
        ymin=0, ymax=1,
        xmin=0, xmax = 1.1,
        xtick={0.001, 0.5, 1},
        xticklabels={0,$x^*$,1},
        ytick=\empty,
        width=5cm,
        height=4cm,
        every axis x label/.style={at={(current axis.right of origin)},anchor=west},
        every axis y label/.style={at={(current axis.north west)},above=2mm},   
        ]
        \addplot[domain=0:0.5, samples=100, blue, thick] {x};
        \addplot[domain=0.5:1, samples=100, blue, thick] {0.5*x+0.25};
        \node[below] at (axis cs:0.8,0.3) {${\mathcal{O}}(\sqrt{c})$};
        \addplot[domain=0.25:0.5, samples=100, dotted, red, thick] {1.4*(x-0.5)+0.5};
        \addplot[domain=0.5:0.75, samples=100, dotted, red, thick] {0.1*(x-0.5)+0.5};
        \draw[densely dotted] (axis cs:0.5,0) -- (axis cs:0.5,0.5);
    \end{axis}
\end{tikzpicture}
\caption{$\alpha = 1$}
\end{subfigure}
    \hspace{0.03\textwidth}
\begin{subfigure}[b]{0.3\textwidth}
    \centering
    \begin{tikzpicture}
    \begin{axis}[
        axis lines = middle,
        xlabel = $x$,
        ylabel = $g(x)$,
        ymin=0, ymax=1,
        xmin=0, xmax = 1.1,
        xtick={0.001, 0.5, 1},
        xticklabels={0,$x^*$,1},
        ytick=\empty,
        width=5cm,
        height=4cm,
        every axis x label/.style={at={(current axis.right of origin)},anchor=west},
        every axis y label/.style={at={(current axis.north west)},above=2mm},
        ]
        \addplot[domain=0:0.6, samples=100, blue, thick] {2*x*(1.2-x)};
        \addplot[domain=0.6:1, samples=100, blue, thick] {0.72};
        \node[below] at (axis cs:0.8,0.3) {$\tilde{\mathcal{O}}(c^{1/3})$};
        \addplot[domain=0.25:0.75, samples=100, dotted, red, thick] {0.7 + 0.4*(x-0.5) - 8*(x-0.5)^2};
        \draw[densely dotted] (axis cs:0.5,0) -- (axis cs:0.5,0.7);
    \end{axis}
\end{tikzpicture}
\caption{$\alpha = 2$}
\end{subfigure}
    \hspace{0.03\textwidth}
\begin{subfigure}[b]{0.3\textwidth}
\begin{tikzpicture}
\centering
    \begin{axis}[
        axis lines = middle,
        xlabel = $x$,
        ylabel = $g(x)$,
        ymin=0, ymax=1,
        xmin=0, xmax = 1.1,
        xtick={0.001, 0.5, 1},
        xticklabels={0,$x^*$,1},
        ytick=\empty,
        width=5cm,
        height=4cm,
        every axis x label/.style={at={(current axis.right of origin)},anchor=west},
        every axis y label/.style={at={(current axis.north west)},above=2mm},   
        ]
        \addplot[domain=0:0.3, samples=100, blue, thick] {x};
        \addplot[domain=0.3:1, samples=100, blue, thick] {0.5*x + 0.15};
        \node[below] at (axis cs:0.8,0.3) {${\tilde{\mathcal{O}}}(1)$};
        \addplot[domain=0.15:0.75, samples=100, dotted, red, thick] {0.5*x+0.15};
        \draw[densely dotted] (axis cs:0.5,0) -- (axis cs:0.5,0.4);
    \end{axis}
\end{tikzpicture}
\caption{$\alpha = \infty$}
\end{subfigure}
\caption{Illustration of \Cref{assmp: condition for upper bound}. The blue lines depict $g(x)$ while the red dotted lines depict $g(x^*)+g'(x^*) (c-x^*)- k_1|x-x^*|^\alpha, x\in [x^*-\varepsilon, x^*+\varepsilon]$, in different cases.} 
\label{fig: illus of assmp ub}
\end{figure}
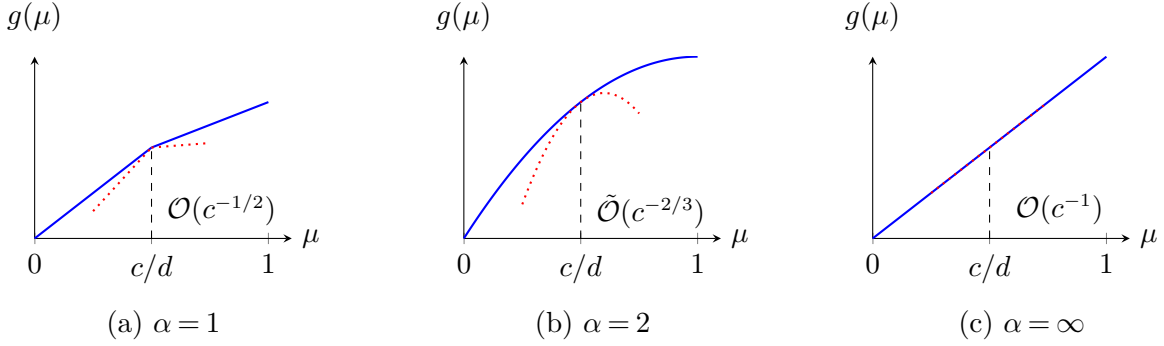

\Cref{assmp: condition for upper bound} can be interpreted as follows. {We first note that because $|x-x^*|\le\eps<1$, the term $|x-x^*|^\alpha$ can be understood to be 0 when $\alpha=\infty$.} The term $g(x^*)+g'(x^*)(x-x^*)$ is a local approximation of $g(x)$ around $x^*$, which is guaranteed to overestimate $g$ since $g$ is concave.
\Cref{assmp: condition for upper bound} is imposing that for a sufficiently small ball of radius $\varepsilon$ around $x^*$, we can upper-bound this difference between $g(x)$ and its approximation using $k_1|x-x^*|^\alpha$:
\begin{align}\label{eqn: condition for upper bound}
g(x^*)+g'(x^*)(x-x^*) - g(x) \le k_1|x-x^*|^\alpha.
\end{align}
The larger the value of the shape parameter $\alpha$ can be made to satisfy \Cref{assmp: condition for upper bound} around a small ball, the faster that the stock-dependent policy will converge to optimality. Here are some prototypical cases and the associated values of $\alpha$ (see \Cref{fig: illus of assmp ub} below for illustration):
\begin{enumerate}
\item If $g$ is linear around $x^*$, then~\eqref{eqn: condition for upper bound} can be satisfied using $\alpha=\infty$ and $k_1 = 0$;
\item As long as $g$ is twice-differentiable at $x^*$, by Taylor's theorem,~\eqref{eqn: condition for upper bound} can be satisfied using $\alpha=2$;
\item Finally, if $g$ is piecewise-linear and non-differentiable at $x^*$, then $\alpha$ must be as small as 1 in order to satisfy~\eqref{eqn: condition for upper bound}.
\end{enumerate}
For more discussions, including examples of intermediate cases for $\alpha$, see \Cref{apx: two assumptions}.

\begin{theorem}\label{thm: upper bound}
Under \Cref{assmp: condition for upper bound}, if the supergradient {further satisfies} $g'(x^*) >0$, then there exists a two-price policy with performance loss $\Tilde{\mathcal{O}}(c^{1/(\alpha+1)})$. Specifically, when $\alpha = \infty$, the performance loss is $\mathcal{O}((\log{c})^2)$.
\end{theorem}

When there does not exist a supergradient $g'(x^*)>0$, i.e.\ we are in the excess supply case discussed in \Cref{pro: excess supply helps}, {there is an even better guarantee---}the optimal static policy achieves a performance loss of $\mathcal{O}(e^{-c})$. As a two-price policy can be reduced to a static policy, this performance guarantee applies. Similarly, when $\alpha = \infty$ and there exists $r >0, \varepsilon >0$ such that $\partial g(x) = \{r\}$ for all $x \in [0,x^* + \varepsilon]$, i.e., the extremely scarce supply case discussed in \Cref{pro: scarce supply helps}, the performance loss of the optimal static policy, $\mathcal{O}(1)$, also applies. The results are summarized in \Cref{tab: compare policy performance}. 

\textbf{Implications of \Cref{thm: upper bound}}.
We return to the prototypical cases for $\alpha$, and provide some concrete examples alongside the performance loss result implied by \Cref{thm: upper bound}.

\begin{enumerate}
\item The reward function $g$ is locally linear at $x^*$: $\alpha = \infty$, performance loss is $\mathcal{O}((\log{c})^2)$. This case arises if the WTP is discrete and there is no valuation at which $F$ equals $x^*$. For example, in \Cref{ex: discrete} with $x^*\neq 0.5$.
\item The reward function $g$ is locally smooth at $x^*$: $\alpha=2$, performance loss is $\tilde{\mathcal{O}}(c^{1/3})$.
{This case can arise if the WTP distribution is continuous; e.g.\ in \Cref{ex: continuous} with $x^*\in (0,0.5)$.
It arises for all}
continuous WTP distributions when the objective is welfare maximization, {and under mild conditions when the objective is revenue maximization.
We note that continuous distributions commonly used for pricing \citep[e.g.][Ch.~7.3.3]{talluri2004theory} would typically fall under this case. We provide further discussion of the conditions in \Cref{apx: two assumptions}.}
\item The reward function $g$ is non-differentiable at $x^*$: $\alpha = 1$, performance loss is ${\mathcal{O}}(\sqrt{c})$. This case can arise if the WTP is discrete and $F$ equals $x^*$ at some valuation.
For example, in \Cref{ex: discrete} with $x^*= 0.5$.
\end{enumerate}

When $\alpha>1$, stock-dependent policies improve upon the performance loss of $\sqrt{c}$, which is the best performance loss of \textit{any} policy for reusable resource allocation known in the literature.
If $g$ is non-differentiable at the fluid admission probability $x^*$, i.e., $\alpha = 1$, then we show in the subsequent \namecref{sec: lower bound} that the performance loss of stock-dependent policies in fact cannot {be improved beyond} $\sqrt{c}$. 

\subsection{Proof of \Cref{thm: upper bound}}\label{apx: upper bound}

When $\alpha = 1$, a direct application of \Cref{pro: failure of fluid pricing} implies a performance loss of $\mathcal{O}(\sqrt{c})$. 
Our proof relies on approximating $g$ from below using a function $\hat{g}$ and analyzing the performance of a {two-price} policy on $\hat{g}$. Specifically, we construct $\hat{g}$ so that is linear within $[x^*-\delta(c), x^*+\delta(c)]$, for some $\delta(c)>0$ to be determined. The two-price policy adopts a low admission probability $x_L = x^*-\delta(c)$, a high admission probability $x_H=x^*+\delta(c)$, and a switching threshold $\tau(c)$. When $\alpha = \infty$, $\delta(c)$ is fixed. When $\alpha\in (1,\infty)$, we design $\delta(c)$ to decrease with $c$, i.e., the difference between the two admission probabilities is shrinking as the problem scale increases. Optimizing $\delta(c)$ and tuning the relationship between $\tau(c)$ and $\delta(c)$ suffices to obtain the {performance loss of $\Tilde{\mathcal{O}}(c^{1/(\alpha+1)})$}. 

We prove the result in four steps. First, we characterize the steady-state distribution of the stock level under a general two-price policy characterized by $x_L < x_H \in (0,1)$ with $x^*-x_L = x_H - x^*$ and switching threshold $\tau$. Second, we tune $x_L, x_H$, and $\tau$ and present a sharper bound on the steady-state distribution. Third, we explicitly define the approximation $\hat g$ and provide a bound on the performance loss of the two-price policy under $\hat g$ based on the steady-state probabilities. In the last step, we connect the performance loss under $\hat g$ with that of the original problem, i.e., $g$, by leveraging \Cref{assmp: condition for upper bound}.

\subsubsection*{Step 1. General Bound on Steady-State Probabilities of Two-Price Policies.}\label{sec: bound probabilities}

Consider a two-price policy parameterized by $\tau,x_L ,x_H$. Let $\pi_0$ denote the steady-state probability that the system runs out of stock, and let $\pi_L = \sum_{j=1}^\tau \pi_j$ and $\pi_H = \sum_{j=\tau+1}^c \pi_j$ denote the steady-state probability that the admission probabilities are low and high, respectively. In this step, we characterized the bounds on the steady-state probabilities. For its proof, see \Cref{apx: bound varsigma}

\begin{lemma}\label{lem: bound varsigma}
For $\tau \leq c/2, 0 < x_L < x_H < 1$ and $x^*-x_L = x_H-x^*$, let $\varsigma_L = {\pi_L}/{\pi_0}$ and $\varsigma_H = {\pi_H}/{\pi_0}$, then the following holds:
\begin{align}
     &\varsigma_L \geq  \left(1-\frac{\tau^2}{2c}\right)  \frac{(x^*/x_L)^\tau - 1}{1-x_L/x^*}\label{eq: lb varsigmaL}, \\
     & \varsigma_L - \varsigma_H \leq \frac{(x^*/x_L)^\tau}{1-x_L/x^*} \left( \frac{2\tau^2}{c} + (x^*/x_H)^\tau \right) \label{eq: up diff vL vH}.
\end{align}
\end{lemma}

\subsubsection*{Step 2. Bound on Steady-State Probabilities under Specific $\tau$ and $\delta(c)$.}

In this step, we consider a two-price policy that sets
\begin{align}\label{eq: two price para}
    \tau(c) = \left \lceil\frac{\log{c}}{\log(1+\delta(c)/x^*)} \right \rceil, \quad \delta(c) = c^{-1/(\alpha+1)} \text{ when } \alpha \in (1,\infty), \quad \delta(c) =  \varepsilon \text{ when }\alpha = \infty. 
\end{align}

We first present two useful lemmas that directly follow from the choice of $\tau(c)$. For the proofs, see \Cref{apx: property of tau} and \Cref{apx: bound tau2c}.
\begin{lemma}\label{lem: property of tau}
Under the two-price policy with parameters specified as \eqref{eq: two price para}, for $c \geq 4$, the following holds
\begin{align}
    &c \leq (1+\delta(c)/x^*)^{\tau(c)}  \leq (1+1/x^*)c \label{eq: equiva tau c},\\
    &\frac{\log{c}}{\delta(c)/x^*} \leq \tau(c) \leq \frac{3\log{c}}{\delta(c)/x^*}, \label{eq: bound tau}\\
    &(1-\delta(c)/x^*)^{-\tau(c)} \geq (1+\delta(c)/x^*)^{\tau(c)} \geq c. \label{eq: minus tau bigger}
\end{align}
\end{lemma}

\begin{lemma}\label{lem: bound tau2c}
For $\alpha \in (1,\infty]$, and any constant $\Gamma > 0$, there exists a $N(\Gamma)$ such that for any $c > N(\Gamma)$, 
    \[\frac{\tau(c)^2}{c} \leq \Gamma.\]
\end{lemma}

Let $c\geq \max\{N(1), 4\}$. By \Cref{lem: bound tau2c}, for $c\geq N(1)$, $\tau(c)^2/c \leq 1$, thus $\tau(c) \leq \sqrt{c}$. As $c\geq 4$, $\sqrt{c} \leq c/2$, then we have $\tau(c) \leq c/2$. Therefore, we can apply \Cref{lem: bound varsigma} inequality \eqref{eq: lb varsigmaL}:

\begin{align}
\varsigma_L \stackrel{(a)}{\geq} \left(1-\frac{\tau(c)^2}{2c}\right) \frac{(1-\delta(c)/x^*)^{-\tau(c)}-1}{\delta(c)/x^*} \stackrel{(b)}{\geq} \frac{1}{2} \frac{(1-\delta(c)/x^*)^{-\tau(c)}-1}{\delta(c)/x^*} \stackrel{(c)}{\geq} \frac{(1-\delta(c)/x^*)^{-\tau(c)}}{4\delta(c)/x^*}, \label{eq: lb varsigmaL sharper}
\end{align}
where $(a)$ follows from replacing $x_L/x^* = (x^*-\delta(c))/x^* = 1-\delta(c)/x^*$; $(b)$ follows $\tau(c)^2/c \leq 1$, thus $(1-\tau(c)^2/c) \geq 1/2$; $(c)$ follows as $y-1 \geq y/2$ for all $y \geq 2$, and by \eqref{eq: minus tau bigger}, we have $(1-\delta(c)/x^*)^{-\tau(c)}\geq c \geq 4 \geq 2$.

Again by \Cref{lem: bound varsigma}, inequality \eqref{eq: up diff vL vH}: 
\begin{align}
    \varsigma_L-\varsigma_H & \leq \frac{(x^*/x_L)^{\tau(c)}}{1-x_L/x^*} \left( \frac{2\tau(c)^2}{c} + (x^*/x_H)^{\tau(c)} \right) \notag \\
    & \stackrel{(d)}{\leq} \frac{ (1-\delta(c)/x^*)^{-\tau(c)}}{\delta(c)/x^*} \left( \frac{2\tau(c)^2}{c} + (1+\delta(c)/x^*)^{-\tau(c)} \right)  \notag\\
    & \stackrel{(e)}{\leq} \frac{ (1-\delta(c)/x^*)^{-\tau(c)} }{\delta(c)/x^*}\left( \frac{2\tau(c)^2}{c} + \frac{1}{c} \right)  \notag \\
    & \stackrel{(f)}{\leq} \frac{ (1-\delta(c)/x^*)^{-\tau(c)} }{\delta(c)/x^*}\frac{3\tau(c)^2}{c},
    \label{eq: ub diff varsigma}
\end{align}
where $(d)$ holds by replacing $x_L/x^* = (x^*-\delta(c))/x^* = 1-\delta(c)/x^*$ and $x_H/x^* = (x^*+\delta(c))/x^* = 1+\delta(c)/x^*$; $(e)$ holds by \eqref{eq: equiva tau c}; $(f)$ holds as $1\leq \tau(c)$.

\subsubsection*{Step 3. Performance Loss of the Two-Price Policy \eqref{eq: two price para} under $\hat{g}$.} 
Consider an approximation of $g$, denoted by $\hat{g}$, which interpolates linearly between $(x^*-\delta(c), x^*+\delta(c))$ and is equal to $g(x)$ otherwise. The function is defined as follows (see \Cref{fig: illus of hat{g}} for an illustration):
$$
\hat{g}(x) = \left\{
\begin{aligned}
&g(x),& 0\leq x\leq x^*-\delta(c),\\
& g(x^*-\delta(c)) + \frac{g(x^*+\delta(c))-g(x^*-\delta(c))}{2\delta(c)}(x-(x^*-\delta(c))), & x^*-\delta(c) < x < x^*+\delta(c),\\
&g(x), & x^*+\delta(c) \leq x \leq 1.
\end{aligned}
\right.
$$

Denote $\hat{g}'(x^*) = ({g(x^*+\delta(c))-g(x^*-\delta(c))})/({2\delta(c)})$. By \Cref{assmp: condition for upper bound}, there exists a supergradient $g'(x^*) >0$, thus $g(x^*+\delta(c))-g(x^*-\delta(c)) >0$, and $\hat{g}'(x^*) >0$.

We proceed to bound the performance loss of the two-price policy specified above as \eqref{eq: two price para} under the reward function $\hat{g}$.

The reward of the two-price policy under $\hat{g}$, denoted as $\hat{\mathcal{R}}$, is
\begin{align}
    \hat{\mathcal{R}} &  =\pi_L \lambda \hat{g}(x_L)+\pi_H \lambda \hat{g}(x_H) \notag \\
    & \stackrel{(a)}{=}\pi_L \lambda (\hat{g}(x^*)-\hat{g}'(x^*)\delta(c)) + \pi_H \lambda (\hat{g}(x^*)+\hat{g}'(x^*)\delta(c))\notag\\
    & \stackrel{(b)}{=} (1-\pi_0) \lambda \hat{g}(x^*) + (\pi_H-\pi_L)\lambda \hat{g}'(x^*)\delta(c), \label{eq: hat{R}}
\end{align}
where $(a)$ holds as $\hat{g}$ is locally linear in $[x^*-\delta(c),x^*+\delta(c)]$; $(b)$ holds as $\pi_0+\pi_L+\pi_H = 1$. 

Then, the performance loss of the two-price policy under $\hat g$ is
\begin{align}\label{eq:regret-hat1}
    \lambda \hat{g}(x^*) - \hat{\mathcal{R}} & \stackrel{(c)}{=} \pi_0 \lambda \hat{g}(x^*)+(\pi_L-\pi_H)\lambda \hat{g}'(x^*)\delta(c) \nonumber\\
    & \stackrel{(d)}{=} \frac{\lambda\hat{g}(x^*)+(\varsigma_L-\varsigma_H)\lambda\hat{g}'(x^*)\delta(c)}{1+\varsigma_L+\varsigma_H}\nonumber\\
    & \stackrel{(e)}{\leq} \frac{\lambda\hat{g}(x^*)+(\varsigma_L-\varsigma_H)\lambda\hat{g}'(x^*)\delta(c)}{\varsigma_L},
\end{align}
where $(c)$ follows from \eqref{eq: hat{R}}; $(d)$ holds as $\pi_0 = 1/(1+\varsigma_L + \varsigma_H)$, $\pi_L = \pi_0 \varsigma_L$, and $\pi_H = \pi_0 \varsigma_H$; $(e)$ holds as $1+\varsigma_H \geq 0$.

We next plug in the bounds of $\varsigma_L, \varsigma_H$. For $c \geq \max\{N(1), 4\}$, by \eqref{eq: lb varsigmaL sharper} and \eqref{eq: ub diff varsigma},

\begin{align}\label{eq:regret-hat2}
     \lambda \hat{g}(x^*)-\hat{\mathcal{R}} & \leq \frac{\lambda\hat{g}(x^*)+ \frac{ (1-\delta(c)/x^*)^{-\tau(c)} }{\delta(c)/x^*} \frac{3\tau(c)^2}{c} \lambda\hat{g}'(x^*)\delta(c)}{\frac{(1-\delta(c)/x^*)^{-\tau(c)}}{4\delta(c)/x^*}} \notag\\
     & = \lambda \hat{g}(x^*) \frac{4\delta(c)/x^*}{(1-\delta(c)/x^*)^{-\tau(c)}} + \frac{12\tau(c)^2}{c}\lambda \hat{g}'(x^*)\delta(c) \notag\\
     & \stackrel{(e)}{\leq} \lambda \hat{g}(x^*)\frac{4\delta(c)}{c x^*} + \frac{12\tau(c)^2}{c}\lambda \hat{g}'(x^*)\delta(c) \notag\\
     & \stackrel{(f)}{\leq}  \left(\lambda \hat{g}(x^*)\frac{4}{x^*} + 12\lambda \hat{g}'(x^*) \right)\frac{\tau(c)^2 \delta(c)}{c} \notag\\
     &  \stackrel{(g)}{\leq} \lambda \left(\hat{g}(x^*)\frac{4}{x^*}  + 12\hat{g}'(x^*) \right) \left(\frac{3\log{c}}{\delta(c)/x^*}\right)^2 \frac{\delta(c)}{c} \notag \\
     & \leq 36  \lambda x^* \left(\hat{g}(x^*) + 3x^*\hat{g}'(x^*) \right) \frac{(\log{c})^2}{\delta(c)c}
\end{align}
where $(e)$ follows from \eqref{eq: minus tau bigger}; $(f)$ follows as $\tau(c)^2 \geq 1$; $(g)$ follows from \eqref{eq: bound tau}.

\subsubsection*{Step 4. Performance Loss of the Two-Price Policy for $g$.} Let $\mathcal{R}$ denote the reward obtained by this two-price policy under the original problem $g$, then for $c \geq \max\{4,N(1)\}$, and $\alpha \in (1,\infty)$ the performance loss is 
\begin{align*}
\lambda g(x^*)-\mathcal{R}& \stackrel{(a)}{\leq} \lambda  g(x^*) - \hat{\mathcal{R}} \\
& \stackrel{(b)}{\leq} \lambda (\hat{g}(x^*)  + k_1\delta(c)^\alpha) - \hat{\mathcal{R}} \\
& \stackrel{(c)}{\leq}  
    \lambda k_1\delta(c)^\alpha + 36  \lambda x^* \left(\hat{g}(x^*) + 3x^*\hat{g}'(x^*) \right) \frac{(\log{c})^2}{\delta(c)c} \\
    & \stackrel{(d)}{\leq}  \lambda \left(k_1 + 36 x^* \left(\hat{g}(x^*) + 3x^*\hat{g}'(x^*) \right) {(\log{c})^2} \right)c^{-\alpha/(1+\alpha)}\\
    & \stackrel{(e)}{=} \tilde{\mathcal{O}}(c^{1/(1+\alpha)}),
\end{align*}
where $(a)$ follows by the definition of $\hat{g}$ because $g \ge \hat g$; $(b)$ follows from 
$$\hat{g}(x^*) = \frac{1}{2}\hat{g}(x^*-\delta(c))  + \frac{1}{2}\hat{g}(x^*+\delta(c))= \frac{1}{2}{g}(x^*-\delta(c))  + \frac{1}{2}{g}(x^*+\delta(c)),$$
by definition and under \Cref{assmp: condition for upper bound}, $\hat{g}(x^*) \geq {g}(x^*) - k_1\delta(c)^\alpha$; $(c)$ follows from \eqref{eq:regret-hat2}; $(d)$ follows from $\delta(c) = c^{-1/(\alpha+1)}$ when $\alpha \in (1,\infty)$; $(e)$ follows as we scale $\lambda$ and $c$ proportionally while keeping $d$ fixed, thus $x^*$ is a fixed constant and $\lambda c^{-\alpha/(1+\alpha)} = \Theta(c^{1/(1+\alpha})$. 

For $\alpha = \infty$, we have $\delta(c) = \varepsilon$, and the performance loss is
\begin{align*}
    \lambda g(x^*)-\mathcal{R} & \stackrel{(e)}{=} \lambda \hat{g}(x^*)  - \hat{\mathcal{R}} \\
    & \stackrel{(f)}{=}  36  \lambda x^* \left(\hat{g}(x^*) + 3x^*\hat{g}'(x^*) \right) \frac{(\log{c})^2}{\delta(c)c}\\
    & \stackrel{(g)}{=} \mathcal{O}((\log{c})^2)
\end{align*}
where $(e)$ follows as $\alpha= \infty$, $\hat{g} = g$ by definition; $(f)$ follows from \eqref{eq:regret-hat2}; $(g)$ follows as we scale $\lambda$ and $c$ proportionally and $\delta(c)$ is a constant.

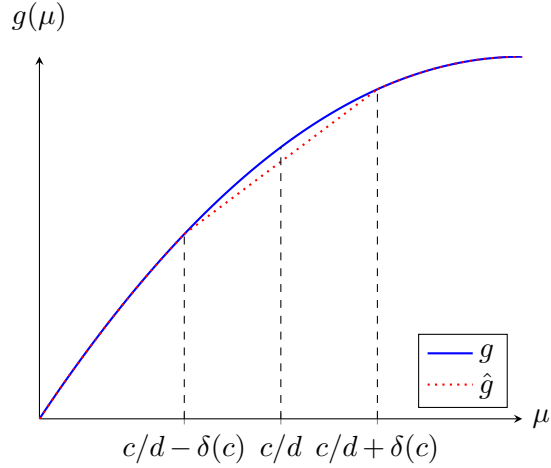
\begin{figure}
    \centering
\begin{tikzpicture}
    \begin{axis}[
        axis lines = middle,
        xlabel = $x$,
        ylabel = $g(x)$,
        ymin=0, ymax=1,
        xmin=0, xmax=1,
        xtick={0.3,0.5,0.7},
        xticklabels={$x^*-\delta(c)$,$x^*$,$x^*+\delta(c)$},
        ytick=\empty,
        width=8cm,
        height=6.4cm,
        legend pos= south east,
        every axis x label/.style={at={(current axis.right of origin)},anchor=west},
        every axis y label/.style={at={(current axis.north west)},above=2mm},
        ]
        \addplot[domain=0:1, samples=100, blue, thick] {x*(2-x)};
        \addlegendentry{${g}$}
        \addplot[domain=0:0.3, samples=100, dotted, red, thick] {x*(2-x)};
        \addplot[domain=0.7:1, samples=100, dotted, red, thick] {x*(2-x)};
        \addplot[domain=0.3:0.7, samples=100, dotted, red, thick] {(0.7*1.3-0.3*1.7)/0.4(x-0.3)+0.3*1.7};
        \addlegendentry{$\hat{g}$}
        \draw[densely dotted] (axis cs:0.5,0) -- (axis cs:0.5,0.75);
        \draw[densely dotted] (axis cs:0.3,0) -- (axis cs:0.3,0.3*1.7);
        \draw[densely dotted] (axis cs:0.7,0) -- (axis cs:0.7,0.7*1.3);
    \end{axis}
\end{tikzpicture}
    \caption{Illustration of $\hat{g}$ and $g$.}
    \label{fig: illus of hat{g}}
\end{figure}

\subsection{Performance Loss Lower Bound for General Stock-dependent Policies}\label{sec: lower bound}
A remaining question is whether more sophisticated stock-dependent policies can improve upon a two-price policy. We show that two-price policies attain the optimal performance loss within this class by providing a matching lower bound up to logarithmic factors on the performance loss of stock-dependent policies under the following condition, which is the analog of \Cref{assmp: condition for upper bound} where we are instead upper-bounding the reward function $g$ using a local approximation. We provide the formal proof in the next subsection.

{
\begin{assumption}\label{assmp: condition for lower bound}
There exists $0< \varepsilon < \min\{x^*, 1-x^*\}, \alpha \in [1,\infty], k_2 >0$, and a supergradient $g'(x^*) \geq 0$ of $g$ at $x^*$ such that $g(x)\leq g(x^*)+g'(x^*)(x-x^*)-k_2|x-x^*|^\alpha$ for all $x\in[x^*-\varepsilon,x^*+\varepsilon]$.
\end{assumption}
}
\Cref{assmp: condition for lower bound} is lower-bounding the difference between $g(x)$ and its approximation:
\[g(x^*)+g'(x^*)(x-x^*) - g(x) \geq k_2|x-x^*|^\alpha.\]
It is illustrated in \Cref{fig: illus of assmp lb} below.

\begin{figure}[!htbp]
\captionsetup{justification=centering}
\begin{subfigure}[b]{0.3\textwidth}
    \centering
    \begin{tikzpicture}
      \begin{axis}[
        axis lines = middle,
        xlabel = $x$,
        ylabel = $g(x)$,
        ymin=0, ymax=1,
        xmin=0, xmax = 1.1,
        xtick={0.001, 0.5, 1},
        xticklabels={0,$x^*$,1},
        ytick=\empty,
        width=5cm,
        height=4cm,
        every axis x label/.style={at={(current axis.right of origin)},anchor=west},
        every axis y label/.style={at={(current axis.north west)},above=2mm},   
        ]
        \addplot[domain=0:0.5, samples=100, blue, thick] {x};
        \addplot[domain=0.5:1, samples=100, blue, thick] {0.5*x+0.25};
        \node[below] at (axis cs:0.8,0.3) {$\Omega(\sqrt{c})$};
        \addplot[domain=0.25:0.5, samples=100, dashed, orange, thick] {0.8*(x-0.5)+0.5};
        \addplot[domain=0.5:0.75, samples=100, dashed, orange, thick] {0.7*(x-0.5)+0.5};
        \draw[densely dotted] (axis cs:0.5,0) -- (axis cs:0.5,0.5);
    \end{axis}
\end{tikzpicture}
\caption{$\alpha =1$}
\end{subfigure}
    \hspace{0.03\textwidth}
\begin{subfigure}[b]{0.3\textwidth}
    \centering
    \begin{tikzpicture}
    \begin{axis}[
        axis lines = middle,
        xlabel = $x$,
        ylabel = $g(x)$,
        ymin=0, ymax=1,
        xmin=0, xmax = 1.1,
        xtick={0.001, 0.5, 1},
        xticklabels={0,$x^*$,1},
        ytick=\empty,
        width=5cm,
        height=4cm,
        every axis x label/.style={at={(current axis.right of origin)},anchor=west},
        every axis y label/.style={at={(current axis.north west)},above=2mm},
        ]
        \addplot[domain=0:0.6, samples=100, blue, thick] {2*x*(1.2-x)};
        \addplot[domain=0.6:1, samples=100, blue, thick] {0.72};
        \node[below] at (axis cs:0.8,0.3) {$\Omega(c^{1/3})$};
        \addplot[domain=0.25:0.75, samples=100, dashed, orange, thick] {0.7 + 0.4*(x-0.5) - (x-0.5)^2};
        \draw[densely dotted] (axis cs:0.5,0) -- (axis cs:0.5,0.7);
    \end{axis}
\end{tikzpicture}
\caption{$\alpha =2$}
\end{subfigure}
    \hspace{0.03\textwidth}
\begin{subfigure}[b]{0.3\textwidth}
\begin{tikzpicture}
\centering
    \begin{axis}[
        axis lines = middle,
        xlabel = $x$,
        ylabel = $g(x)$,
        ymin=0, ymax=1,
        xmin=0, xmax = 1.1,
        xtick={0.001, 0.5, 1},
        xticklabels={0,$x^*$,1},
        ytick=\empty,
        width=5cm,
        height=4cm,
        every axis x label/.style={at={(current axis.right of origin)},anchor=west},
        every axis y label/.style={at={(current axis.north west)},above=2mm},   
        ]
        \addplot[domain=0:0.3, samples=100, blue, thick] {x};
        \addplot[domain=0.3:1, samples=100, blue, thick] {0.5*x + 0.15};
        \node[below] at (axis cs:0.8,0.3) {$\Omega(\log{c})$};
        \addplot[domain=0.2:0.75, samples=100, dashed, orange, thick] {0.5*x+0.15};
        \draw[densely dotted] (axis cs:0.5,0) -- (axis cs:0.5,0.4);
    \end{axis}
\end{tikzpicture}
\caption{$\alpha =\infty$}
\end{subfigure}
\caption{Illustration of \Cref{assmp: condition for lower bound}. The blue lines depict $g(x)$ while the orange dashed lines depict $g(x^*)+g'(x^*) (x-x^*)- k_2|x-x^*|^\alpha, x\in [x^*-\varepsilon, x^*+\varepsilon]$, in different cases.}
\label{fig: illus of assmp lb}
\end{figure}

\begin{theorem}\label{thm: lower bound}
Under \Cref{assmp: condition for lower bound}, if $\alpha \in [1, \infty)$, the performance loss of the optimal stock-dependent policy is ${{\Omega}}(c^{1/(\alpha+1)})$; if $\alpha = \infty$, there are three subcases ($\partial g(x)$ denote the set of supergradients of $g$ at $x\in [0,1]$):  
\begin{itemize}
        \item if $\partial g(x^*) = \{0\}$ (extremely scarce supply), the performance loss is $\Omega(e^{-c})$;
        \item if $\partial g(x) = \{r\}, \forall x \in [0,x^*+\varepsilon], \varepsilon >0$ (excess supply), the performance loss is $\Omega(1)$;
        \item otherwise, the performance loss is $\Omega(\log{c})$.
\end{itemize}
\end{theorem}

\textbf{Implications of \Cref{thm: lower bound}}. 
If a reward function satisfies \Cref{assmp: condition for upper bound} and \Cref{assmp: condition for lower bound} with the same $\alpha$, then both \Cref{thm: upper bound} and \Cref{thm: lower bound} apply. In this case, we obtain a tight convergence rate for the optimal stock-dependent policy with matching upper and lower bounds up to logarithmic factors. This holds for the three examples shown in Figures \ref{fig: illus of assmp lb}.
Recall that the proof of the performance loss upper bound was based on analyzing a two-price policy. This implies that a two-price policy attains the same asymptotic performance loss as optimal stock-dependent policies.
{Technically, the} values of $\alpha$ can be different for the two assumptions if the local curvatures of the reward function $g$ is different to the left and right of $x^*$, but this would be uncommon. {See Appendix \ref{apx: example of mismatch alpha} for an example.}

\subsection{Proof of \Cref{thm: lower bound}} \label{apx: lower bound}
We first show that when $\alpha = 1$, the performance loss of the optimal stock-dependent policy is lower bounded by $\Omega(\sqrt{c})$, and also depends on the supergradient values $\partial g(x^*)$. The proof relies on bounding the optimal objective value of the optimization problem \eqref{OPTSD}. Specifically, we begin by relaxing this primal problem, then we formulate its dual and construct a feasible solution to the dual. By weak duality, this provides an upper bound to the performance of stock-dependent policies, thus a lower bound to the performance loss. 

The performance loss for $\alpha \in (1,\infty)$ builds on the result for $\alpha = 1$. We first approximate the original reward function $g$ from above using a piecewise-linear function $\tilde{g}$ that has a kink at $x^*$. Thus by construction, the reward function $\tilde{g}$ satisfies \Cref{assmp: condition for lower bound} with $\alpha = 1$, allowing us to apply the performance loss bound derived earlier. 
We then tune the construction of $\tilde{g}$ to balance the trade-off between the approximation error and the performance loss on $\tilde{g}$ (i.e., via $\partial \tilde{g}(x^*)$ ) to obtain a $\Omega(c^{1/(\alpha+1)})$ performance loss bound on the original reward function $g$.
Finally, we consider $\alpha=\infty$ in both its general case and two corner cases, following a similar proof strategy.

\subsubsection{$\alpha = 1$ case.}\label{thm: degenerate case}
We first show the following performance loss results for a reward function that is upper bounded by two linear functions with different slope values.
\begin{lemma}\label{lem: sqrt c lower bound R r}
Suppose there exists $R>r>0$ such that $g(x) \leq \min\{Rx+g(x^*) - R x^*, r(x) + g(x^*) - rx^*\}$ for $x\in [0,1]$. For $c \geq 1/(R-r)$, the performance loss of the optimal stock-dependent policy under $g$ is lower bounded by $\lambda \min\{1, r/2\}x^* \sqrt{(R-r)/c}$.
\end{lemma}
\proof{Proof.}
As $g(x) \leq \min\{Rx+g(x^*) - R x^*, r(x) + g(x^*) - rx^*\}$, by \Cref{thm: optsd to lpsd}, we can formulate the following relaxed version of \eqref{LPSD} (where the constraints $\pi_j \geq (c-j+1)/(\lambda d) \pi_{j-1}, j \in [c]$ are relaxed), which provides an upper bound on the performance of the optimal stock-dependent policy:

\begin{align}
\textsf{LP}^* = \max_{\bm{\pi}\geq \bm{0}}\quad & \sum_{j=1}^c \lambda g_j \notag
\\ \text{s.t.}\quad &\sum_{j=0}^c \pi_j= 1, & && \text{(dual variable $\zeta$)}\notag
\\ &g_j \le \left(g(x^*) - Rx^*\right)\pi_j+R\frac{c-j+1}{\lambda d} \pi_{j-1}, &\forall j \in [c], && \text{(dual variable $\alpha_j$)}\notag
\\ &g_j \le \left(g(x^*) - rx^*\right)\pi_j+r\frac{c-j+1}{\lambda d} \pi_{j-1}, &\forall j \in [c]. && \text{(dual variable $\beta_j$)} \label{deg primal}
\end{align}

The dual problem of the optimization problem \eqref{deg primal} is
\begin{align}
\min\quad &\zeta \notag
\\ \text{s.t.}\quad &\alpha_j+\beta_j= \lambda, & \forall j \in [c], && \text{(primal variable $g_j$)} \notag
\\ &\zeta \ge \lambda g(x^*)-Rx^*\alpha_j-rx^*\beta_j+R\frac{c-j}{\lambda d}\alpha_{j+1}+r\frac{c-j}{\lambda d}\beta_{j+1}, &\forall j >0, && \text{(primal variable $\pi_j$)} \notag
\\ &\zeta \ge Rx^*\alpha_1+rx^*\beta_1, & && \text{(primal variable $\pi_0$)} \notag
\\ &\alpha_j,\beta_j \ge 0,  &\forall j \in [c] \label{deg dual}, \notag \\
    & \alpha_{c+1}= \beta_{c+1} = 0.
\end{align}

We present a feasible solution to the dual problem \eqref{deg dual}. The proof can be found in \Cref{apx: feasible solution to dual alpha 1}.
\begin{proposition}\label{pro: feasible solution to dual alpha 1}
For $c \geq 1/(R-r)$, the following is a feasible solution to the dual problem \eqref{deg dual}:
    \[\zeta = \lambda g(x^*) - \lambda \min\{1, r/2\} x^*\sqrt{(R-r)/c}, \alpha_j = \lambda \max\{1-j/\sqrt{c(R-r)}, 0\}, \beta_j = \lambda - \alpha_j, \forall j \in [c].\]
\end{proposition}

Therefore, let $\mathcal{R}^*$ denote the performance of the optimal stock-dependent policy under $g$, then we have  
\begin{align*}
    \mathcal{R}^* \stackrel{(a)}{\leq} \textsf{LP}^* \stackrel{(b)}{\leq}\zeta = \lambda g(x^*) - \lambda \min\{1, r/2\} x^*\sqrt{(R-r)/c}.
\end{align*}
where $(a)$ follows from \Cref{thm: optsd to lpsd}; $(b)$ follows from weak duality. The lower bound of the performance loss for the optimal stock-dependent policy follows.
\Halmos 
\endproof

Note for any $g$ that is non-differentiable at $x^*$, we can always find a piecewise linear function that upper bounds $g$ with a kink at $x^*$ with the same $\FLU$ value, i.e., $g(x) \leq \bar{g}(q) = \min\{Rx + g(x^*)-Rx^*, rx + g(x^*)-rx^*\}$ for some $R>r\geq 0$. See \Cref{fig: illus of bar{g}} for an illustration. Therefore, the above bound applies to reward functions that satisfy \Cref{assmp: condition for lower bound} with $\alpha = 1$.

\begin{figure}
    \centering
        \begin{tikzpicture}
    \begin{axis}[
        axis lines = middle,
        xlabel = $x$,
        ylabel = $g(x)$,
        ymin=0, ymax=1,
        xmin=0, xmax = 1.1,
        xtick={0.001, 0.5, 1},
        xticklabels={0,$x^*$,1},
        ytick=\empty,
        width=8cm,
        height=6.4cm,
        legend pos= south east,
        every axis x label/.style={at={(current axis.right of origin)},anchor=west},
        every axis y label/.style={at={(current axis.north west)},above=2mm},
        ]
        \draw[densely dotted] (axis cs:0.5,0) -- (axis cs:0.5,0.5);
        \addplot[domain=0.5:1, samples=100, blue, thick] {0.3*x+0.35};
        \addlegendentry{$g$}
        \addplot[domain=0:0.5, samples=100, red, dotted, thick] {0.5*(x-0.5)+0.5};
        \addlegendentry{$\bar{g}$}
        \addplot[domain=0.5:1, samples=100, red, dotted, thick] {0.4*x+0.3};
        \addplot[domain=0:0.5, samples=100, blue, thick] {x*(1.5-x)};
\end{axis}
\end{tikzpicture}
    \caption{Illustration of $\bar{g}$ and $g$, $\alpha = 1$. {Unlike the first picture in \Cref{fig: illus of assmp lb}, here we illustrate on a curve $g$ that is curved to the left of $x^*$ and straight to the right of $x^*$.}}
    \label{fig: illus of bar{g}}
\end{figure}
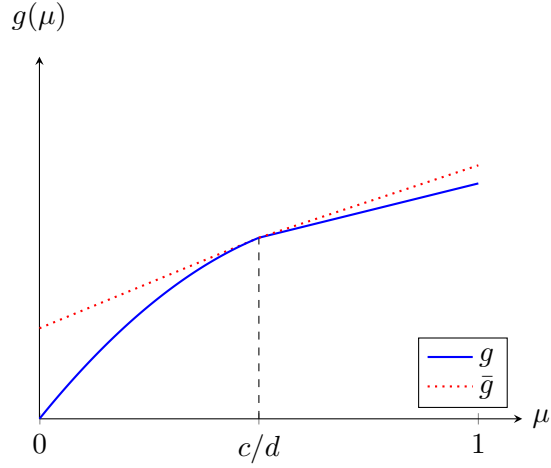

\subsubsection{$\alpha \in (1,\infty)$ case.}
Define $h(x) = g(x^*)+g'(x^*)(x-x^*)-k_2|x-x^*|^\alpha$. By \Cref{assmp: condition for lower bound}, we have $h(x)\geq g(x)$. We consider a piecewise-linear approximation of $g$, denoted by $\tilde{g}$, defined as (see \Cref{fig: illus of tilde{g}} for illustration):
$$
\tilde{g}(x) = \left\{
\begin{aligned}
& h'(x^*-\eta(c))(x-(x^*-\eta(c))) + h(x^* - \eta(c)), & 0\leq x\leq x^*,\\
& g'(x^*)(x-x^*) + h(x^*-\eta(c)) + h'(x^*-\eta(c))\eta(c), & x^* < x \leq 1.\\
\end{aligned}
\right.
$$
where $\eta(c)>0$ will be determined later to optimize the bound.

We first show $\tilde{g}(x) \geq g(x)$ for all $x \in [0,1]$:
\begin{itemize}
    \item For $0\leq x\leq x^*$, 
\[h'(x^*-\eta(c))(x-(x^*-\eta(c))) + h(x^* - \eta(c)) \stackrel{(a)}{\geq} h(x) \stackrel{(b)}{\geq} g(x). \]
where $(a)$ follows from the concavity of $h$; $(b)$ follows from \Cref{assmp: condition for lower bound}.
    \item For $x^* <x \leq 1$, 
     \[g(x) \stackrel{(c)}{\leq} g'(x^*)(x-x^*) + g(x^*) \stackrel{(d)}{\leq}  g'(x^*)(x-x^*) + h(x^*-\eta(c))+h'(x^*-\eta(c))\eta(c), \]
    where $(c)$ follows from the concavity of $g$; $(d)$ follows from definition of $h$:
\begin{align*}
    h(x^*-\eta(c))+h'(x^*-\eta(c))\eta(c) & = g(x^*) - g'(x^*)\eta(c) - k_2\eta(c)^\alpha +(g'(x^*) + \alpha k_2\eta(c)^{\alpha-1})\eta(c) \\
    & = g(x^*) + (\alpha-1)k_2\eta(c)^\alpha.
\end{align*}
 
\end{itemize}

Let $\mathcal{R}^*$ and $\tilde{\mathcal{R}}^*$ denote the performance of the optimal stock-dependent policy under $g$ and $\tilde{g}$ respectively. The performance loss of the original problem under the optimal stock-dependent policy is
\begin{align*}
    \lambda  g(x^*)-\mathcal{R}^*   &  \stackrel{(a)}{\geq} \lambda g(x^*)-\tilde{\mathcal{R}}^* \\
     & = \lambda g(x^*) - \lambda \tilde{g}(x^*) + \lambda \tilde{g}(x^*) - \tilde{\mathcal{R}}^* \\
     & \stackrel{(b)}{\geq} - \lambda (\alpha-1)k_2\eta(c)^\alpha + \lambda \tilde{g}(x^*) - \tilde{\mathcal{R}}^* \\
    & \stackrel{(c)}{=} - \lambda (\alpha-1)k_2\eta(c)^\alpha + \lambda \min\{1,g'(x^*)/2\}x^*\sqrt{\frac{\alpha k_2\eta(c)^{\alpha-1}}{c}}\\
    & \stackrel{(d)}{=} \Omega(c^{1/(\alpha+1)}),
\end{align*}
where $(a)$ follows as $\tilde{g}(x) \geq g(x)$ for all $x\in \mathbb{R}^+$ thus $\mathcal{R}^* \leq \tilde{\mathcal{R}}^*$; $(b)$ follows from the definition of $\tilde{g}$: $\tilde{g}(x^*)=h(x^*-\eta(c))+h'(x^*+\eta(c))\eta(c) = g(x^*)+(\alpha-1) k_2 \eta(c)^\alpha$; $(c)$ follows from \Cref{lem: sqrt c lower bound R r} where $r = g'(x^*)$, and $(R-r) = h'(x^*-\eta(c)) - g'(x^*) = \alpha k_2\eta(c)^{\alpha-1}$ for $\tilde{g}$, and the condition $c \geq 1/(R-r)$ holds for $c$ large enough when we choose $\eta(c)^* = \Theta(c^{-1/(\alpha+1)})$ in $(d)$; $(d)$ follows as we scale $\lambda$ and $c$ proportionally while keeping $d$ fixed, thus $x^*$ is a fixed constant and $\lambda c^{-\alpha/(1+\alpha)} = \Theta(c^{1/(1+\alpha})$; therefore optimizing over $\eta(c)$, we have $\eta(c)^* = \Theta(c^{-1/(\alpha+1)})$ and the lower bound follows.

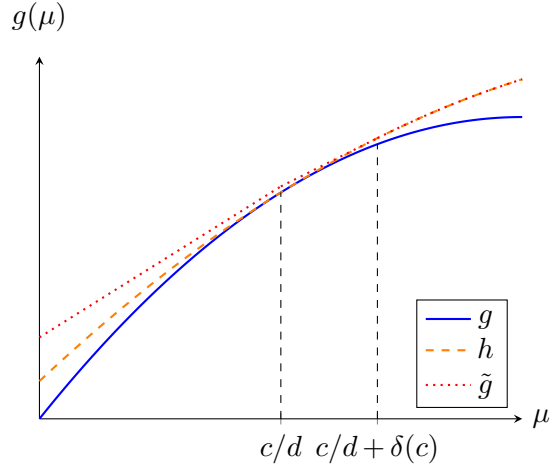
\begin{figure}
    \centering
\begin{tikzpicture}
    \begin{axis}[
        axis lines = middle,
        xlabel = $x$,
        ylabel = $g(x)$,
        ymin=0, ymax=1.2,
        xmin=0, xmax=1,
        xtick={0.001, 0.3, 0.5, 1},
        xticklabels={0, $x^* - \eta(c)$,$x^*$, 1},
        ytick=\empty,
        width=8cm,
        height=6.4cm,
        legend pos= south east,
        every axis x label/.style={at={(current axis.right of origin)},anchor=west},
        every axis y label/.style={at={(current axis.north west)},above=2mm},
        ]
        \addplot[domain=0:0.6, samples=100, blue, thick] {2*x*(1.2-x)};
        \addlegendentry{$g$}
         \addplot[domain=0.1:0.75, samples=100, dashed, orange, thick] {0.7 + 0.4*(x-0.5) - (x-0.5)^2};
         \addlegendentry{$h$}
          \addplot[domain=0.6:1, samples=100, blue, thick] {0.72};
        \addplot[domain=0:0.5, samples=100, dotted, red, thick] {0.8*(x-0.3) + 0.58};
        \addplot[domain=0.5:1, samples=100, dotted, red, thick] {0.4*(x-0.5) + 0.74};
        \addlegendentry{$\tilde{g}$}
        \draw[densely dotted] (axis cs:0.5,0) -- (axis cs:0.5,0.74);
        \draw[densely dotted] (axis cs:0.3,0) -- (axis cs:0.3,0.58);
    \end{axis}
\end{tikzpicture}
    \caption{Illustration of $\tilde{g}$, $h$, and $g$.}
    \label{fig: illus of tilde{g}}
\end{figure}

\subsubsection{$\alpha=\infty$ case.} See \Cref{apx: lower bound infty}.

\subsection{On the Power of Two Prices}\label{sec: two-price}

In this section, we discuss 1) the intuition of why a minimally dynamic policy with two prices can give us the optimal performance loss within the class of stock-dependent policies; and 2) how the local shape of reward function $g$ affects the shape parameter $\alpha$ and hence the performance of a two-price policy.
The following corollary is a direct result of combining \Cref{thm: upper bound} and \Cref{thm: lower bound}, which implies that a little dynamicity can go a long way.

\begin{corollary}
Under Assumptions~\ref{assmp: condition for upper bound} and~\ref{assmp: condition for lower bound}, the optimal performance loss $\tilde{\Theta}(c^{1/(\alpha+1)})$ of stock-dependent policies can be achieved by a two-price policy.
\end{corollary}

\textbf{Intuition on the benefits of dynamicity}. Our analysis of the performance of the fluid policy suggests that the stock-out probability $\pi_0$ is a key driver of performance. By dynamically adjusting the admission probability based on the stock level, our two-price policy reduces this loss from stock-outs.

{We further illustrate this intuition by making a concrete comparison assuming that the usage duration is exponentially distributed (with mean $d$).}
Recall that for a stock-dependent policy, when the stock level is $j>0$ the rate at which the stock level decreases is $\lambda x_j$ (the arrival rate times the admission probability) while the rate at which the stock level increases is $(c-j)/d$ {(because there are $c-j$ units in use, each of which returns at the rate of $1/d$ assuming exponential distributions)}. Let $S_t$ denote the number of units of stock available at time $t$. Therefore, the drift is $\lim_{h \downarrow 0} \mathbb E[S_{t+h} - S_t | S_t = j] / h = (c-j)/d -\lambda x_j$. In particular, the fluid policy induces a drift of 
$(c-j)/d- \lambda x^* = (c-j)/d - c/d = -j/d$ on the stock level process---a drift that is negative irrespective of $j$, pushing the stock level to zero (i.e., stock-out) and leading to revenue loss. On the contrary, a two-price policy, when properly parameterized, induces a mean reversion process that concentrates the stock level around the switching threshold: when the stock level is lower than the switching threshold $\tau$, it sets $x_L < c/(\lambda d)$ so that the system experiences a positive drift and runs out of stock less often; otherwise, it sets $x_H > c/(\lambda d)$ to collect rewards at a faster rate by admitting consumers more aggressively. 

More specifically, the parameter $\delta(c)$ controls the deviation from the fluid policy: $x_L = x^* - \delta(c)$ and $x_H = x^* + \delta(c)$. Setting $\delta(c)$ too small reduces the policy's ability to adjust the admission probabilities effectively in response to the system state, leading to a higher probability of stock-outs ($\pi_0$). In the extreme case where we set $\delta(c)=0$, the two-price policy becomes the fluid policy, with loss from stock-outs being $\Theta(c^{-1/2})$. On the other hand, setting $\delta(c)$ too large increases the difference between $g(x^*)$ and $g(x_L) (g(x_H))$, resulting in a higher loss from concavity. The tuning of $\delta(c)$ aims to balance these two types of losses. The switching threshold $\tau$ determines when the policy switches between the high and low admission probabilities based on the current stock level. As argued above, it represents the center of a mean reversion process, thus it should be set large enough to ensure a low probability of the stock-out, but not too large that an excessive number of units are idle.

As an interesting example, when the WTP distribution $F$ is discrete and never equals $x^*$ (i.e.\ $g$ is locally linear at $x^*$), the fluid policy randomizes between two prices to have a static admission probability $x^*$. By contrast, our two-price policy distinguishes between these two prices and dynamically chooses between them based on the stock level.

\textbf{Intuition on how the local shape of reward function affects performance loss.} 
By \Cref{thm: lower bound}, when $\alpha=1$ (reward function is locally non-differentiable), the optimal performance loss of stock-dependent policies is $\Omega(\sqrt{c})$, no better than the fluid policy. We now explain why this is the hard case and how, in general, the local shape of the reward function drives performance loss.

Let $\pi_L$ and $\pi_H$ denote the steady-state probability that the system is under low and high admission probabilities respectively. By definition, we have $\pi_0 + \pi_L +\pi_H = 1$. The performance loss of the two-price policy can be decomposed into two parts:
\begin{align}
    g(x^*) - (\pi_L g(x_L) + \pi_H g(x_H)) = \underbrace{\pi_0 g(x^*)}_{\text{Loss from stock-outs}} + \underbrace{\pi_L \left(g(x^*)-g(x_L)\right) - \pi_H \left(g(x_H)-g(x^*)\right)}_{\text{Loss from concavity}}. \label{two-price regret}
\end{align}
First note that under the two-price policy with parameters specified in \eqref{eq: two price para}, we would have $\pi_L/\pi_H \to 1$ as $c$ and $\lambda$ go to infinity.
Recall from our sketch of \Cref{thm: upper bound} that $x_L=x^*-\delta(c)$ and $x_H=x^*+\delta(c)$ for some $\delta(c)>0$.  By concavity of $g$, it will be the case that $\big(g(x^*)-g(x_L)\big)-\big(g(x_H)-g(x^*)\big)\ge0$.  The magnitude of this difference affects how well we can control the "loss from concavity" by calibrating $\pi_L$ and $\pi_H$ (via controlling $\tau(c)$ and $\delta(c))$.  In particular, this difference is 0 when $g$ is locally linear ($\alpha=\infty$), leading to the performance loss of $\mathcal{O}((\log{c})^2)$.  When $g$ is differentiable but not linear at $x^*$, this difference shrinks sublinearly in $\delta(c)$, with the parameter $\alpha$ affecting the sublinear rate (and hence factoring into our performance loss).  Finally, when $g$ is non-differentiable at $x^*$, this difference shrinks linearly in $\delta(c)$ and prevents us from achieving better than $\mathcal{O}(\sqrt{c})$ performance loss.

Recall that one example of a locally non-differentiable reward function is when the WTP distribution $F$ is discrete and equals $x^*$ at some valuation. In the context of non-reusable resource allocation, this situation is known as a degenerate case, and it is known that the fluid benchmark leads to an unavoidable square-root loss in this case \citep{bumpensanti_re-solving_2018, vera_bayesian_2020}. Our findings highlight that degeneracy remains a challenging case for reusable resources: when the usage durations are exponential, the optimal dynamic policy (DP) lies in our stock-dependent policy class and still attains $\Omega(\sqrt{c})$ performance loss.
Therefore, achieving a better performance loss than $\Omega(\sqrt{c})$ would require not only a different policy but a different benchmark for reusable resources.
We defer further investigation to future research.

\section{Extensions}\label{sec: extension}
In the {base} model, we implicitly assume that the usage duration distribution conditional on the admission probability $x$ remains constant with mean $d$. In this section, we demonstrate how this assumption can be relaxed by extending our results to a model with multiple customer classes and multiple resource types. Specifically, we show that by assigning dedicated resources to each customer class based on the fluid relaxation and implementing two-price policies for each class separately, we can also achieve $o(\sqrt{c})$ performance loss. This highlights the power of two prices in more general settings.

\subsection{Extension to Multiple Customer Classes}\label{sec: multiple customer class}
There is a single reusable resource with $c$ units of stock to serve $m$ classes of customers. Each class $i\in [m]$ is associated with a reward function $g_i$ and a usage duration distribution $G_i$ with mean $d_i$. Customers from each class $i$ arrive following an independent Poisson process with rate $\lambda_i$. We assume that the decision-maker sets personalized (type-specific) admission probabilities $\{x_i(t)\}_{i\in [m]} \in [0,1]$ over a continuous time, and when there is no stock available, $x_i(t) = 0$ for all classes $i\in [m]$. The asymptotic regime can be similarly defined as before, i.e., scale the arrival rates of each customer class $\{\lambda_i\}_{i \in [m]}$ and the total number of stock $c$ proportionally, while keeping the usage duration distributions and means $\{d_i\}_{i\in [m]}$ fixed.

The fluid relaxation is defined as follows. 
\begin{align*}
    \mathsf{FLU}(m)=\max_{x_i \in [0,1], i \in [m]}\ & \sum_{i\in [m]} \lambda_i g_i(x_i) \notag \\
    \text{s.t. }& \sum_{i\in [m]}\lambda_i x_i d_i \leq c.
\end{align*}
Similar as before, we assume the resource is scarce in that  $c/(\sum_{i\in[m]}\lambda_i d_i)<1$ as the problem is essentially trivial otherwise.

The model with multiple customer classes introduced complexities not present in the base model. As usage duration is correlated with WTP, the decision-maker faces a dual trade-off: between serving current and future customers, and among different customer classes. For example, suppose class 1's WTP first-order stochastically dominates that of class 2, and its duration distribution $G_1$ is first-order stochastically dominated by $G_2$. Admitting a customer of group 1 might yield a higher reward, but can also incur a longer service time, leading to a higher risk of missing future customers. Nevertheless, we show how to extend our results for the case of multiple customer classes by assigning resources to each customer class and implementing a two-price policy separately for each customer class.  

\begin{proposition}\label{pro: performance loss for multiple customer class}
If the reward function $g_i$ for each customer class $i\in [m]$ satisfies \Cref{assmp: condition for upper bound} with $\alpha_i$, there exists a policy with performance loss $\tilde{O}(c^{1/(1+\underline{\alpha})})$, where $\underline{\alpha} = \min_{i\in [m]} \alpha_i$. More specifically, this policy assigns dedicated resources to each customer class based on $\mathsf{FLU}(m)$, and implements a two-price policy for each class separately. 
\end{proposition}
\proof{Proof.}
Let $\{x_i^*\}_{i\in [m]}$ denote the optimal solution to $\mathsf{FLU}(m)$, i.e., $\mathsf{FLU}(m) = \sum_{i\in [m]} \lambda_i g_i(x_i^*)$. Consider a policy that assigns dedicated resource level $c_i := \lfloor\lambda_i x_i^* d_i\rfloor$ to each customer class $i\in [m]$. Under this policy, each customer class can only use resources from their dedicated assignment and a customer class would stock out if none of their resources are available, even if other classes have resources available. 

By the feasibility of $\{x_i^*\}_{i\in [m]}$, we have $\sum_{i\in [m]} \lambda_i x_i^* d_i \leq c$, thus the dedicated resource assignment is feasible. Therefore, we reduce the problem to $m$ decoupled single-customer-class problems with resources $c_i$, whose fluid relaxation benchmark is 
\begin{align*}
    \mathsf{FLU}_i=\max_{x_i \in [0,1]}\ &  \lambda_i g_i(x_i) \\
    \text{s.t. }&\lambda_i x_i d_i \leq c_i.
\end{align*}
Again by the feasibility of $\{x_i^*\}_{i\in [m]}$, we have $x_i^* \leq 1$, thus $\mathsf{FLU}_i = \lambda_i g_i(c_i/(\lambda_i d_i))$ for $i\in [m]$. 

Therefore, by \Cref{thm: upper bound}, there exists a two-price policy for each class that yields a reward of $\lambda_i g_i(c_i/(\lambda_i d_i)) - \tilde{\mathcal{O}}(c_i^{1/(1+\alpha_i)})$. Note that 
\[\lambda_i g_i\left(\frac{c_i}{\lambda_i d_i}\right) \stackrel{(a)}{\geq} \lambda_i g_i\left (x^*-\frac{1}{\lambda_i d_i}\right) \stackrel{(b)}{\geq}  \lambda_i g(x_i^*) - \frac{g_i'(0)}{d_i}, \]
where $(a)$ holds as $\lambda_i x_i^* d_i - c_i \leq 1$; $(b)$ holds by the concavity of $g$.

Putting things together, the performance loss of the policy is 
\begin{align*}
    &\mathsf{FLU}(m) - \sum_{i\in [m]}\left( \lambda_i g_i(c_i/(\lambda_i d_i)) - \tilde{\mathcal{O}}(c_i^{1/(1+\alpha_i)})\right) \\
    & \quad \quad \leq  \mathsf{FLU}(m) - \sum_{i\in [m]} \left(\lambda_i g(x_i^*) - \frac{g_i'(0)}{d_i} - \tilde{\mathcal{O}}(c_i^{1/(1+\alpha_i)})\right)\\
    & \quad \quad = \sum_{i\in [m]}  \frac{g_i'(0)}{d_i}   + \tilde{\mathcal{O}}(c_i^{1/(1+\alpha_i)})\\
    & \quad \quad= \tilde{O}(c^{1/(1+\underline{\alpha})}). \Halmos
\end{align*}

\endproof

\subsection{Extension to Multiple Resource Types}
The performance loss result as stated in \Cref{pro: performance loss for multiple customer class} can also be extended to a network revenue management model. Suppose there are $n$ resource types, each consisting of $c_j, j\in [n]$ units of stock, serving $m$ classes of customer. If a customer enters service, then $a_{i,j} \ge 0$ units of resource $j$ are used. The asymptotic regime is to scale the arrival rates of each customer class $\{\lambda_i\}_{i \in [m]}$ and the number of stock $\{c_j\}_{j\in [n]}$ proportionally, while keeping the usage duration distributions and means $\{d_i\}_{i\in [m]}$ fixed. In this model, the fluid relaxation is defined as 
\begin{align*}
    \mathsf{FLU}(m,n)=\max_{x_i \in [0,1], i \in [m]}\ & \sum_{i\in [m]} \lambda_i g_i(x_i) \notag \\
    \text{s.t. }& \sum_{i\in [m]} \lambda_i x_i a_{i,j} d_i \leq c_j, \forall j \in [n].
\end{align*}

We consider the following policy: let $\{x_i^*\}_{i\in [m]}$ denote the optimal solution to $\mathsf{FLU}(m,n)$. For each resource type $j$, we dedicate $\lfloor \lambda_i x_i^* a_{i,j} d_i \rfloor$ units to customer class $i$. Then for each customer type $i$, {the total number that can be served at one time (using $\lfloor \lambda_i x_i^* a_{i,j} d_i \rfloor$ of each resource $j$) is}
\begin{align} \label{eqn:12938}
\min_{j\in [n], a_{i,j} >0} \left\lfloor \frac{\lfloor \lambda_i x_i^* a_{i,j} d_i \rfloor}{a_{i,j}}\right\rfloor.
\end{align}
We run a two-price policy for each customer class $i$ with the total stock {given by expression~\eqref{eqn:12938}}. Then the results from \Cref{sec: multiple customer class} apply, and we will obtain the same performance loss.

\begin{proposition}\label{pro: performance loss for network model}
If the reward function $g_i$ for each customer class $i\in [m]$ satisfies \Cref{assmp: condition for upper bound} with $\alpha_i$, there exists a policy with performance loss $\tilde{O}(c^{1/(1+\underline{\alpha})})$, where $\underline{\alpha} = \min_{i\in [m]} \alpha_i$. More specifically, this policy assigns dedicated resources to each customer class based on $\mathsf{FLU}(m,n)$, and implements a two-price policy for each class separately. 
\end{proposition}

In the above models, we assume each customer class consumes a fixed combination of resources, which allows us to decouple the problem and manage each class independently by assigning dedicated resources. Our decomposition does not work when product substitution is allowed, for example, when customers choose products according to a multinomial logit choice model. In that case, as products share resources that have fixed supply, the decision-maker must jointly decide what products to offer and what prices to quote. We leave this extension to future work. 

\section{Numerical Examples}\label{sec: numerics}
In this section, we numerically examine the performance of the optimal\footnote{Solved \eqref{OPTSD} numerically, following \cite{beck2003mirror}.} stock-dependent policy (denoted as SDOPT in the legend in Figures \ref{fig: compare convergence rate} and \ref{fig: compare reward proportion}), our two-price policy (denoted as TwoPriceOPT when the parameters are optimized and TwoPrice otherwise), a static policy with optimized admission probabilities\footnote{Optimized using the Sequential Least Squares Programming (\textsf{SLSQP}) algorithm in \textsf{scipy.optimize.minimize}.} parameter (denoted as StaticOPT), and the fluid policy that sets the admission probability to be $x^*$ whenever there are units of stock available; denoted as FluidPolicy\footnote{Codes for the simulation can be found \href{https://github.com/wenxinzhang0/Dynamic-Pricing-for-Reusable-Resources}{here}.}. 

We use two sets of performance metrics: the performance loss compared to $\FLU$ for large $c$, and the proportion of reward obtained compared to $\FLU$ for small $c$. The first metric is the focus of our theoretical results and the numerical examples corroborate them---optimal stock-dependent policies perform better than the fluid policy except for $g$ non-differentiable at $x^*$, and a two-price policy can achieve the optimal performance loss of stock-dependent policies. The second metric complements our analysis and illustrates how these policies perform when $c$ is small.

\subsection{Asymptotic Performance}
\begin{figure}[!htbp]
    \centering
    \begin{subfigure}[b]{0.49\linewidth}
        \begin{tikzpicture}[scale = 0.9]
    \begin{axis}[
        xlabel={No. of Supply Units},
        ylabel={Difference from Target},
        grid=major,
        legend entries={SDOPT, StaticOPT, FluidPolicy},
        legend pos=north west,
        cycle list name=color list,
        xmode=log,
        ymode=log,
        xtick={1000, 2000, 5000, 10000, 50000},
        xticklabels={1000, 2000, 5000, 10000, 50000}
    ]

    \addplot+[mark=triangle*, color=red] coordinates {
        (1000, 2000.0 - 1952.0716448415947)
        (2000, 4000.0 - 3931.8688156598455)
        (5000, 10000.0 - 9891.780048546492)
        (10000, 20000.0 - 19846.599213410056)
        (50000, 100000.0 - 99655.92256855163)
    };
    \addplot+[mark=diamond*, color=black] coordinates {
        (1000, 2000.0 - 1950.3761643959406)
        (2000, 4000.0 - 3929.476769193205)
        (5000, 10000.0 - 9888.006415377964)
        (10000, 20000.0 - 19841.268731311116)
        (50000, 100000.0 - 99644.02256303414)
    };
    
    \addplot+[mark=square*, color=blue] coordinates {
        (1000, 2000.0 - 1950.3761647076792)
        (2000, 4000.0 - 3929.4767698809305)
        (5000, 10000.0 - 9888.006417214945)
        (10000, 20000.0 - 19841.26873502389)
        (50000, 100000.0 - 99644.0225798464)
    };

    \end{axis}
\end{tikzpicture}
        
        \caption{Locally non-differentiable $g$}
        \label{fig: degenerate g}
    \end{subfigure}
    \begin{subfigure}[b]{0.49\linewidth}
    \begin{tikzpicture}[scale = 0.9]
    \begin{axis}[
        xlabel={No. of Supply Units},
        ylabel={Difference from Target},
        grid=major,
        legend entries={SDOPT, TwoPrice, StaticOPT, FluidPolicy},
        legend pos=north west,
        cycle list name=color list,
        xmode=log,
        ymode=log,
        xtick={1000, 2000, 5000, 10000},
        xticklabels={1000, 2000, 5000, 10000}
    ]

    \addplot+[mark=triangle*, color=red] coordinates {
        (1000, 1500.0 - 1475.7289654697286)
        (2000, 3000.0 - 2968.2515841207023)
        (5000, 7500.0 - 7454.903203839796)
        (10000, 15000.0 - 14940.886344063769)
        (50000, 75000.0 - 74885.26526334464)
    };
    \addplot+[mark=o, color=orange] coordinates {
        (1000, 1500.0 - 1465.0083279285398)
        (2000, 3000.0 - 2952.2670262616857)
        (5000, 7500.0 - 7431.513756187085)
        (10000, 15000.0 - 14907.337992939185)
        (50000, 75000.0 - 74804.6027776065)
    };

    \addplot+[mark=diamond*, color=black] coordinates {
        (1000, 1500.0 - 1462.9096273587866)
        (2000, 3000.0 - 2947.276524524947)
        (5000, 7500.0 - 7416.254800613539)
        (10000, 15000.0 - 14881.289624330548)
        (50000, 75000.0 - 74733.755666975)
    };
    \addplot+[mark=square*, color=blue] coordinates {
        (1000, 1500.0 - 1462.782123530759)
        (2000, 3000.0 - 2947.1075774106926)
        (5000, 7500.0 - 7416.0048129112065)
        (10000, 15000.0 - 14880.951551267912)
        (50000, 75000.0 - 74733.0169348843)
    };
    \end{axis}
\end{tikzpicture}
        
        \caption{Locally smooth $g$}
        \label{fig: smooth g}
    \end{subfigure}
    \begin{subfigure}[b]{0.49\linewidth}
        \begin{tikzpicture}[scale = 0.9]
    \begin{axis}[
        xlabel={No. of Supply Units},
        ylabel={Difference from Target},
        grid=major,
        legend entries={SDOPT, TwoPrice, StaticOPT, FluidPolicy },
        legend pos=north west,
        cycle list name=color list,
        xmode=log,
        ymode=log,
        xtick={1000, 2000, 5000, 10000, 50000},
        xticklabels={1000, 2000, 5000, 10000, 50000}
    ]
    \addplot+[mark=triangle*, color=red] coordinates {
        (1000, 2100.0 - 2085.5872904903654)
        (2000, 4200.0 - 4184.395086004221)
        (5000, 10500.0 - 10482.762444482547)
        (10000, 21000.0 - 20981.58512471086)
        (50000, 105000.0 - 104978.8342575474)
    };
    \addplot+[mark=o, color=orange] coordinates {
        (1000, 2100.0 - 2084.715337068959)
        (2000, 4200.0 - 4183.0397628581395)
        (5000, 10500.0 - 10481.353264328905)
        (10000, 21000.0 - 20979.844241451905)
        (50000, 105000.0 - 104978.04459714513)
    };
    \addplot+[mark=diamond*, color=black] coordinates {
        (1000, 2100.0 - 2049.031840031701)
        (2000, 4200.0 - 4127.493201815741)
        (5000, 10500.0 - 10384.758803031636)
        (10000, 21000.0 - 20836.59628322627)
        (50000, 105000.0 - 104633.33903074022)
    };
    \addplot+[mark=square*, color=blue] coordinates {
        (1000, 2100.0 - 2047.8949729430633)
        (2000, 4200.0 - 4125.9506083749775)
        (5000, 10500.0 - 10382.406738075693)
        (10000, 21000.0 - 20833.33217177508)
        (50000, 105000.0 - 104626.22370883872)
    };
    \end{axis}
\end{tikzpicture}
        \caption{Locally linear $g$}
        \label{fig: linear g}
    \end{subfigure}
    \caption{Comparing the performance of SDOPT, TwoPrice, StaticOPT, and FluidPolicy in terms of their performance loss when $c$ is large. }
    \label{fig: compare convergence rate}
\end{figure}
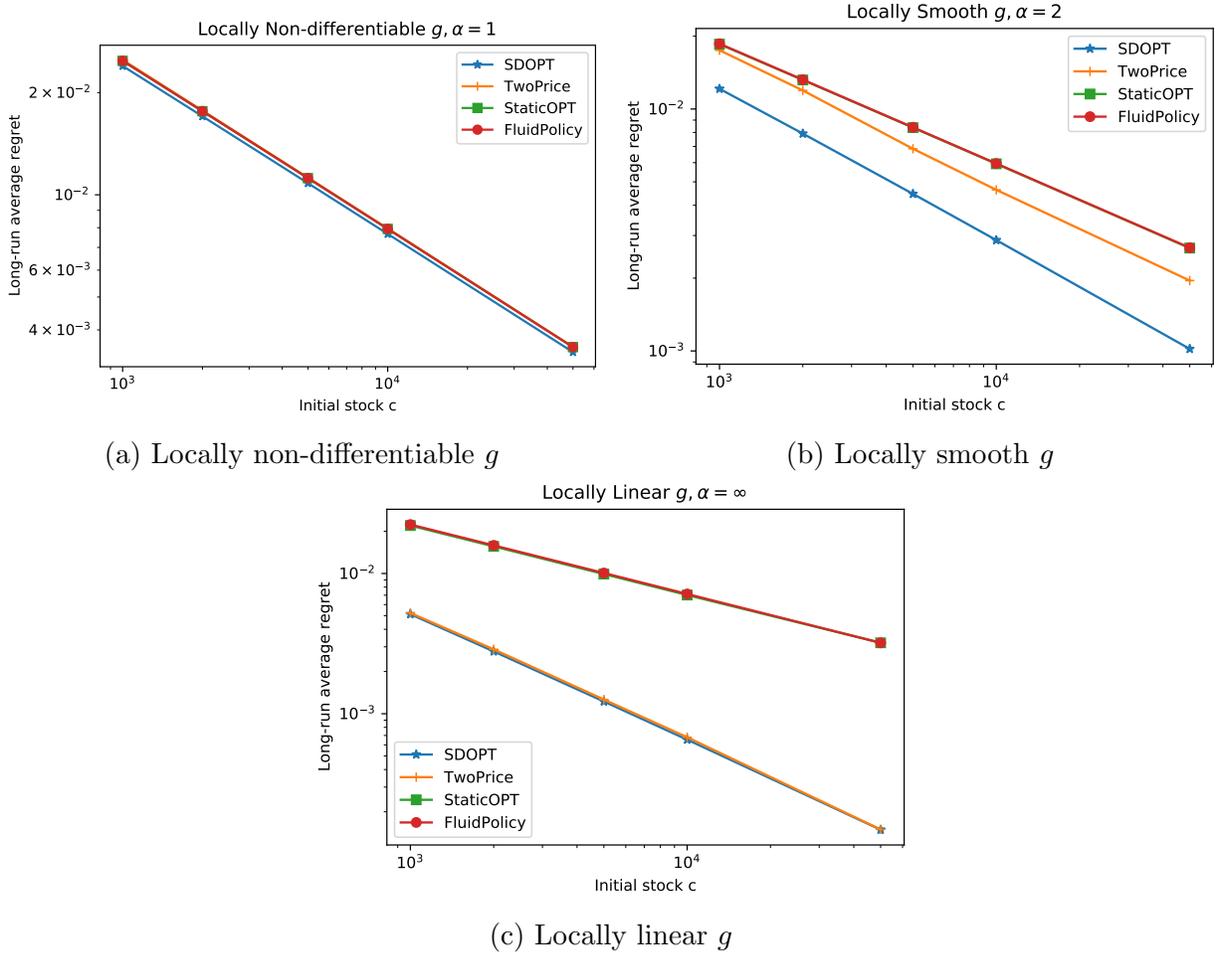

In this part, we numerically illustrate performance loss in some examples. Specifically, for the two-price policy used (denoted as TwoPrice in \Cref{fig: compare convergence rate}), the parameters $x_L,x_H,\tau$ are not necessarily optimized. 

Fix $x^* = c/(\lambda d) = 1/2$. For locally non-differentiable $g$, we consider \Cref{ex: discrete}; for locally smooth $g$, we consider a Uniform$[1,2]$ WTP distribution; 
for locally linear $g$, we consider a pricing problem where WTP is $3,2$ or $1$ w.p. $0.2$, $0.4$, and $0.4$ respectively.

{In \Cref{fig: compare convergence rate}, we plot the performance loss of the four policies with log scaling on both axes: the slopes of the lines are the performance loss scaling. The flatter the line is, the slower the performance loss scales.} We can see for locally non-differentiable $g$ (\Cref{fig: degenerate g}), all three policies show the same performance. For locally smooth $g$ (\Cref{fig: smooth g}), the performance losses of a two-price policy and optimal stock-dependent policy are close and smaller than static policies, but the constants of their performance losses are different. We can find similar results for locally linear $g$ (\Cref{fig: linear g}), but this time the performance loss of a two-price policy and the optimal stock-dependent policy are almost the same. In fact, {on this example}, the optimal stock-dependent policy obtained by solving LP \eqref{LPSD} is a two-price policy. The differences in the performances are due to the fact that the parameters of the two-price policy we used in simulation are not optimized. In (b) and (c), the performance of StaticOPT and FluidPolicy are very close to each other.

\subsection{Small Stock}
\begin{figure}[!htbp]
    \centering
    \begin{subfigure}[b]{0.49\textwidth}
        \begin{tikzpicture}[scale = 0.8]
    \begin{axis}[
        xlabel={No. of Initial Stocks ($c$)},
        ylabel={Proportion of Revenue Gained},
        grid=major,
        legend entries={SDOPT, TwoPriceOPT, StaticOPT, FluidPolicy},
        legend pos=south east, 
        cycle list name=color list,
    ]
    
    \addplot+[mark=triangle*, color=red] coordinates {
        (20, 0.8566663050896396)
        (40, 0.8991422395059463)
        (60, 0.9169055982381861)
        (80, 0.9293900566336813)
        (100, 0.9381217185341018)
    };
    
    \addplot+[mark=o, color=orange] coordinates {
        (20, 0.8526702587308088)
        (40, 0.8944257086574562)
        (60, 0.9117638429473186)
        (80, 0.9238216966845113)
        (100, 0.9323988393843622)
    };
    
    \addplot+[mark=diamond*, color=black] coordinates {
        (20, 0.84390784)
        (40, 0.88586659)
        (60, 0.90494716)
        (80, 0.91692794)
        (100, 0.92534056)
    };
    
    \addplot+[mark=square*, color=blue] coordinates {
        (20, 0.8411080384580303)
        (40, 0.8838440156892446)
        (60, 0.9037331879907949)
        (80, 0.9158812942047665)
        (100, 0.9242995472891383)
    };

    \end{axis}
\end{tikzpicture}
        \caption{Locally non-differentiable $g$}
        \label{fig: degenerate_reward_propotion}
    \end{subfigure}
    \begin{subfigure}[b]{0.49\textwidth}

       \begin{tikzpicture}[scale = 0.8]
    \begin{axis}[
        xlabel={No. of Initial Stocks ($c$)},
        ylabel={Proportion of Revenue Gained},
        grid=major,
        legend entries={SDOPT, TwoPriceOPT, StaticOPT, FluidPolicy},
        legend pos=south east, 
        cycle list name=color list,
    ]
    
    \addplot+[mark=triangle*, color=red] coordinates {
        (20, 0.8727303019582192)
        (40, 0.9088206478140445)
        (60, 0.9256978713347231)
        (80, 0.9357242315259822)
        (100, 0.9440408585227119)
    };
    
    \addplot+[mark=o, color=orange] coordinates {
        (20, 0.8701714025810336)
        (40, 0.9056793541177968)
        (60, 0.9223131168577302)
        (80, 0.9324646986914314)
        (100, 0.9409048853140354)
    };
    
    \addplot+[mark=diamond*, color=black] coordinates {
        (20, 0.86264324)
        (40, 0.89749506)
        (60, 0.91420668)
        (80, 0.92396319)
        (100, 0.93270096)
    };
    
    \addplot+[mark=square*, color=blue] coordinates {
        (20, 0.8411080384580303)
        (40, 0.8838440156892446)
        (60, 0.9037331879907949)
        (80, 0.9158812942047665)
        (100, 0.9242995472891383)
    };

    \end{axis}
\end{tikzpicture}

        \caption{Locally smooth $g$}
        \label{fig: smooth_reward_propotion}
    \end{subfigure}
    \begin{subfigure}[b]{0.49\textwidth}
        \begin{tikzpicture}[scale = 0.8]
    \begin{axis}[
        xlabel={No. of Initial Stocks ($c$)},
        ylabel={Proportion of Revenue Gained},
        grid=major,
        legend entries={SDOPT, TwoPriceOPT, StaticOPT, FluidPolicy},
        legend pos=south east, 
        cycle list name=color list,
    ]
    
    \addplot+[mark=triangle*, color=red] coordinates {
        (20, 0.8679965966958315)
        (40, 0.9091112529075269)
        (60, 0.9288607252853193)
        (80, 0.9397773910306786)
        (100, 0.9479066627724034)
    };
    
    \addplot+[mark=o, color=orange] coordinates {
        (20, 0.8645147008846646)
        (40, 0.9060762301345165)
        (60, 0.9254706186474759)
        (80, 0.9371159233940913)
        (100, 0.9452548776671893)
    };
    
    \addplot+[mark=diamond*, color=black] coordinates {
        (20, 0.85304348)
        (40, 0.89211109)
        (60, 0.91202646)
        (80, 0.9220475)
        (100, 0.93129226)
    };
    
    \addplot+[mark=square*, color=blue] coordinates {
        (20, 0.8411080384580303)
        (40, 0.8838440156892448)
        (60, 0.9037331879907949)
        (80, 0.9158812942047665)
        (100, 0.9242995472891383)
    };

    \end{axis}
\end{tikzpicture}

        \caption{Locally linear $g$}
        \label{fig: linear_reward_propotion}
    \end{subfigure}
    \caption{Comparing the performance of SDOPT, TwoPriceOPT, StaticOPT, and FluidPolicy in terms of their proportion of optimum obtained when $c$ is small.}
    \label{fig: compare reward proportion}
\end{figure}
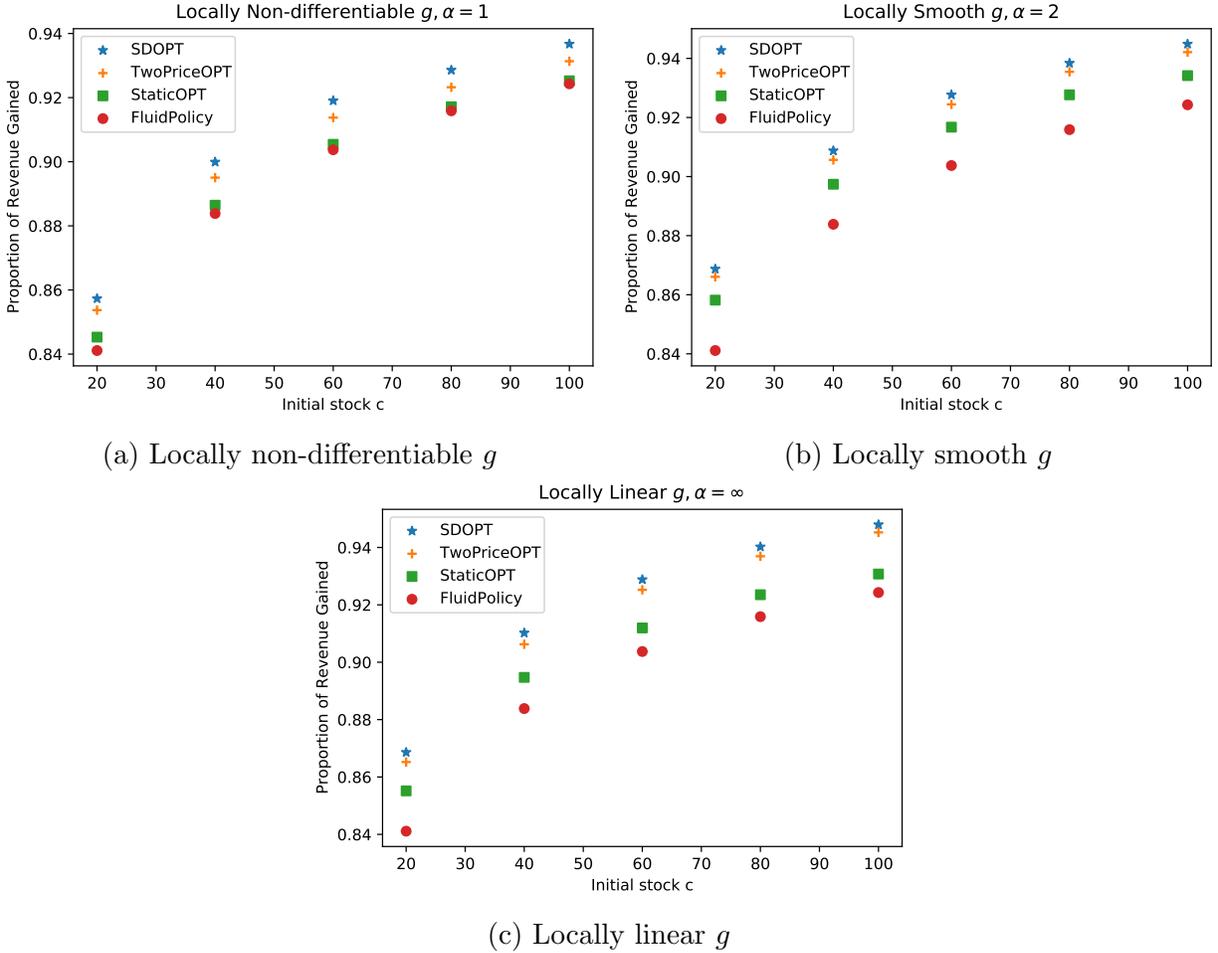

Our theoretical analysis characterizes the performance of static policies and stock-dependent policies in terms of their asymptotic performance loss scaling (compared with $\FLU$) as $c, \lambda$ grows large proportionally. In this part, we provide empirical results to illustrate the performance of the aforementioned policies in terms of their proportion of reward (compared with $\FLU$) when $c$ and $\lambda$ is small. 
In \Cref{fig: compare reward proportion}, we use a two-price policy with optimized parameters, denoted by TwoPriceOPT.\footnote{Optimized using the Sequential Least Squares Programming (\textsf{SLSQP}) algorithm in \textsf{scipy.optimize.minimize}.}
Note that the choice of reward function would affect the proportion of reward gained. Therefore, we randomly generate 100 instances and report the average {of the performance ratios} here. More specifically, for $\alpha = 1$ and $\alpha = \infty$, we generate $6$ and $5$ customer types respectively, where each type has an equal arrival rate and a WTP generated randomly from $1,\ldots,10$ without repetition. (Recall that $x^*=1/2$, which is why having 6 and 5 customer types would lead to $\alpha=1$ and $\alpha=\infty$, respectively.) For $\alpha = 2$, we consider WTP distributed uniformly over $[a,b]$, where $a,b$ are randomly generated parameters from $1,\ldots,10$ that satisfy $a<b$.

We focus on values of $c$ from 20 to 100 and illustrate how the proportion garnered from the benchmark increases with $c$.
As the initial stock $c$ increases, the proportions gained by all policies increase, with marginal gains decreasing. Specifically, when $g$ is locally non-differentiable (\Cref{fig: degenerate_reward_propotion}) and $c=20$, the fluid policy obtains 84.1\% of the reward of $\FLU$, the optimal static policy obtains 84.3\%, while the optimal two-price policy obtains 85.2\% and the optimal stock-dependent policy obtains 85.6\%. When $c=100$, the four policies obtain 92.4\%, 92.5\%, 93.2\%, and 93.7\% respectively.
Even though this is the case where $\alpha=1$ and all four policies have the same asymptotic convergence rate according to our theory, our two-price and stock-dependent policies still noticeably improve empirical performance.
These observations carry through to the $\alpha=2$ and $\alpha=\infty$ cases, in which our policies would also achieve a faster asymptotic convergence rate.

Finally, we point out that most of the gains beyond FluidPolicy come after going from optimizing one price (StaticOPT) to optimizing two prices (TwoPriceOPT).
This is consistent with our theoretical finding that a little dynamicity goes a long way.
We note that the derivation of our stock-dependent and two-price policy classes arose from the desire to analyze tight convergence rates, instead of settling with static policies that already achieved sub-linear performance loss in $c$.
This shows the value of our refined theoretical analysis.

\section{Conclusion}\label{sec: conclusion}
In this paper, we investigate the problem of reusable resource allocation with general reward structures and usage durations in the infinite horizon setting. We propose a class of stock-dependent policies that only consider how many units of stock are available when making decisions, and compare its performance with static policies, benchmarked against $\FLU$, the fluid relaxation problem. Our analysis relies on an insensitivity result for loss networks with state-dependent arrival rates. This property allows us to obtain the steady-state distributions of stock level under any stock-dependent policy in closed form, and more importantly, the corresponding long-run average reward. We can also formulate a convex optimization problem to derive the optimal policy.

We conclude by discussing some extensions that our paper does not handle.
One important extension of our model is to allow state-dependent usage durations---how long customers stay in the system depends on how congested the system is. {We note that} the detailed balance equation \eqref{balance equations} we used for the performance analyses can be extended to this general case \citep{brumelle1978generalization}, so extending the rest of our analysis would be an interesting future direction.

Another important extension {to consider} is to allow for other supply-side dynamics including acquisition of new units and maintenance of service units which are quite common in real-world applications.
Finally, our model assumes the service time distribution and reward function are known to the decision-maker, but in reality, both need to be learned. 
How to do pricing while learning is a practically relevant and theoretically interesting future research direction, and our results provide important insights on this.
Indeed, a two‑price policy is fully specified by just three parameters---the low price, the high price, and the switching threshold---which suggests that a data-driven algorithm can search this low‑dimensional space directly (e.g., via simulation on data samples) rather than fitting a high‑complexity model or policy class. This will be more data-efficient and reduces the risk of overfitting. 
Moreover, to understand how much our policy class improves over the fluid policy, we only need to know the local shape of the reward function at the fluid admission probability $x^*$, instead of the entire function. 
All in all, the implications of our simple class of two-price policies for learning remain to be explored.





\bibliographystyle{informs2014} 
\bibliography{bibliography} 

\begin{thebibliography}{46}
\providecommand{\natexlab}[1]{#1}
\providecommand{\url}[1]{\texttt{#1}}
\providecommand{\urlprefix}{URL }

\bibitem[{Baek \protect\BIBand{} Ma(2022)}]{baek_bifurcating_2022}
Baek J, Ma W (2022) Bifurcating {Constraints} to {Improve} {Approximation} {Ratios} for {Network} {Revenue} {Management} with {Reusable} {Resources}. \emph{Operations Research} opre.2022.2282, ISSN 0030-364X, 1526-5463, \urlprefix\url{http://dx.doi.org/10.1287/opre.2022.2282}.

\bibitem[{Baek \protect\BIBand{} Wang(2023)}]{baek2023leveraging}
Baek J, Wang S (2023) Leveraging reusability: Improved competitive ratio of greedy for reusable resources. \emph{arXiv preprint arXiv:2304.03377} .

\bibitem[{Balseiro et~al.(2021)Balseiro, Besbes, \protect\BIBand{} Pizarro}]{balseiro_survey_2021}
Balseiro S, Besbes O, Pizarro D (2021) Survey of {Dynamic} {Resource} {Constrained} {Reward} {Collection} {Problems}: {Unified} {Model} and {Analysis}. {SSRN} {Scholarly} {Paper} 3963265, Social Science Research Network, Rochester, NY.

\bibitem[{Balseiro et~al.(2022)Balseiro, Besbes, \protect\BIBand{} Castro}]{balseiro2022mechanism}
Balseiro SR, Besbes O, Castro F (2022) Mechanism design under approximate incentive compatibility. \emph{Operations Research} .

\bibitem[{Balseiro et~al.(2019)Balseiro, Brown, \protect\BIBand{} Chen}]{balseiro2019dynamic}
Balseiro SR, Brown DB, Chen C (2019) Dynamic pricing of relocating resources in large networks. \emph{ACM SIGMETRICS Performance Evaluation Review} 47(1):29--30.

\bibitem[{Banerjee et~al.(2022)Banerjee, Freund, \protect\BIBand{} Lykouris}]{banerjee2022pricing}
Banerjee S, Freund D, Lykouris T (2022) Pricing and optimization in shared vehicle systems: An approximation framework. \emph{Operations Research} 70(3):1783--1805.

\bibitem[{Banerjee et~al.(2015)Banerjee, Riquelme, \protect\BIBand{} Johari}]{banerjee2015pricing}
Banerjee S, Riquelme C, Johari R (2015) Pricing in ride-share platforms: A queueing-theoretic approach. \emph{Available at SSRN 2568258} .

\bibitem[{Beck \protect\BIBand{} Teboulle(2003)}]{beck2003mirror}
Beck A, Teboulle M (2003) Mirror descent and nonlinear projected subgradient methods for convex optimization. \emph{Operations Research Letters} 31(3):167--175.

\bibitem[{Benjaafar \protect\BIBand{} Shen(2022)}]{benjaafar2022pricing}
Benjaafar S, Shen X (2022) Pricing in on-demand (and one-way) vehicle sharing networks. \emph{Available at SSRN 3998297} .

\bibitem[{Besbes et~al.(2019)Besbes, Elmachtoub, \protect\BIBand{} Sun}]{besbes2019static}
Besbes O, Elmachtoub AN, Sun Y (2019) Static pricing: Universal guarantees for reusable resources. \emph{Proceedings of the 2019 ACM Conference on Economics and Computation}, 393--394.

\bibitem[{Besbes et~al.(2022)Besbes, Kanoria, \protect\BIBand{} Kumar}]{besbes2022multi}
Besbes O, Kanoria Y, Kumar A (2022) The multi-secretary problem with many types. \emph{Proceedings of the 23rd ACM Conference on Economics and Computation}, 1146--1147.

\bibitem[{Brumelle(1978)}]{brumelle1978generalization}
Brumelle SL (1978) A generalization of erlang's loss system to state dependent arrival and service rates. \emph{Mathematics of Operations Research} 3(1):10--16.

\bibitem[{Bumpensanti \protect\BIBand{} Wang(2018)}]{bumpensanti_re-solving_2018}
Bumpensanti P, Wang H (2018) A {Re}-solving {Heuristic} with {Uniformly} {Bounded} {Loss} for {Network} {Revenue} {Management}. Technical Report arXiv:1802.06192, arXiv.

\bibitem[{Chen et~al.(2023)Chen, Lei, \protect\BIBand{} Jasin}]{chen2023real}
Chen Q, Lei Y, Jasin S (2023) Real-time spatial--intertemporal pricing and relocation in a ride-hailing network: Near-optimal policies and the value of dynamic pricing. \emph{Operations Research} .

\bibitem[{Chen et~al.(2017)Chen, Levi, \protect\BIBand{} Shi}]{chen2017revenue}
Chen Y, Levi R, Shi C (2017) Revenue management of reusable resources with advanced reservations. \emph{Production and Operations Management} 26(5):836--859.

\bibitem[{Delong et~al.(2022)Delong, Farhadi, Niazadeh, \protect\BIBand{} Sivan}]{delong2022online}
Delong S, Farhadi A, Niazadeh R, Sivan B (2022) Online bipartite matching with reusable resources. \emph{Proceedings of the 23rd ACM Conference on Economics and Computation}, 962--963.

\bibitem[{Doan et~al.(2020)Doan, Lei, \protect\BIBand{} Shen}]{doan2020pricing}
Doan XV, Lei X, Shen S (2020) Pricing of reusable resources under ambiguous distributions of demand and service time with emerging applications. \emph{European Journal of Operational Research} 282(1):235--251.

\bibitem[{Elmachtoub \protect\BIBand{} Shi(2023)}]{elmachtoub2023power}
Elmachtoub AN, Shi J (2023) The power of static pricing for reusable resources. \emph{arXiv preprint arXiv:2302.11723} .

\bibitem[{Feng et~al.(2020)Feng, Niazadeh, \protect\BIBand{} Saberi}]{feng_near-optimal_2020}
Feng Y, Niazadeh R, Saberi A (2020) Near-{Optimal} {Bayesian} {Online} {Assortment} of {Reusable} {Resources}. \emph{SSRN Electronic Journal} ISSN 1556-5068, \urlprefix\url{http://dx.doi.org/10.2139/ssrn.3714338}.

\bibitem[{Gallego \protect\BIBand{} Van~Ryzin(1994)}]{gallego1994optimal}
Gallego G, Van~Ryzin G (1994) Optimal dynamic pricing of inventories with stochastic demand over finite horizons. \emph{Management science} 40(8):999--1020.

\bibitem[{Gans \protect\BIBand{} Savin(2007)}]{gans2007pricing}
Gans N, Savin S (2007) Pricing and capacity rationing for rentals with uncertain durations. \emph{Management Science} 53(3):390--407.

\bibitem[{Gong et~al.(2021)Gong, Goyal, Iyengar, Simchi-Levi, Udwani, \protect\BIBand{} Wang}]{gong2021online}
Gong XY, Goyal V, Iyengar GN, Simchi-Levi D, Udwani R, Wang S (2021) Online assortment optimization with reusable resources. \emph{Management Science} .

\bibitem[{Goyal et~al.(2021)Goyal, Iyengar, \protect\BIBand{} Udwani}]{goyal_asymptotically_2021}
Goyal V, Iyengar G, Udwani R (2021) Asymptotically {Optimal} {Competitive} {Ratio} for {Online} {Allocation} of {Reusable} {Resources}. Technical Report arXiv:2002.02430, arXiv, \urlprefix\url{http://arxiv.org/abs/2002.02430}, arXiv:2002.02430 [cs] type: article.

\bibitem[{Hartline(2013)}]{hartline2013mechanism}
Hartline JD (2013) Mechanism design and approximation. \emph{Book draft. October} 122(1).

\bibitem[{Huo \protect\BIBand{} Cheung(2022)}]{huo2022online}
Huo T, Cheung WC (2022) Online reusable resource allocations with multi-class arrivals. \emph{Available at SSRN 4320423} .

\bibitem[{Jasin \protect\BIBand{} Kumar(2012)}]{jasin_re-solving_2012}
Jasin S, Kumar S (2012) A {Re}-{Solving} {Heuristic} with {Bounded} {Revenue} {Loss} for {Network} {Revenue} {Management} with {Customer} {Choice}. \emph{Mathematics of Operations Research} 37(2):313--345.

\bibitem[{Kanoria \protect\BIBand{} Qian(2022)}]{kanoria_blind_2022}
Kanoria Y, Qian P (2022) Blind {Dynamic} {Resource} {Allocation} in {Closed} {Networks} via {Mirror} {Backpressure}. Technical Report arXiv:1903.02764, arXiv, \urlprefix\url{http://dx.doi.org/10.48550/arXiv.1903.02764}, arXiv:1903.02764 [math] type: article.

\bibitem[{Key(1990)}]{key1990optimal}
Key PB (1990) Optimal control and trunk reservation in loss networks. \emph{Probability in the Engineering and Informational Sciences} 4(2):203--242.

\bibitem[{Kim \protect\BIBand{} Randhawa(2018)}]{kim2018value}
Kim J, Randhawa RS (2018) The value of dynamic pricing in large queueing systems. \emph{Operations Research} 66(2):409--425.

\bibitem[{Lei \protect\BIBand{} Jasin(2020)}]{lei2020real}
Lei Y, Jasin S (2020) Real-time dynamic pricing for revenue management with reusable resources, advance reservation, and deterministic service time requirements. \emph{Operations Research} 68(3):676--685.

\bibitem[{Levi \protect\BIBand{} Radovanovi{\'c}(2010)}]{levi2010provably}
Levi R, Radovanovi{\'c} A (2010) Provably near-optimal lp-based policies for revenue management in systems with reusable resources. \emph{Operations Research} 58(2):503--507.

\bibitem[{Miller(1969)}]{miller1969queueing}
Miller BL (1969) A queueing reward system with several customer classes. \emph{Management Science} 16(3):234--245.

\bibitem[{Morrison(2010)}]{morrison2010optimal}
Morrison JA (2010) Optimal trunk reservation for an overloaded link. \emph{Operations research letters} 38(6):499--501.

\bibitem[{Owen \protect\BIBand{} Simchi-Levi(2018)}]{owen2018price}
Owen Z, Simchi-Levi D (2018) Price and assortment optimization for reusable resources. \emph{Available at SSRN 3070625} .

\bibitem[{Papier \protect\BIBand{} Thonemann(2010)}]{papier2010capacity}
Papier F, Thonemann UW (2010) Capacity rationing in stochastic rental systems with advance demand information. \emph{Operations research} 58(2):274--288.

\bibitem[{Paschalidis \protect\BIBand{} Tsitsiklis(2000)}]{paschalidis2000congestion}
Paschalidis IC, Tsitsiklis JN (2000) Congestion-dependent pricing of network services. \emph{IEEE/ACM transactions on networking} 8(2):171--184.

\bibitem[{Rusmevichientong et~al.(2020)Rusmevichientong, Sumida, \protect\BIBand{} Topaloglu}]{rusmevichientong_dynamic_2020}
Rusmevichientong P, Sumida M, Topaloglu H (2020) Dynamic {Assortment} {Optimization} for {Reusable} {Products} with {Random} {Usage} {Durations}. \emph{Management Science} 66(7):2820--2844, ISSN 0025-1909, 1526-5501, \urlprefix\url{http://dx.doi.org/10.1287/mnsc.2019.3346}.

\bibitem[{Rusmevichientong et~al.(2023)Rusmevichientong, Sumida, Topaloglu, \protect\BIBand{} Bai}]{rusmevichientong2023revenue}
Rusmevichientong P, Sumida M, Topaloglu H, Bai Y (2023) Revenue management with heterogeneous resources: Unit resource capacities, advance bookings, and itineraries over time intervals. \emph{Operations Research} .

\bibitem[{Talluri et~al.(2004)Talluri, Van~Ryzin, \protect\BIBand{} Van~Ryzin}]{talluri2004theory}
Talluri KT, Van~Ryzin G, Van~Ryzin G (2004) \emph{The theory and practice of revenue management}, volume~1 (Springer).

\bibitem[{Teicher(1955)}]{teicher1955inequality}
Teicher H (1955) An inequality on poisson probabilities. \emph{The Annals of Mathematical Statistics} 26(1):147--149.

\bibitem[{Tijms(2003)}]{tijms2003first}
Tijms HC (2003) \emph{A first course in stochastic models} (John Wiley and sons).

\bibitem[{Vera \protect\BIBand{} Banerjee(2020)}]{vera_bayesian_2020}
Vera A, Banerjee S (2020) The {Bayesian} {Prophet}: {A} {Low}-{Regret} {Framework} for {Online} {Decision} {Making}. Technical Report arXiv:1901.05028, arXiv, \urlprefix\url{http://arxiv.org/abs/1901.05028}, arXiv:1901.05028 [cs, math] type: article.

\bibitem[{Waserhole \protect\BIBand{} Jost(2016)}]{waserhole2016pricing}
Waserhole A, Jost V (2016) Pricing in vehicle sharing systems: Optimization in queuing networks with product forms. \emph{EURO Journal on Transportation and Logistics} 5(3):293--320.

\bibitem[{Xie et~al.(2022)Xie, Gurvich, \protect\BIBand{} K{\"u}{\c{c}}{\"u}kyavuz}]{xie2022dynamic}
Xie X, Gurvich I, K{\"u}{\c{c}}{\"u}kyavuz S (2022) Dynamic allocation of reusable resources: Logarithmic regret in overloaded networks. \emph{preprint available online} .

\bibitem[{Xu \protect\BIBand{} Li(2013)}]{xu2013dynamic}
Xu H, Li B (2013) Dynamic cloud pricing for revenue maximization. \emph{IEEE Transactions on Cloud Computing} 1(2):158--171.

\bibitem[{Zhang \protect\BIBand{} Cheung(2022)}]{zhang2022online}
Zhang X, Cheung WC (2022) Online resource allocation for reusable resources. \emph{arXiv preprint arXiv:2212.02855} .

\end{thebibliography}

\ECSwitch
\ECDisclaimer
\ECHead{E-Companion}

\begin{APPENDICES}
\crefalias{section}{appendix}
\crefalias{subsection}{appendix}

\section{Examples of Reward Functions}\label{apx: extension on g}
The performance loss guarantees hold for any reward function $g$ that is concave and non-decreasing with $g(0)=0$. We have presented examples of reward functions when the objective is revenue maximization or reward maximization, here we present some examples of different payment collection schemes and cost structures.

\begin{example}[One-Time Fee Payment]
If customers pay a fixed price $p$ when admitted into the system, then the reward rate with admission probability $x$ is $g(x) = px$.
\end{example}

\begin{example}[Rate Payment]
If customers pay a fixed price $u$ per time in the system, then the reward rate with admission probability $x$ is $g(x) = dux$ as the expected usage duration is $d$. 
\end{example}

If the customers pay a one-time fixed price $p$ when admitted into the system and a fixed price $u$ per time in the system, then the reward function is $g(x) = px + dux$.

\begin{example}[Adding Operating Cost]
If admitting a new customer into the system incurs a fixed cost of $c_0>0$, then we can modify the original reward function $g(x)$ to be the upper concave envelope of $g(x)- x c_0$. The modified function, denoted as $g_{c_0}(x)$, is still concave and non-decreasing with $g_{c_0}(0) = 0$. Moreover, the reward functions $g_{c_0}(x)$ and $g(x)$ satisfy \Cref{assmp: condition for upper bound} and \Cref{assmp: condition for lower bound} with the same local curvature parameters $\alpha$, thus the same performance guarantees apply.
\end{example}

{
\section{Example Illustrating the Difference Between Assumptions \ref{assmp: condition for upper bound} and \ref{assmp: condition for lower bound}}\label{apx: example of mismatch alpha}

Local curvatures of $g$ can be different to the left and right of the fluid admission probability $c/(\lambda d)$. Consider $g(x) = x$ for $x\in [0,0.5]$,  $g(x) = 2x-x^2-0.25$ for $x \in [0.5,1]$, both with $c/(\lambda d) = 0.5$, then $g$ satisfies Assumption 1 with $\alpha = 1$ and Assumption 2 with $\alpha = \infty$. See Figure \ref{fig: illus mismatch A12} for illustration. 
    
\begin{figure}
    \centering
\begin{tikzpicture}
    \begin{axis}[
        axis lines = middle,
        xlabel = $\mu$,
        ylabel = $g(\mu)$,
        ymin=0, ymax=1,
        xmin=0, xmax=1,
        xtick={0.5},
        xticklabels={$c/d$},
        ytick=\empty,
        width=8cm,
        height=6.4cm,
        legend style={at={(1.05,0.2)}, anchor=south east}, 
        every axis x label/.style={at={(current axis.right of origin)},anchor=west},
        every axis y label/.style={at={(current axis.north west)},above=2mm},
        ]
        \addplot[domain=0:0.5, samples=100, blue, thick, forget plot] {x};
        \addplot[domain=0.5:1, samples=100, blue, thick] {2*x-x*x-0.25};
        \addlegendentry{${g}$};

        \addplot[domain=0.3:0.5, samples=100, dotted, red, thick, forget plot] {2*x-0.5};
        \addplot[domain=0.5:0.8, samples=100, dotted, red, thick] {0.5*x + 0.25};
        \addlegendentry{$\alpha = 1$ (A1)};
    
        \addplot[domain=0.2:0.8, samples=100, dashed, orange, thick] {x};
        \addlegendentry{$\alpha = \infty$ (A2)};
        \draw[dashed] (axis cs:0.5,0) -- (axis cs:0.5,0.5);
    \end{axis}
\end{tikzpicture}
    \caption{Illustration of a $g(\cdot)$ that has different $\alpha$ for Assumption 1 and 2.}
    \label{fig: illus mismatch A12}
\end{figure}

}

\section{Proof of \Cref{lem: FLU benchmark}}\label{apx: FLU benchmark}
\proof{Proof.}
{
Let $x_t^\sigma$ denote the admission probability under the online policy $\sigma$ at time $t$, for all $t\in [0,T]$, which is 0 if the system is out of stock at time $t$. Define
\[x^\sigma := \liminf_{T\to\infty}\frac{1}{T}\int_0^Tx_t^\sigma \, dt.\]

Then, by Little's Law and the feasibility of the online policy $\sigma$, $\lambda \cdot x^\sigma \cdot  d \leq c$.

The long-run average reward obtained by policy $\sigma$ is 
\[\liminf_{T\to\infty}\frac{1}{T}\int_{t=0}^T \lambda g(x^\sigma_t)\,dt \stackrel{(a)}{\leq}  \liminf_{T\to\infty}\lambda g\left(\frac{1}{T}\int_0^T x^\sigma_t \,dt\right) \stackrel{(b)}{=} \lambda g\left(\liminf_{T\to\infty}\frac{1}{T}\int_0^T  x^\sigma_t \,dt\right) = \lambda g(x^\sigma) \leq g(x^*),\]
where $(a)$ follows from Jensen's inequality; $(b)$ follows from the continuity of $g$.
}
\Halmos\endproof

\section{Analysis of the Optimization Problem \eqref{OPTSD}}

\subsection{Proof of \Cref{pro: optsd convex}}\label{apx: optsd convex}
\proof{Proof.}
Define $f_j(y) = \lambda g(({c-j+1})y/(\lambda {d}))$. As $g$ is concave, $f_{j}$ is also concave.

For any $\alpha\in [0,1]$ and $\bm{z}, \bm{y}\in \mathbb{R}^{c+1}_+$ such that $\sum_{j=0}^c z_{j}=\sum_{j=0}^c y_{j}=1$, by \Cref{pro: optimal pis positive}, we can assume $z_j, y_j >0, \forall j \in [0:c]$ without loss. Let $\Pi(\bm{\pi}) = \sum_{j=1}^c \pi_j f_j({\pi_{j-1}}/{\pi_j})$ denote the objective function:
\begin{align*}
\Pi(\alpha\bm{z}+(1-\alpha)\bm{y})& = \sum_{j=1}^c(\alpha z_{j}
 +(1-\alpha) y_{j})f_j\left(\frac{\alpha z_{j-1}
 +(1-\alpha) y_{j-1}}{\alpha z_{j}
 +(1-\alpha) y_{j}}\right)\\
 & = \sum_{j=1}^c(\alpha z_{j}
 +(1-\alpha) y_{j})f_j\left(\frac{\alpha z_{j}}{\alpha z_{j}
 +(1-\alpha) y_{j}}\frac{z_{j-1}}{z_{j}}+\frac{(1-\alpha)y_{j}}{\alpha z_{j}
 +(1-\alpha) y_{j}}\frac{y_{j-1}}{y_{j}}\right)\\
 & \stackrel{(a)}{\geq} \sum_{j=1}^c(\alpha z_{j}
 +(1-\alpha) y_{j})\left(\frac{\alpha z_{j}}{\alpha z_{j}
 +(1-\alpha) y_{j}}f_j\left(\frac{z_{j-1}}{z_{j}}\right)+\frac{(1-\alpha)y_{j}}{\alpha z_{j}
 +(1-\alpha) y_{j}} f_j\left(\frac{y_{j-1}}{y_{j}}\right)\right)\\
 & = \alpha\sum_{j=1}^c z_{j} f_j\left(\frac{z_{j-1}}{z_{j}}\right)+(1-\alpha)\sum_{j=1}^c y_{j} f_j\left(\frac{y_{j-1}}{y_{j}}\right) \\
 & = \alpha \Pi(\bm{z}) +(1-\alpha)\Pi(\bm{y}),
\end{align*}
where $(a)$ follows from the concavity of $f_j$.\Halmos \endproof

\begin{proposition}\label{pro: optimal pis positive}
In the optimal solution to \eqref{OPTSD}, $\pi_j^* >0$.
\end{proposition}
\proof{Proof.}
{
First note that if $\pi_j=0$ in a feasible solution, by \Cref{pro: balance equations}, we must have $\pi_{j-1}=0$. Following this argument, we would have $\pi_\ell=0$ for all $\ell\leq j$. 

We prove $\pi_j^* >0$ for $j = 0,\ldots,c$ by contradiction. Suppose not, denote $k=\arg\max\{j: \pi^*_j = 0\}$, $0\leq k<c$. Here $k \neq c$ as $\sum_{j=0}^c \pi_j = 1$. 
Then $x^*_{k+1}=0, \pi^*_{k+1} \neq 0$. If $k = c-1$, the objective value of the original solution is zero, whereas letting $x_c = 1$ yields a nonnegative objective value, which contradicts the optimality of the original solution. We proceed to show that if $k \leq c-2$, there exists a solution $\bm{\pi}, \bm{x}$ that improves the objective value of $\bm{\pi}^*, \bm{x}^*$, where $\pi_k \neq 0, x_{k+1} \neq 0$. More specifically, the new solution needs to satisfy the following constraints:

\begin{align}
& \pi_k + \pi_{k+1} = \pi_{k+1}^*, \label{eq:EC1_1} \\
& \pi_{k+1}x_{k+1} = \frac{c-k}{\lambda d} \pi_k, \label{eq:EC1_2} \\
& \pi_{k+2}^*x_{k+2} = \frac{c-k-1}{\lambda d}\pi_{k+1}, \label{eq:EC1_3} \\
& \pi_j = \pi^*_{j}, \quad j \geq k+2,\notag \\
& x_j = x_j^*, \quad j \geq k+3, \notag \\
& \pi_{k+1} g(x_{k+1}) + \pi_{k+2}^* g(x_{k+2}) - \pi_{k+2}^* g(x_{k+2}^*) > 0, \label{eq:EC1_6}
\end{align}
where the first five equations imply $\bm{\pi}$ and $\bm{x}$ is feasible and the last inequality implies the newly constructed solution improves over $\bm{\pi}^*, \bm{x}^*$. 

Combining constraint \eqref{eq:EC1_3} with the last inequality \eqref{eq:EC1_6}, the inequality can be written as
\begin{align}
    \frac{x_{k+2}}{(c-k-1)/(\lambda d)}g(x_{k+1}) > g(x_{k+2}^*) - g(x_{k+2}). \label{eq: last ineq reform}
\end{align}

We next show that there exists $x_{k+1}<x_{k+2}$ such that inequality \eqref{eq: last ineq reform} holds. 

\begin{align*}
    \frac{x_{k+2}}{(c-k-1)/(\lambda d)}g(x_{k+1}) & \stackrel{(a)}{\geq} \frac{x_{k+1}}{(c-k-1)/d}g(x_{k+2}) >\frac{x_{k+1}}{(c-k)/(\lambda d)}g(x_{k+2}) \\
    &  \stackrel{(b)}{=} \left( \frac{x_{k+2}^*}{x_{k+2}}-1\right)g(x_{k+2}) \stackrel{(c)}{\geq} g(x_{k+2}^*)-g(x_{k+2}),
\end{align*}
where $(a)$ follows from \Cref{lem: f(y)/y} as $g$ is concave with $g(0)=0$ and $0<x_{k+1}<x_{k+2}$; $(b)$ holds as: combining \eqref{eq:EC1_1} and \eqref{eq:EC1_2} we have 
\[\pi_{k+1} = \frac{(c-k)/(\lambda d)}{(c-k)/(\lambda d) + x_{k+1}} \pi_{k+1}^*,\]
thus
\[x_{k+2} = \frac{c-k-1}{\lambda d} \frac{\pi_{k+1}}{\pi_{k+2}^*}  = \frac{c-k-1}{d} \frac{\pi_{k+1}^*}{\pi_{k+2}^*} \frac{(c-k)/(\lambda d)}{(c-k)/(\lambda d) + x_{k+1}} = x_{k+2}^* \frac{(c-k)/(\lambda d)}{(c-k)/(\lambda d) + x_{k+1}} < x_{k+2}^*,\]
which implies 
\[\frac{x_{k+1}}{(c-k)/(\lambda d)} = \frac{x_{k+2}^*}{x_{k+2}}-1; \]
$(c)$ follows from \Cref{lem: f(y)/y} and $0< x_{k+2} < x_{k+2}^*$ again.

Given $x_{k+1}$ and $x_{k+2}$ that satisfy inequality \eqref{eq: last ineq reform}, $\pi_k$ and $\pi_{k+1}$ are uniquely defined. Therefore, we show that we can construct a new solution that improves over $\bm{x}^*, \bm{\pi}^*$, which contradicts the optimality of the original solution. \Halmos \endproof

\begin{lemma}\label{lem: f(y)/y}
Suppose $f: [0,1]\mapsto \mathbb{R}^+$ is concave with $f(0)=0$, then $f(y)/y$ is non-increasing in $y \in (0,1]$.
\end{lemma}
\proof{Proof.}
By concavity, for $\alpha\in [0,1]$, $f(\alpha y + (1-\alpha)0) \geq \alpha f(y) + (1-\alpha)f(0) = \alpha f(y)$. 
Thus for $0< y'<y$, let $\alpha = y'/y$, we have $f(y') \geq y' f(y)/y$.
\Halmos\endproof
}

\subsection{Proof of \Cref{thm: optsd to lpsd}}\label{apx: optsd to lpsd}
\proof{Proof.}
Suppose $g(x) = \min_{k\in[m]}\{a_k+b_k x\}$, \eqref{OPTSD} can be written as 
\begin{align}
    \max_{\bm{\pi}\geq \bm{0}} & \sum_{j=1}^c\pi_j g(x_j) \notag\\
   \text{s.t. }  &\pi_jx_j = \pi_{j-1}\frac{c-j+1}{\lambda d}, \forall j \in [c], \notag\\
   & 0 \leq x_j \leq 1, \forall j \in [c], \notag\\
   & \sum_{j=0}^c\pi_j = 1, \notag\\
   & g(x_j) \leq a_k+b_kx_j, \forall k \in [m],  \forall j \in [c]. \label{OPTSD (finite)}
\end{align}

For any feasible solution $(\{\pi_j\}_{j=0}^c, \{x_j\}_{j=1}^c, \{g(x_j)\}_{j=0}^c)$ to \eqref{OPTSD (finite)}, we must have 
    $$\pi_jg(x_j) \leq a_k\pi_j+b_kx_j\pi_j = a_k\pi_j+b_k\frac{c-j+1}{\lambda  d}\pi_{j-1}.$$
    Therefore, define $g_j = \pi_j g(x_j), \forall  j \in [c]$, then 
    $(\{\pi_j\}_{j=0}^c, \{g_j\}_{j=1}^c)$ is feasible to \eqref{LPSD} with the same objective value. 

For any feasible solution $(\{\pi_j\}_{j=0}^c, \{g_j\}_{j=1}^c)$ to \eqref{LPSD}, we show how to construct a feasible solution to \eqref{OPTSD (finite)} with the same objective value. As $\pi_j \geq ({c-j+1})\pi_{j-1}/(\lambda d), \forall j =1,\ldots, c$, if $\pi_k >0$, then $\pi_j >0$ for all $j >k$. Define $k = \max\{j: \pi_j = 0\}$, then $k<c$ as $\sum_{j=0}^c \pi_j = 1$ and $g_j \leq 0$ for $j =1,\ldots, k$. Let $x_j = 0, g(x_j) = 0, \forall j = 1,\ldots, k$ and 
    $$x_j = \frac{c-j+1}{\lambda  d}\frac{\pi_{j-1}}{\pi_j}, g(x_j) = \frac{g_j}{\pi_j}, \forall j >k,$$
    then $x_j \leq 1$, $g(x_j) \leq a_k + b_kx_j$, for all $j,k$.  
\Halmos \endproof

\section{Insufficiency of Static Policies}\label{apx: static policy}

\subsection{Proof of \Cref{pro: failure of fluid pricing}}\label{apx: failure of fluid pricing}
\proof{Proof.}
The performance loss of {the} fluid policy is
\begin{align} \label{eqn:9803}
\lambda g(x^*)-\sum_{j=1}^c\pi_j \lambda g(x^*) = \pi_0 \lambda g(x^*).
\end{align}
Thus the performance loss scaling solely depends on $\pi_0$, the probability that no units of stock are available. 
We derive:
\[\frac{1}{\pi_0} \stackrel{(a)}{=} \sum_{i=0}^c \frac{c!}{c^i(c-i)!} \stackrel{(b)}{\geq} \sum_{i=0}^c \frac{\sqrt{2\pi c}(c/e)^{c}}{c^i(c-i)!}e^{\frac{1}{12c+1}} \geq \sqrt{2\pi c} \sum_{i=0}^c \frac{c^{c-i}e^{-c}}{(c-i)!} \stackrel{(c)}{=}  \sqrt{2\pi c} \sum_{i=0}^c \frac{c^{i}e^{-c}}{i!} \stackrel{(d)}{\geq} \frac{1}{e}\sqrt{2\pi c},\]
\[\frac{1}{\pi_0} \stackrel{(a)}{=} \sum_{i=0}^c \frac{c!}{c^i(c-i)!} \stackrel{(b)}{\leq} \sum_{i=0}^c \frac{\sqrt{2\pi c}(c/e)^{c}}{c^i(c-i)!}e^{\frac{1}{12c}} \leq \sqrt{2\pi c} \sum_{i=0}^c \frac{c^{c-i}e^{-c}}{(c-i)!}\cdot 2 \stackrel{(c)}{=}  2\sqrt{2\pi c} \sum_{i=0}^c \frac{c^{i}e^{-c}}{i!} \stackrel{(d)}{\leq} 2\sqrt{2\pi c},\]
where $(a)$ follows from \Cref{pro: balance equations}; $(b)$ follows from Stirling's Formula; 
$(c)$ follows from reindexing $i$; and $(d)$ follows from the fact that $(c^{i}e^{-c})/{i!}$ can be seen as $\pr(X=i)$ where $X$ is a Poisson random variable with parameter $c$ (which implies that $\sum_{i=0}^c (c^{i}e^{-c})/{i!} = \mathbb{P}(X \leq c) \leq 1$; and the lower bound follows from \cite{teicher1955inequality} who shows that if $X$ is a Poisson random variable with mean $x \geq 0$, then $\mathbb{P}( X \leq \lfloor x \rfloor) > 1/e$).  {This shows that $\pi_0=\Theta(1/\sqrt{c})$, which completes the proof of \Cref{pro: failure of fluid pricing} after substituting into~\eqref{eqn:9803} and recalling that $\lambda=\Theta(c)$.}
\Halmos \endproof 

\subsection{Two Special Cases}\label{apx: 2special cases}
We first present two useful lemmas in bounding the factorial terms. 
\begin{lemma}\label{lem: bound product term}
 For $j\in [c]$, the following holds
 $$1 - \frac{j^2}{2c} \leq \frac{c!}{c^j(c-j)!} \leq 1.$$
\end{lemma}
\proof{Proof.} By explicitly writing out the factorial, we obtain
\[\frac{c!}{c^j(c-j)!} = \frac{c}{c}\frac{c-1}{c}\cdots \frac{c-j+1}{c} = 1\cdot\left(1-\frac{1}{c}\right)\cdots \left(1-\frac{j-1}{c}\right) \geq 1- \sum_{i=1}^{j-1} \frac{i}{c} = 1 - \frac{j(j-1)}{2c} \geq 1 - \frac{j^2}{2c}.\Halmos\]
\endproof

\begin{lemma}\label{lem: sum j2rj}
For $n\in \mathbb{Z}^+, z \in (0,1)$, $\sum_{j=0}^n z^j j^2 \leq \frac{z(1+z)}{(1-z)^3}$.
\end{lemma}
\proof{Proof.}
We first show 
\begin{align}
    \sum_{j=0}^n z^j j^2 = \frac{z(z+1-((1-z)^2n^2+2(1-z)n+z+1)z^n)}{(1-z)^3}. \label{eq: sum zjj2}
\end{align}
Then, the result follows. 

Prove by induction. For $n=1$, the left-hand side is $z$, while the right-hand side is 
$$\frac{z(z+1-((1-z)^2n^2+2(1-z)n+z+1)z^n)}{(1-z)^3} = z.$$ 

Suppose we have \eqref{eq: sum zjj2} holds for $n$, then 
\begin{align*}
    \sum_{j=0}^{n+1} z^j j^2 & = \frac{z(z+1-((1-z)^2n^2+2(1-z)n+z+1)z^n)}{(1-z)^3} + z^{n+1}(n+1)^2 \\
    & = \frac{z(z+1-((1-z)^2(n+1)^2+2(1-z)((n+1)+z+1)z^{n+1})}{(1-z)^3}.\Halmos
\end{align*}
\endproof

\subsubsection{Proof of \Cref{pro: scarce supply helps}}

\begin{figure}
    \centering
        \begin{tikzpicture}
    \begin{axis}[
        axis lines = middle,
        xlabel = $x$,
        ylabel = $g(x)$,
        ymin=0, ymax=1,
        xmin=0, xmax = 1.1,
        xtick={0.001, 0.5, 0.7, 1},
        xticklabels={0,$x^*$, $x^*+\varepsilon$,1},
        ytick=\empty,
        width=8cm,
        height=6.4cm,
        every axis x label/.style={at={(current axis.right of origin)},anchor=west},
        every axis y label/.style={at={(current axis.north west)},above=2mm},
        ]
        \addplot[domain=0:0.7, samples=100, blue, thick] {x};
        \addplot[domain=0.7:1, samples=100, blue, thick] {0.7};
        \draw[densely dotted] (axis cs:0.5,0) -- (axis cs:0.5,0.5);
        \draw[densely dotted] (axis cs:0.7,0) -- (axis cs:0.7,0.7);
    \end{axis}
\end{tikzpicture}
    \caption{Illustration of Extremely Scarce Supply.}
    \label{fig: scarce supply}
\end{figure}

\proof{Proof.}
See \Cref{fig: scarce supply} for illustration. If we set $x^*+\varepsilon$ (i.e., $c/(\lambda d) + \varepsilon$) to be the static admission probability, we have: 
\[ 1 =  \sum_{j=0}^c \pi_j \stackrel{(a)}{=} \pi_0 \sum_{j=0}^c \frac{c!}{(c+\varepsilon \lambda d)^j(c-j)!} \stackrel{(b)}{\geq} \pi_0 \sum_{j=0}^{c}z^j (1-\frac{j^2}{2c})\stackrel{(c)}{\geq} \pi_0 \left(\frac{1-z^{c}}{1-z} -\frac{1}{2c} \frac{z(1+z)}{(1-z)^3}\right),\]
where $(a)$ follows from \Cref{pro: balance equations}:
\[\pi_j \lambda \left(\frac{c}{\lambda d} + \varepsilon\right) = \pi_{j-1} \frac{c-j-1}{d} \implies \pi_j = \frac{c}{c+\varepsilon \lambda d} \cdot \frac{c-1}{c+\varepsilon \lambda d} \cdot \cdots \cdot  \frac{c-j+1}{c+\varepsilon \lambda d} = \frac{c!}{(c+\varepsilon \lambda d)^j(c-j)!};\]
$(b)$ follows from Lemma \Cref{lem: bound product term} and let $z = c/(c+\varepsilon \lambda d) < 1$; $(c)$ follows from \Cref{lem: sum j2rj}.

Thus 
$$\pi_0 \leq \left(\frac{1-z^{c}}{1-z} -\frac{1}{2c} \frac{z(1+z)}{(1-z)^3}\right)^{-1} = 1-z+\mathcal{O}(c^{-1}).$$ 

 As $g(x) = rx$ for $x\in [0,x^*+\varepsilon]$, the performance loss of this static policy is 
\[\lambda g(x^*) - (1-\pi_0)  \lambda g(x^*+\varepsilon) = \lambda r x^* - \lambda r(x^*+\varepsilon)  + (1-z +\mathcal{O}(c^{-1}))\lambda r(x^*+\eps)  =  \mathcal{O}(1). \Halmos 
\]
\endproof

\subsubsection{Proof of \Cref{pro: excess supply helps}}

\begin{figure}
    \centering
        \begin{tikzpicture}
    \begin{axis}[
        axis lines = middle,
        xlabel = $x$,
        ylabel = $g(x)$,
        ymin=0, ymax=1,
        xmin=0, xmax = 1.1,
        xtick={0.001, 0.5, 0.7, 1},
        xticklabels={0,$q$, $x^*$,1},
        ytick=\empty,
        width=8cm,
        height=6.4cm,
        every axis x label/.style={at={(current axis.right of origin)},anchor=west},
        every axis y label/.style={at={(current axis.north west)},above=2mm},
        ]
        \addplot[domain=0:0.5, samples=100, blue, thick] {x*(2-x)};
        \addplot[domain=0.5:1, samples=100, blue, thick] {0.75};
        \draw[densely dotted] (axis cs:0.5,0) -- (axis cs:0.5,0.75);
        \draw[densely dotted] (axis cs:0.7,0) -- (axis cs:0.7,0.75);
    \end{axis}
\end{tikzpicture}
    \caption{Illustration of Excess Supply.}
    \label{fig: excess supply}
\end{figure}
\proof{Proof.}
Define $q = \min\{x: 0\in \partial g(x)\}$, then $q<x^*$ and $g(x^*) = g(q)$. See \Cref{fig: excess supply} for illustration. 
Define $\pi_0(q)$ to be the stationary probability of running out of capacity under a static policy that sets $x_j = q, \forall j \in [c]$.

\begin{align*}
    \frac{1}{\pi_0(q)} & =\sum_{j=0}^c \frac{c!}{(\lambda dq)^j(c-j)!} \stackrel{(a)}{\geq} \sum_{j=0}^c \frac{\sqrt{2\pi c}(c/e)^c}{(\lambda dq)^j(c-j)!}e^{\frac{1}{12c}}
    \stackrel{(b)}{\geq} \sqrt{2\pi c }\sum_{j=0}^c \frac{c^{c-j}e^{-c}}{(\lambda dq/c)^j(c-j)!}   
    \stackrel{(c)}{=} \sqrt{2\pi c }\sum_{j=0}^c \frac{c^{j}e^{-c}}{(\lambda dq/c)^{c-j}j!}\\
    & \geq \sqrt{2\pi c }\sum_{j=0}^{\lfloor c/2\rfloor} \frac{c^{j}e^{-c}}{(\lambda dq/c)^{c-j}j!}    
    \stackrel{(d)}{\geq} \sqrt{2\pi c }\sum_{j=0}^{\lfloor c/2\rfloor} \frac{c^{j}e^{-c}}{(\lambda dq/c)^{j}j!}     
    = \sqrt{2\pi c }e^{-\frac{\lambda dq-c}{\lambda dq/c}}\sum_{j=0}^{\lfloor c/2\rfloor} \frac{\frac{c}{\lambda dq/c}^{j}e^{-\frac{c}{\lambda dq/c}}}{j!},
\end{align*}
where $(a)$ follows from Stirling's Formula; $(b)$ follows from $e^{1/12c} \geq 1$ for all $c\geq 1$; $(c)$ follows from re-indexing $j$; $(d)$ follows as $\lambda dq/c <1$, thus $(\lambda dq/c)^{c-j} < (\lambda dq/c)^j$ for $j \in [0, \lfloor {c}/{2} \rfloor]$.

As $$\sum_{j=0}^{\lfloor c/2\rfloor} \frac{\frac{c}{\lambda dq/c}^{j}\exp({-\frac{c}{\lambda dq/c}})}{j!}$$ can be seen as the probability of $\mathbb{P}(X = j)$ where $X$ is a Poisson random variable with parameter $c/(\lambda dq/c)$, thus the sum is bounded by 1. Therefore we have the probability of running out of stock is
\begin{align}
    \pi_0(q) = \mathcal{O}(e^{-c}). \label{bound pi0 excess supply}
\end{align}

Therefore, the performance loss is
\begin{align*}
   \lambda g(x^*)- \lambda g(x)(1-\pi_0(q)) \stackrel{(a)}{=} \lambda g(x)\pi_0(q) \stackrel{(b)}{=} \mathcal{O}(c/e^{c}),
\end{align*}
where $(a)$ follows as $g(x^*) = g(q)$; $(b)$ follows from \eqref{bound pi0 excess supply}.
\Halmos \endproof

\subsection{Proof of Proposition \ref{pro: failure of static pricing}}\label{apx: lower bound statis pricing}
\proof{Proof.}
Excluding the two special cases discussed in \Cref{apx: 2special cases} implies that 
there exist $r,b >0$ such that $g(x) \leq rx+b, \forall x \in [0,1]$ and $g(x^*) = rx^*+b$.

The proof relies on showing a feasible solution to the dual problem of the optimization problem yields an objective value of $\lambda g(x^*) - \Theta(\sqrt{c})$. 

First, by \Cref{thm: optsd to lpsd} we can analyze the following LP to obtain an upper bound on the performance of the optimal static policy---the last constraint restricts the policy to set admission probability $q$ for any stock level in $[c]$, and the constraints $\pi_j \geq \frac{c-j+1}{\lambda  d}\pi_{j-1}, \forall j \in [c]$ in \eqref{LPSD} is relaxed.

\begin{align*}
    \max_{\bm{\pi}\geq \bm{0}} & \sum_{j=1}^c \lambda g_j\\
   \text{s.t. }  & \sum_{j=0}^c\pi_j = 1, & &&\text{(dual variable $\zeta$)}\\
   & g_j \leq r\frac{c-j+1}{\lambda d}\pi_{j-1} + b \pi_j, &\forall j \in [c], \\
   & \pi_{j}q = \frac{c-j+1}{\lambda d}\pi_{j-1}, & \forall j \in [c]. &&\text{(dual variable $\beta_j$)}
\end{align*}

By moving the second set of the constraints to the objective, i.e., 
\[\max_{\bm{\pi}\geq \bm{0}} \sum_{j=1}^c \lambda g_j =  \max_{\bm{\pi}\geq \bm{0}} \sum_{j=1}^c \lambda r({c-j+1})/(\lambda {d}) \pi_{j-1} + \lambda b \pi_j,\]
we can write the dual of the above LP as:
\begin{align*}
\min\quad &\zeta 
\\ \text{s.t.}\quad &\zeta \ge \lambda r \frac{c-j}{\lambda d}+\lambda b -\beta_j q + \beta_{j+1} \frac{c-j}{\lambda  d}, &\forall j \in [c-1],\\ 
& \zeta \ge \lambda r\frac{c}{\lambda  d}+\beta_1 \frac{c}{\lambda d}, \\
&\zeta \ge \lambda b-\beta_cq.
\end{align*}

Second, we present different feasible solutions depending on whether $q \geq x^*$.

\subsubsection{First case: $q \geq x^*$.\ } 
Consider the following solution
\begin{align}\label{eq: lower bound static pricing q big}
    \beta_j = \lambda \frac{b}{x^*} \max\left\{1-\frac{j}{\sqrt{c}}, 0\right\}, \forall c \in [j], \zeta = \lambda g(x^*) - \min\{b, r x^*/2\} \frac{\lambda}{\sqrt{c}}.
\end{align}

\noindent\textbf{Constraint $\zeta \ge \lambda r \frac{c-j}{\lambda d}+\lambda b -\beta_j q + \beta_{j+1} \frac{c-j}{\lambda  d}$}

The right-hand-side of the constraint can be written as 

\begin{align*}
    \lambda r \frac{c-j}{\lambda d}+\lambda b -\beta_j q + \beta_{j+1} \frac{c-j}{\lambda  d} & = \lambda r \frac{c}{\lambda d} - \lambda r \frac{j}{\lambda d}+ \lambda b - \beta_j q + \beta_{j+1} \frac{c}{\lambda d} - \beta_{j+1} \frac{j}{\lambda d}\\
    & \stackrel{(a)}{=}  \lambda g(x^*) - \beta_j q + \beta_{j+1} \frac{c}{\lambda d} - (\lambda r + \beta_{j+1}) \frac{j}{\lambda d}\\
    & \stackrel{(b)}{\leq} \lambda g(x^*) - (\beta_j - \beta_{j+1}) x^* - (\lambda r + \beta_{j+1}) \frac{j}{\lambda d},
\end{align*}
where $(a)$ holds as $\lambda g(x^*)  =\lambda rx^* + \lambda b, x^* = c/(\lambda d)$; $(b)$ holds as $q \geq x^* = c/(\lambda d), \beta \geq 0$.

\begin{itemize}
    \item When $j < \lfloor \sqrt{c} \rfloor$, $\beta_j - \beta_{j+1} = (\lambda b/x^*)/\sqrt{c}$, thus 
    \[\lambda g(x^*) - (\beta_j - \beta_{j+1}) x^* - (\lambda r + \beta_{j+1}) \frac{j}{\lambda d} \leq \lambda g(x^*) - (\beta_j - \beta_{j+1}) x^* = \lambda g(x^*) - \frac{\lambda b }{\sqrt{c}} \leq \zeta.\]
    where the first inequality holds as $\beta_j \geq 0$.
    \item When $j \geq \lfloor \sqrt{c} \rfloor$. $\beta_{j+1} = 0$, thus 
    \begin{align*}
        \lambda g(x^*) - (\beta_j - \beta_{j+1}) x^* - (\lambda r + \beta_{j+1}) \frac{j}{\lambda d} &\stackrel{(a)}{\leq}   \lambda g(x^*) - (\lambda r + \beta_{j+1}) \frac{j}{\lambda d} \\
        & = \lambda g(x^*) - \lambda r \frac{j}{\lambda d} \\
        & \stackrel{(b)}{\leq}   \lambda g(x^*) - \lambda r \frac{\lfloor \sqrt{c} \rfloor}{\lambda d} \\
        & = \lambda g(x^*) - \lambda r \frac{\lfloor \sqrt{c} \rfloor}{\sqrt{c}} \frac{\sqrt{c}}{c}\frac{c}{\lambda d} \\
        & \stackrel{(c)}{\leq} \lambda g(x^*) - \frac{\lambda r x^*}{2\sqrt{c}}\\
        &  \leq \zeta,
    \end{align*}
    where $(a)$ holds as $\beta_j \geq 0$; $(b)$ holds as $j \geq \lfloor \sqrt{c} \rfloor$; $(c)$ holds as $\lfloor \sqrt{c} \rfloor/\sqrt{c} \geq 1/2$ for $c \geq 1$ and $x^* = c/(\lambda d)$.
\end{itemize}

\noindent\textbf{Constraint $\zeta \ge \lambda r\frac{c}{\lambda  d}+\beta_1 \frac{c}{\lambda d}$}

The right-hand-side of the constraint:

\[r\frac{c}{\lambda  d}+\beta_1 \frac{c}{\lambda d} = r\frac{c}{\lambda  d}+ \lambda \frac{b}{x^*} (1-\frac{1}{\sqrt{c}}) \frac{c}{\lambda d} = g(x^*) - \frac{\lambda b}{\sqrt{c}} \leq \zeta.\]

\noindent\textbf{Constraint $\zeta \ge \lambda b-\beta_cq$}

As $\beta_c = 0$, the constraint holds.

Therefore, the solution stated in \eqref{eq: lower bound static pricing q big} is a feasible solution. Therefore by weak duality, the performance loss of the optimal static policy is lower bounded by $\Omega(\sqrt{c})$ when choosing a static price $q \geq x^*$.

\subsubsection{Second case: $q < x^*$.\ } 
Consider the following solution
\begin{align}\label{eq: lower bound static pricing q small}
    \beta_j = \frac{\lambda r}{2}\max\left\{-\frac{j}{\sqrt{c}}, - \frac{\lfloor \sqrt{c} \rfloor}{\sqrt{c}}\right\}, \forall c \in [j], \zeta = \lambda g(x^*) - \frac{\lambda r x^*}{4\sqrt{c}}.
\end{align}

\noindent\textbf{Constraint $\zeta \ge \lambda r \frac{c-j}{\lambda d}+\lambda b -\beta_j q + \beta_{j+1} \frac{c-j}{\lambda  d}$}

The right-hand-side of the constraint can be written as 

\begin{align*}
    \lambda r \frac{c-j}{\lambda d}+\lambda b -\beta_j q + \beta_{j+1} \frac{c-j}{\lambda  d} & = \lambda r \frac{c}{\lambda d} - \lambda r \frac{j}{\lambda d}+ \lambda b - \beta_j q + \beta_{j+1} \frac{c}{\lambda d} - \beta_{j+1} \frac{j}{\lambda d}\\
    & \stackrel{(a)}{=}  \lambda g(x^*) - \beta_j q + \beta_{j+1} \frac{c}{\lambda d} - (\lambda r + \beta_{j+1}) \frac{j}{\lambda d}\\
    & \stackrel{(b)}{\leq} \lambda g(x^*) - (\beta_j - \beta_{j+1}) x^* - (\lambda r + \beta_{j+1}) \frac{j}{\lambda d},
\end{align*}
where $(a)$ holds as $\lambda g(x^*)  =\lambda rx^* + \lambda b, x^* = c/(\lambda d)$; $(b)$ holds as $q < x^* = c/(\lambda d), \beta \leq 0$.

\begin{itemize}
    \item When $j < \lfloor \sqrt{c} \rfloor$, $\beta_j - \beta_{j+1} = \lambda r/(2\sqrt{c})$, thus 
    \[\lambda g(x^*) - (\beta_j - \beta_{j+1}) x^* - (\lambda r + \beta_{j+1}) \frac{j}{\lambda d} \leq \lambda g(x^*) - (\beta_j - \beta_{j+1}) x^* = \lambda g(x^*) - \frac{\lambda r x^* }{2\sqrt{c}} \leq \zeta.\]
    where the first equality holds as $\beta_{j+1} \geq -\lambda r \lfloor \sqrt{c} \rfloor/\sqrt{c} \geq -\lambda r$.
    \item When $j \geq \lfloor \sqrt{c} \rfloor$. $\beta_{j+1} = \beta_j = -\lambda r\lfloor \sqrt{c} \rfloor/(2\sqrt{c})$, thus 
    \begin{align*}
        \lambda g(x^*) - (\beta_j - \beta_{j+1}) x^* - (\lambda r + \beta_{j+1}) \frac{j}{\lambda d} & = \lambda g(x^*) - (\lambda r + \beta_{j+1}) \frac{j}{\lambda d}\\
        & = \lambda g(x^*) - \left(\lambda r -\lambda r \frac{\lfloor \sqrt{c} \rfloor}{2\sqrt{c}} \right) \frac{j}{\lambda d} \\
        & \stackrel{(a)}{\leq}  \lambda g(x^*) - \left(\lambda r -\lambda r \frac{\lfloor \sqrt{c} \rfloor}{2\sqrt{c}} \right) \frac{\lfloor \sqrt{c} \rfloor}{\lambda d} \\
        & = \lambda g(x^*) - \lambda r\left( 1-\frac{\lfloor \sqrt{c} \rfloor}{2\sqrt{c}} \right) \frac{\lfloor \sqrt{c} \rfloor}{\sqrt{c}} \frac{\sqrt{c}}{c}\frac{c}{\lambda d}\\
        & \stackrel{(b)}{\leq}  \lambda g(x^*) - \frac{\lambda r x^*}{4\sqrt{c}}\\
        & = \zeta.
    \end{align*}
    where $(a)$ holds as $j \geq \lfloor \sqrt{c} \rfloor$; $(b)$ holds as $\lfloor \sqrt{c} \rfloor/\sqrt{c} \geq 1/2$ for $c \geq 1$ and $x^* = c/(\lambda d)$.
\end{itemize}

\noindent\textbf{Constraint $\zeta \ge \lambda r\frac{c}{\lambda  d}+\beta_1 \frac{c}{\lambda d}$}

As $\beta_1 <0$, the constraint holds.

\noindent\textbf{Constraint $\zeta \ge \lambda b-\beta_cq$}

The right-hand-side of the constraint:
\[\lambda b + \frac{\lambda r}{2}\frac{\lfloor \sqrt{c}\rfloor}{\sqrt{c}}q \leq \lambda b + \frac{\lambda r}{2}\frac{\lfloor \sqrt{c}\rfloor}{\sqrt{c}}x^* \leq \lambda b + \frac{\lambda r x^*}{2} < \zeta.\]

Therefore, the solution stated in \eqref{eq: lower bound static pricing q small} is a feasible solution. Therefore by weak duality, the performance loss of the optimal static policy is lower bounded by $\Omega(\sqrt{c})$ when choosing a static price $q < x^*$.
\Halmos\endproof

\section{Two-Price Policies}
\subsection{Proof of \Cref{lem: bound varsigma}}\label{apx: bound varsigma}
\proof{Proof.}
By \Cref{pro: balance equations}, the steady-state distribution under the two-price policy parameterized by $\tau \leq c/2, 0 <x_L < x_H < 1$ would satisfy
\[ 1 = \sum_{j=0}^c \pi_j = \pi_0\Bigg(1+\underbrace{\sum_{j=1}^\tau \frac{c!}{(\lambda d x_L)^j(c-j)!}}_{\varsigma_L}+\underbrace{\sum_{j=\tau+1}^c\frac{c!}{(\lambda d x_L)^\tau(\lambda d x_H)^{j-\tau}(c-j)!}}_{\varsigma_H}\Bigg). \] 
We now proceed to bound $\varsigma_L$ and $\varsigma_H$.  

\begin{align*}
    \varsigma_L & = \sum_{j=1}^\tau \frac{c!}{(\lambda d x_L)^j(c-j)!} = \sum_{j=1}^\tau \left(\frac{c}{\lambda d x_L}\right)^j \frac{c!}{c^j(c-j)!} \stackrel{(a)}{\leq}  \sum_{j=1}^\tau \left(\frac{c}{\lambda d x_L}\right)^j \\
    & = \sum_{j=1}^\tau (x^*/x_L)^j  \stackrel{(c)}{=} \frac{(x^*/x_L)^\tau - 1}{1-x_L/x^*} \leq \frac{(x^*/x_L)^\tau}{1-x_L/x^*},
\end{align*}
and
\begin{align*}
\varsigma_L &\stackrel{(a)}{\geq} \sum_{j=1}^\tau \left(\frac{c}{\lambda d x_L}\right)^j 
 \left(1-\frac{j^2}{2c}\right) \stackrel{(b)}{\geq} \left(1-\frac{\tau^2}{2c}\right) \sum_{j=1}^\tau  \left(\frac{c}{\lambda d x_L}\right)^j \\
 & = \left(1-\frac{\tau^2}{2c}\right) \sum_{j=1}^\tau  (x^*/x_L)^j \stackrel{(c)}{=} \left(1-\frac{\tau^2}{2c}\right)  \frac{(x^*/x_L)^\tau - 1}{1-x_L/x^*},
\end{align*}
where $(a)$ follows from \Cref{lem: bound product term}; $(b)$ holds as $j \leq \tau$ for $j \in [\tau]$ ; $(c)$ follows from geometric sum.

Similarly, for the second term we have
\begin{align*}
    \varsigma_H & = \sum_{j=\tau+1}^c\frac{c!}{(\lambda dx_L)^\tau(\lambda dx_H)^{j-\tau}(c-j)!} \\
    & = \sum_{j=\tau+1}^c  \left(\frac{c}{\lambda dx_L}\right)^\tau \left(\frac{c}{\lambda dx_H}\right)^{j-\tau} \frac{c!}{c^j(c-j)!} \\ 
    & \stackrel{(d)}{\geq}  \left(\frac{c}{\lambda dx_L}\right)^\tau \sum_{j=\tau+1}^{2\tau}\left(\frac{c}{\lambda dx_H}\right)^{j-\tau} \frac{c!}{c^j(c-j)!} \\ 
    & \stackrel{(e)}{\geq}  \left(\frac{c}{\lambda dx_L}\right)^\tau 
 \sum_{j=\tau+1}^{2\tau} \left(\frac{c}{\lambda dx_H}\right)^{j-\tau} \left(1-\frac{j^2}{2c}\right) \\
    & \stackrel{(f)}{\geq} \left(\frac{c}{\lambda dx_L}\right)^\tau  \left(1-\frac{2\tau^2}{c}\right)\sum_{j=\tau+1}^{2\tau}   \left(\frac{c}{\lambda dx_H}\right)^{j-\tau}  \\
    & = (x^*/x_L)^\tau  \left(1-\frac{2\tau^2}{c}\right) \sum_{j=\tau+1}^{2\tau}  (x^*/x_H)^{j-\tau}\\
    &\stackrel{(g)}{=} (x^*/x_L)^\tau    \left(1-\frac{2\tau^2}{c}\right) \frac{1-(x^*/x_H)^\tau}{x_H/x^* - 1}, 
\end{align*}
where $(d)$ follows from $\tau \leq c/2$; $(e)$ follows from \Cref{lem: bound product term}; $(f)$ holds as $j \geq \tau$ for $\tau+1 \leq j\leq 2\tau$; $(g)$ follows from geometric sum.

Putting everything together, 
\begin{align*}
    \varsigma_L-\varsigma_H & \leq \frac{(x^*/x_L)^\tau - 1}{1-x_L/x^*} - (x^*/x_L)^\tau \left(1-\frac{2\tau^2}{c}\right) \frac{1-(x^*/x_H)^\tau}{x_H/x^* - 1} \\
    & \leq \frac{(x^*/x_L)^\tau }{1-x_L/x^*} - (x^*/x_L)^\tau \left(1-\frac{2\tau^2}{c}\right) \frac{1-(x^*/x_H)^\tau}{x_H/x^* - 1}\\
    & \stackrel{(h)}{=} \frac{(x^*/x_L)^\tau }{1-x_L/x^*} - (x^*/x_L)^\tau \left(1-\frac{2\tau^2}{c}\right) \frac{1-(x^*/x_H)^\tau}{1-x_L/x^*} \\
    & = \frac{(x^*/x_L)^\tau }{1-x_L/x^*} \left( 1 - (1-2\tau^2/c) (1-(x^*/x_H)^\tau)\right)\\
    & \stackrel{(i)}{\leq}  \frac{(x^*/x_L)^\tau }{1-x_L/x^*} \left(\frac{2\tau^2}{c}+(x^*/x_H)^\tau \right),
\end{align*}
where $(h)$ holds as $x^*-x_L = x_H - x^*$; $(i)$ holds as $(1-(1-x)(1-y))\leq x+y$ for $x,y\geq 0$. 
\Halmos\endproof

\subsection{Proof of \Cref{lem: property of tau}}\label{apx: property of tau}
\proof{Proof.}
\begin{itemize}
    \item To see \eqref{eq: equiva tau c}. 
    By definition of $\tau$,  
\[(1+\delta(c)/x^*)^{\tau} = (1+\delta(c)/x^*)^{\frac{\log{c}}{\log(1+\delta(c)/x^*)} + \omega} = c \cdot (1+\delta(c)/x^*)^\omega, \]
with 
\[\omega = \tau -  \frac{\log{c}}{\log(1+\delta(c)/x^*)} \in [0,1).\]

By definition, $\delta(c) \in (0,1)$, thus
\[ (1+\delta(c)/x^*)^\omega \geq 1^\omega = 1,\quad (1+\delta(c)/x^*)^\omega \leq (1+1/x^*)^\omega \leq 1+1/x^*. \]
    \item To see \eqref{eq: bound tau}, note that for $y >0$, the following inequality holds
    \[ \frac{y}{1+y} \leq \log(1+y) \leq y.\]
    This is because 
    \[\frac{y}{y+1} = \int_0^y \frac{1}{y+1} \, dt \leq \int_0^y \frac{1}{t+1} \, dt = \log{(y+1)} = \int_0^y \frac{1}{t+1} \, dt  \leq \int_0^y 1 \, dt  = y.\]
    Replacing $y$ with $\delta(c)/x^*$ which is positive by definition, we have 
    \[\frac{\delta(c)/x^*}{1+\delta(c)/x^*} \leq \log(1+\delta(c)/x^*) \leq \delta(c)/x^*.\]
    Therefore, 
    \[\tau \geq \frac{\log{c}}{\log(1+\delta(c)/x^*)} \geq \frac{\log{c}}{\delta(c)/x^*},\]
    and 
    \[\tau \leq \frac{\log{c}}{\log(1+\delta(c)/x^*)} +1 \leq \frac{\log{c} (1+\delta(c)/x^*)}{\delta(c)/x^*} + 1 = \frac{\log{c}}{\delta(c)/x^*} + \log{c} + 1  \leq \frac{3\log{c}}{\delta(c)/x^*},\]
    where the last inequality holds as 1) $\delta(c) \leq x^* \implies log(c) \leq log(x)/(\delta(c)/x^*)$ and 2) $c\geq 4 \implies 1 \leq \log{c}$.
    
    \item To see \eqref{eq: minus tau bigger}, note that 
    \[1 \geq (1-(\delta(c)/x^*)^2)^\tau = (1-\delta(c)/x^*)^\tau (1+\delta(c)/x^*)^\tau \iff (1-\delta(c)/x^*)^{-\tau} \geq (1+\delta(c)/x^*)^\tau.\]
\end{itemize}

\Halmos\endproof

\subsection{Proof of \Cref{lem: bound tau2c}}\label{apx: bound tau2c}
\proof{Proof.}
By \eqref{eq: bound tau} and $\delta(c) = c^{-1/(\alpha+1)}$ when $\alpha \in (1,\infty)$, 
\begin{align*}
   \lim_{c\to\infty} \frac{\tau^2}{c} \leq \frac{(\log{c}x^* c^{1/(\alpha+1)} + \log{c} + 1)^2}{c} = 0.
\end{align*}

When $\alpha = \infty$, $\delta(c) = \varepsilon$
\[\lim_{c\to\infty} \frac{\tau^2}{c} \leq \frac{(\log{c}x^*/\varepsilon + \log{c} + 1)^2}{c} = 0.\]
\Halmos\endproof

\section{Proof of Lower Bound}

\subsection{Proof of \Cref{pro: feasible solution to dual alpha 1}}\label{apx: feasible solution to dual alpha 1}
\proof{Proof.}
Note that the feasible solution we use is similar to the one used in the proof for $\Omega(\sqrt{c})$ performance loss for static policies in \Cref{apx: lower bound statis pricing}.

By construction, we have that $\alpha_j ,\beta_j \ge 0$ for $j \in [c]$. To prove feasibility, we need to show that the dual constraints hold. Let $\epsilon = R-r$, which is independent of $c$. 

\noindent{\textbf{Constraint $\zeta \ge Rx^*\alpha_1+rx^*\beta_1$.}}

As $c \geq 1/(R-r)$, $\alpha_1 = \lambda (1-1/\sqrt{c \epsilon})$, we have 
\[Rx^* \alpha_1 + r x^* \beta_1 \stackrel{(a)}{=} \lambda r x^* +  \epsilon x^* \alpha_1 = \lambda r x^* + \lambda \epsilon x^* - \frac{\lambda \epsilon x^*}{\sqrt{c\epsilon}}  = \lambda R x^* - \lambda x^* \sqrt{\frac \epsilon c} \stackrel{(b)}{\leq }  \lambda g(x^*) - \lambda x^*\sqrt{\frac \epsilon c}  = \zeta, \]
where $(a)$ holds as $\alpha_1 + \beta_1 = \lambda$ by definition; $(b)$ holds as $g(x) \leq Rx + g(x^*) - R x^*$ for all $x\in [0,1]$ by definition, thus $0 = g(0) \leq g(x^*)- R x^*$.

\noindent\textbf{Constraint $\zeta \ge \lambda g(x^*)-Rx^*\alpha_j-rx^*\beta_j+R\frac{c-j}{\lambda d}\alpha_{j+1}+r\frac{c-j}{\lambda d}\beta_{j+1}$}
The right-hand-side of the constraint can be written as:
\begin{align*}
    & \lambda g(x^*)-Rx^*\alpha_j-rx^*\beta_j+R\frac{c-j}{\lambda d}\alpha_{j+1}+r\frac{c-j}{\lambda d}\beta_{j+1}  \\
    & \quad \stackrel{(a)}{=} \lambda g(x^*)-(r+\epsilon)x^*\alpha_j-rx^*\beta_j+(r+\epsilon)\frac{c-j}{\lambda d}\alpha_{j+1}+r\frac{c-j}{\lambda d}\beta_{j+1} \\
    & \quad \stackrel{(b)}{=}  \lambda g(x^*) - \lambda r x^* - \epsilon x^* \alpha_j + \lambda r \frac{c-j}{\lambda d} + \epsilon \frac{c-j}{\lambda d} \alpha_{j+1}\\
    & \quad \stackrel{(c)}{=}  \lambda  g(x^*) - \lambda r \frac{j}{\lambda d} + \epsilon x^* (\alpha_{j+1} - \alpha_j) - \epsilon \frac{j}{\lambda d} \alpha_{j+1},
\end{align*}
where $(a)$ holds as $R  = r+\epsilon$ by definition; $(b)$ holds as $\alpha_j + \beta_j = \lambda$ by definition; $(c)$ holds as $x^* = c/(\lambda d)$ by definition.

We proceed to show the solution is feasible by discussing the following three cases.
\begin{itemize}
    \item $j < \lfloor \sqrt{c\epsilon} \rfloor$: $\alpha_{j+1} - \alpha_j = - \lambda/\sqrt{c\epsilon}$. Thus
    \[\lambda  g(x^*) - \lambda r \frac{j}{\lambda d} + \epsilon x^* (\alpha_{j+1} - \alpha_j)- \epsilon \frac{j}{\lambda d} \alpha_{j+1} \leq \lambda g(x^*) +  \epsilon x^* (\alpha_{j+1} - \alpha_j) \leq \lambda g(x^*) - \lambda x^* \sqrt{\frac{\epsilon}{c}}.\]
    \item $ j \geq \lfloor \sqrt{c\epsilon} \rfloor$: $\alpha_{j+1} = 0$. Thus
    \begin{align*}
        \lambda  g(x^*) - \lambda r \frac{j}{\lambda d} + \epsilon x^* (\alpha_{j+1} - \alpha_j)- \epsilon \frac{j}{\lambda d} \alpha_{j+1} &=  \lambda  g(x^*) - \lambda r \frac{\lfloor \sqrt{c\epsilon} \rfloor}{\lambda d} - \epsilon x^*  \alpha_j\\
        & \leq \lambda g(x^*) - \lambda r \frac{\lfloor \sqrt{c\epsilon} \rfloor}{\lambda d} \\
        & = \lambda g(x^*) - \lambda r \frac{\lfloor \sqrt{c\epsilon} \rfloor}{\sqrt{c\epsilon} } \frac{\sqrt{c\epsilon} }{c} \frac{c}{\lambda d} \\
        & \leq \lambda g(x^*) - \lambda r x^* \frac{1}{2} \sqrt{\frac{\epsilon}{c}} \\
        & \leq \zeta.
    \end{align*}
    where the second inequality holds as $c \geq 1/(R-r) \implies c \epsilon >1$, thus 
    ${\lfloor \sqrt{c\epsilon} \rfloor}/{\sqrt{c\epsilon} } \geq 1/2$.
\end{itemize}
\Halmos\endproof

\subsection{Proof of Lower Bound {for} $\alpha = \infty$}\label{apx: lower bound infty}

\begin{figure}
    \centering

    \begin{subfigure}[b]{0.3\textwidth}
        \centering
\begin{tikzpicture}
    \begin{axis}[
        axis lines = middle,
        xlabel = $x$,
        ylabel = $g(x)$,
        ymin=0, ymax=1,
        xmin=0, xmax=1.1,
        xtick={0.001, 0.2, 1},
        xticklabels={0, $x^*$, 1},
        ytick=\empty,
        width=7cm,
        height=5.5cm,
        legend style={legend columns=2},
        every axis x label/.style={at={(current axis.right of origin)}, anchor=west},
        every axis y label/.style={at={(current axis.north west)}, above=2mm},
        ]
        \addplot[domain=0:0.3, samples=100, blue, thick, opacity=0.8] {x};
        \addlegendentry{$g$}
        \addplot[domain=0:0.45, samples=100, red, thick, dashed, opacity=0.8] {x};
        \addlegendentry{$r_1x$}
        \addplot[domain=0.1:0.9, samples=100, red, thick, dotted, opacity=0.8] {0.3 + 0.5*(x - 0.3)};
        \addlegendentry{$r_2x + b_2$}
        \addplot[domain=0.55:1, samples=100, red, thick, dash dot, opacity=0.8] {0.5};
        \addlegendentry{$b_3$}
        \addplot[domain=0.3:0.7, samples=100, blue, thick, opacity=0.8] {0.3 + 0.5*(x - 0.3)};
        \addplot[domain=0.7:1, samples=100, blue, thick, opacity=0.8] {0.5};
        \draw[densely dotted] (axis cs:0.2,0) -- (axis cs:0.2,0.2);
    \end{axis}
\end{tikzpicture}
        \caption{Case 1}
        \label{fig: g infity case 1}
    \end{subfigure}
    \hfill
    \begin{subfigure}[b]{0.3\textwidth}
        \centering
\begin{tikzpicture}
    \begin{axis}[
        axis lines = middle,
        xlabel = $x$,
        ylabel = $g(x)$,
        ymin=0, ymax=1,
        xmin=0, xmax=1.1,
        xtick={0.001, 0.5, 1},
        xticklabels={0, $x^*$, 1},
        ytick=\empty,
        width=7cm,
        height=5.5cm,
        legend style={legend columns=2},
        every axis x label/.style={at={(current axis.right of origin)}, anchor=west},
        every axis y label/.style={at={(current axis.north west)}, above=2mm},
        ]
        \addplot[domain=0:0.3, samples=100, blue, thick, opacity=0.8] {x};
        \addlegendentry{$g$}
        \addplot[domain=0:0.45, samples=100, red, thick, dashed, opacity=0.8] {x};
        \addlegendentry{$r_1x$}
        \addplot[domain=0.1:0.9, samples=100, red, thick, dotted, opacity=0.8] {0.3 + 0.5*(x - 0.3)};
        \addlegendentry{$r_2x + b_2$}
        \addplot[domain=0.55:1, samples=100, red, thick, dash dot, opacity=0.8] {0.5};
        \addlegendentry{$b_3$}
        \addplot[domain=0.3:0.7, samples=100, blue, thick, opacity=0.8] {0.3 + 0.5*(x - 0.3)};
        \addplot[domain=0.7:1, samples=100, blue, thick, opacity=0.8] {0.5};
        \draw[densely dotted] (axis cs:0.5,0) -- (axis cs:0.5,0.4);
    \end{axis}
\end{tikzpicture}
        \caption{Case 2}
        \label{fig: g infity case 2}
    \end{subfigure}
    \hfill
    \begin{subfigure}[b]{0.3\textwidth}
        \centering
\begin{tikzpicture}
    \begin{axis}[
        axis lines = middle,
        xlabel = $x$,
        ylabel = $g(x)$,
        ymin=0, ymax=1,
        xmin=0, xmax=1.1,
        xtick={0.001, 0.8, 1},
        xticklabels={0, $x^*$, 1},
        ytick=\empty,
        width=7cm,
        height=5.5cm,
        legend style={legend columns=2},
        every axis x label/.style={at={(current axis.right of origin)}, anchor=west},
        every axis y label/.style={at={(current axis.north west)}, above=2mm},
        ]
        \addplot[domain=0:0.3, samples=100, blue, thick, opacity=0.8] {x};
        \addlegendentry{$g$}
        \addplot[domain=0:0.45, samples=100, red, thick, dashed, opacity=0.8] {x};
        \addlegendentry{$r_1x$}
        \addplot[domain=0.1:0.9, samples=100, red, thick, dotted, opacity=0.8] {0.3 + 0.5*(x - 0.3)};
        \addlegendentry{$r_2x + b_2$}
        \addplot[domain=0.55:1, samples=100, red, thick, dash dot, opacity=0.8] {0.5};
        \addlegendentry{$b_3$}
        \addplot[domain=0.3:0.7, samples=100, blue, thick, opacity=0.8] {0.3 + 0.5*(x - 0.3)};
        \addplot[domain=0.7:1, samples=100, blue, thick, opacity=0.8] {0.5};
        \draw[densely dotted] (axis cs:0.8,0) -- (axis cs:0.8,0.5);
    \end{axis}
\end{tikzpicture}
        \caption{Case 3}
        \label{fig: g infity case 3}
    \end{subfigure}
    \caption{Illustrations of $g$ under three cases, $\alpha = \infty$. In this specific example, $g$ and its approximation $\acute{g}$ equals each other, so we plot functions $r_1x, r_2x + b_2, b_3$.}
    \label{fig: illus of g 3 cases}
\end{figure}

\proof{Proof.}
The proof is similar to the one of \Cref{thm: degenerate case}: upper bound $g$ using a piecewise linear $\acute{g}$, and use primal-dual theory to upper bound the reward rate of the relaxed problem. 

When $g$ satisfies \Cref{assmp: condition for lower bound} with $\alpha = \infty$, there must exist $r_1, r_2, b_2, b_3 >0, r_1 > r_2 >0$ and $b_3 > b_2 >0$ such that 
$$
g(x) \leq \acute{g}(x) = \min\{r_1 x, r_2 x + b_2,b_3\}, \forall x\in [0,1], g(x^*) = \acute{g}(x^*), \text{ and } x^*\neq \frac{b_2}{r_1-r_2}, x^* \neq \frac{b_3-b_2}{r_2}.
$$
{We illustrate this in \Cref{fig: illus of g 3 cases} in an example where $\acute{g}$ and $g$ are equal, showing the three possibilities for where $x^*$ can lie.}

The last condition restricts the piecewise linear function $\acute{g}(x)$ to be differentiable at $x^*$, as the non-differentiable case has already been discussed in \Cref{thm: degenerate case}.

Therefore, by \Cref{thm: optsd to lpsd}, we can analyze the following LP to obtain an upper bound on the optimal reward under any stock-dependent policy.
\begin{align*}
    \max_{\bm{\pi}\geq \bm{0}} & \sum_{j=1}^c \lambda g_j\\
   \text{s.t. }  & \sum_{j=0}^c\pi_j = 1,\\
   & g_j \leq r_1\frac{c-j+1}{\lambda d}\pi_{j-1}, &\forall j \in [c], \\
   & g_j \leq r_2\frac{c-j+1}{\lambda d}\pi_{j-1} + b_2 \pi_j, &\forall j \in [c], \\
   & g_j \leq b_3\pi_j, &\forall j \in [c].
\end{align*}
Its dual problem is
\begin{align*}
\min\quad &\zeta 
\\ \text{s.t.}\quad &\alpha_j+\beta_j + \gamma_j= \lambda, & \forall j \in [c], \\ 
&\zeta \ge \beta_jb_2+ \gamma_j b_3 +\alpha_{j+1}r_1\frac{c-j}{\lambda d} + \beta_{j+1}r_2\frac{c-j}{\lambda d}, &\forall j \in [c-1],\\ 
& \zeta \ge \beta_c b_2 + \gamma_c b_3, \\
&\zeta \ge \alpha_1r_1x^* + \beta_1 r_2 x^*,\\ 
&\alpha_j,\beta_j, \gamma_j \ge 0, &\forall j \in [c].
\end{align*}

\noindent\textbf{Case 1: $g(x^*) = r_1x^*$.}

See \Cref{fig: g infity case 1} for an illustration. 
This is the corner case when $\partial g(x^*) = \{r\}$ for $x \in [0,x^* + \varepsilon]$ for some constant $r>0, \varepsilon >0$, and the optimal static policy achieves a performance loss of $\mathcal{O}(1)$ (\Cref{pro: scarce supply helps}). We proceed to show that it is in fact tight.

Consider the following solution:
\begin{align}\label{eq: infinity non-flat case feasible solution}
    \zeta = \lambda r_1 x^* - \frac{r_1}{2d},  \alpha_j = \lambda - \frac{r_1}{2db_2}, \beta_j = \frac{r_1}{2db_2}, \gamma_j = 0, \forall j \in [c].
\end{align}
To see this is feasible, 
\begin{itemize}
    \item \textbf{Constraint} $\zeta \ge \beta_j b_2+\alpha_{j+1}r_1{(c-j)}/{(\lambda d)} + \beta_{j+1} r_2 ({c-j})/({\lambda d})$ for $j \in [c-1]$: 
    
    The right-hand-side of the constraint can be written as 
    \begin{align*}
        \beta_j b_2+\alpha_{j+1}r_1{(c-j)}/{(\lambda d)} + \beta_{j+1} r_2 \frac{c-j}{\lambda d} & \stackrel{(a)}{=} \beta b_2 + (\lambda -\beta) r_1 \frac{c-j}{\lambda d} + \beta r_2 \frac{c-j}{\lambda d}\\
        & \stackrel{(b)}{=} \lambda r_1 x^* - \lambda r_1 \frac{j}{\lambda d} - \beta (r_1-r_2)\frac{c-j}{\lambda d} + \beta b_2 \\
        & \stackrel{(c)}{\leq} \lambda r_1 x^* - \lambda r_1 \frac{j}{\lambda d} + \beta b_2 \\
        & \stackrel{(d)}{\leq} \lambda r_1 x^* - \lambda r_1 \frac{1}{\lambda d} + \beta b_2\\
        & \stackrel{(e)}{\leq} \lambda r_1 x^* - \frac{r_1}{d} +\frac{r_1}{2d} \\
        & = \zeta.
    \end{align*}
    where $(a)$ holds as $\alpha_j + \beta_j = \lambda$, $\beta_j$ is a constant for $j\in [c]$; $(b)$ holds as $x^* = c/(\lambda d)$ by definition; $(c)$ holds as $\beta >0$, $r_1 -r_2 >0$, $c-j \geq 0$; $(d)$ holds as $j \geq 1$; $(e)$ holds as $\beta = r_1/(2 d b_2)$.

    \item \textbf{Constraint} $\zeta \ge \beta_c b_2 + \gamma_c b_3$:

     The right-hand-side of the constraint can be written as
    \[\beta_c b_2 + \gamma_c b_3= \beta_c b_2 = \frac{r_1}{2 d} < \zeta.\]

    \item \textbf{Constraint} $\zeta \ge \alpha_1r_1x^* + \beta_1 r_2 x^*$:
    The right-hand-side of the constraint can be written as 
    \[ \alpha_1r_1x^* + \beta_1 r_2 x^* = (\lambda-\beta) r_1 x^* + \beta r_2 x^* = \lambda r_1 x^* - \beta x^*(r_1 - r_2) < \zeta,\]
    where the last inequality holds as $r_1 > r_2$ and $\beta >0$.
\end{itemize}
Therefore, the solution stated in \eqref{eq: infinity non-flat case feasible solution} is a feasible solution. Therefore the performance loss of the optimal stock-dependent policy is lower bounded by 
\[\frac{r_1}{2 d}= \Omega(1).\]

\noindent\textbf{Case 2: $g(x^*) = r_2x^* + b_2$.}

See \Cref{fig: g infity case 2} for an illustration. 
Consider the following solution:
\begin{align}\label{eq: infinity middle case feasible solution}
    & \zeta = \lambda (r_2 x^* + b_2) - \min\left\{ K, \frac{r_2}{d}\right\} \log{c}, K = \frac{\lambda}{c} \frac{b_2}{x^*(r_1-r_2)}(x^*(r_1-r_2) - b_2) \text{ is a constant,} \\
    & \alpha_1 = \lambda \frac{b_2}{x^*(r_1-r_2)} \frac{\log{c}}{c}, \alpha_j b_2 - \alpha_{j+1}x^*(r_1-r_2) = K \log{c}, \forall j \leq  \log{c}, \alpha_j = 0, \forall j > \log{c} + 1,  \notag\\
    & \alpha_j + \beta_j = \lambda,  \gamma_j = 0, \forall j \in [c]. \notag
\end{align}
To see this is feasible, we first show that $0\leq \alpha_j  \leq \lambda$ for all $j \in [c]$. As $\alpha_j b_2 - \alpha_{j+1}x^*(r_1-r_2) = K \log{c}, \forall j \leq  \log{c}$, this implies
\[\alpha_{j+1} =  \underbrace{\frac{b_2}{x^*(r_1-r_2)}}_{\eta}\alpha_j - \underbrace{\frac{K \log{c}}{x^*(r_1-r_2)}}_{\omega},\]
divide both sides by $\eta^{j+1}$,
\[\frac{\alpha_{j+1}}{\eta^{j+1}} =  \frac{\alpha_j}{\eta^j} - \frac{\omega}{\eta^{j+1}},\]
thus we have for $\tau \leq \log{c}+1$,
\begin{align}\label{eq: alpha tau func of alpha 1}
\frac{\alpha_{\tau}}{\eta^{\tau}} =  \frac{\alpha_1}{\eta} - \sum_{j=2}^\tau \frac{\omega}{\eta^{j}} = \frac{\alpha_1}{\eta} - \omega \frac{\eta^2 - \eta^{\tau+1}}{1-\eta},
\end{align}
where the last equality holds from geometric sum.

Note that as $g(x^*) = r_2x^* + b_2$, we have $x^* > b_2/(r_1-r_2)$, thus $x^*(r_1-r_2)>b_2$, i.e., $\eta <1, \omega >0$, therefore $\alpha_j$ are decreasing in $j$. Therefore, 
\[\alpha_j \leq \alpha_1 = \lambda \frac{b_2}{x^*(r_1-r_2)} \frac{\log{c}}{c} \leq \lambda\frac{\log{c}}{c} \leq \lambda, \]
where the last inequality holds as $\max_{c} \log{c}/c = 1/e <1$.

To see $\alpha_{j} \geq 0$, first by definition of $\eta, \omega, K$ and $\alpha_1$, 
\[ \frac{\omega}{1-\eta} = \frac{K\log{c}}{x^*(r_1-r_2) - b_2}  ={\lambda} \frac{b_2}{x^*(r_1-r_2)}\frac{\log{c}}{c} = \alpha_1. \]
Therefore, divided both sides by $\eta$:
\[\frac{\alpha_1}{\eta} = \frac{\omega}{1-\eta}\frac{1}{\eta} \geq \frac{\omega}{1-\eta} \eta^2 \geq \omega \frac{\eta^2 -  \eta^{\tau+1}}{1-\eta}.\]
where the last two inequalities hold as $0<\eta<1$. By \eqref{eq: alpha tau func of alpha 1}, this implies that $\alpha_j \geq 0$ for all $j\in [c]$.

We proceed to show the other constraints hold.
\begin{itemize}
    \item \textbf{Constraint} $\zeta \ge \beta_j b_2+\alpha_{j+1}r_1{(c-j)}/{(\lambda d)} + \beta_{j+1} r_2 ({c-j})/({\lambda d})$ for $j \in [c-1]$: 
    
    The right-hand-side of the constraint can be written as 
    \begin{align*}
        \beta_j b_2+\alpha_{j+1}r_1\frac{c-j}{\lambda d} + \beta_{j+1} r_2 \frac{c-j}{\lambda d} & {=} \beta_j b_2 + \alpha_{j+1} r_1 \frac{c}{\lambda d} + \beta_{j+1} r_2 \frac{c}{\lambda d} - \left(\alpha_{j+1} r_1 +  \beta_{j+1} r_2\right) \frac{j}{\lambda d} \\       
        & \stackrel{(a)}{\leq} \beta_j b_2 + \alpha_{j+1} r_1 x^* + \beta_{j+1} r_2 x^* - \lambda r_2 \frac{j}{\lambda d}\\ 
        & \stackrel{(b)}{=} (\lambda - \alpha_j) b_2 + \alpha_{j+1} r_1 x^* + (\lambda - \alpha_{j+1}) r_2 x^* - \lambda r_2 \frac{j}{\lambda d} \\
        & =  \lambda (r_2 x^* + b_2) - \alpha_j b_2 + \alpha_{j+1}x^*(r_1-r_2)- \lambda r_2 \frac{j}{\lambda d},
    \end{align*}
    where $(a)$ holds as $\alpha_j + \beta_j = \lambda$ and $r_1 > r_2$; $(b)$ holds as $\alpha_j + \beta_j = \lambda$.
    \begin{itemize}
        \item If $j \leq \log(c)$: 
        \begin{align*}
            \lambda (r_2 x^* + b_2) - \alpha_j b_2 + \alpha_{j+1}x^*(r_1-r_2)- \lambda r_2 \frac{j}{\lambda d} & \stackrel{(a)}{\leq} \lambda (r_2 x^* + b_2) - \alpha_j b_2 + \alpha_{j+1}x^*(r_1-r_2)\\
            & \stackrel{(d)}{=} \lambda (r_2 x^* + b_2) - K \log{c} \\
            & {\leq}\zeta.
        \end{align*}
       where $(a)$ follows as $\lambda r_2 j/(\lambda d) \geq 0$; $(b)$ follows from the definition of the $\alpha_j$ sequence.

        \item If $j > \log{c}$:
                \begin{align*}
            \lambda (r_2 x^* + b_2) - \alpha_j b_2 + \alpha_{j+1}x^*(r_1-r_2)- \lambda r_2 \frac{j}{\lambda d} & \stackrel{(a)}{\leq} \lambda (r_2 x^* + b_2) + \alpha_{j+1}x^*(r_1-r_2)- \lambda r_2 \frac{j}{\lambda d} \\
            & \stackrel{(b)}{=} \lambda (r_2 x^* + b_2) - \lambda r_2 \frac{j}{\lambda d}\\
            & \stackrel{(c)}{\leq} \lambda (r_2 x^* + b_2) -  r_2 \frac{\log{c}}{d} \\
            & \leq \zeta, 
        \end{align*}
        where $(a)$ holds as $\alpha_j b_2 \geq 0$; $(b)$ holds as $\alpha_j = 0$ for $j > \log{c} + 1$; $(c)$ holds as $j > \log{c}$.
    \end{itemize}

    \item \textbf{Constraint} $\zeta \ge \beta_c b_2 + \gamma_c b_3$:

     The right-hand-side of the constraint can be written as
    \[\beta_c b_2 + \gamma_c b_3= \beta_c b_2 = \lambda b_2 < \zeta.\]

    \item \textbf{Constraint} $\zeta \ge \alpha_1r_1x^* + \beta_1 r_2 x^*$:
    The right-hand-side of the constraint can be written as 
    \[ \alpha_1r_1x^* + \beta_1 r_2 x^* =  \alpha_1 r_1 x^* + (\lambda -\alpha_1) r_2 x^* = \lambda r_2 x^* + \alpha_1 x^*(r_1 - r_2),\]
    where the last inequality holds as $\alpha_1 x^*(r_1-r_2) < \lambda b_2$.
\end{itemize}
Therefore, the solution stated in \eqref{eq: infinity middle case feasible solution} is a feasible solution. Therefore the performance loss of the optimal stock-dependent policy is lower bounded by 
\[\min\left\{ K, \frac{r_2}{d}\right\} \log{c} = \Omega(\log{c}).\]

\noindent\textbf{Case 3: $g(x^*) = b_3$.}

See \Cref{fig: g infity case 3} for an illustration. 
This is the corner case when $\partial g(x^*) = \{0\}$, and the optimal static policy achieves a performance loss of $\mathcal{O}(e^{-c})$ (\Cref{pro: excess supply helps}). We proceed to show that it is in fact tight.

Consider the following solution to the dual: 
\begin{align}\label{eq: infinity flat case feasible solution}
\zeta = \lambda b_3 - (b_3 - b_2) \beta_c, \beta_1 = \frac{\lambda b_3}{2 r_2 x^*}, \beta_{j+1} = \frac{b_3-b_2}{2 r_2 x^*}\beta_j, \forall j\in [c-1], \beta_j + \gamma_j = \lambda, \alpha_j = 0, \forall j\in [c].    
\end{align}

To see this is feasible, 
\begin{itemize}
    \item \textbf{Constraint} $\zeta \ge \beta_jb_2+\gamma_jb_3+\beta_{j+1}r_2{(c-j)}/{(\lambda d)}$ for $j \in [c-1]$: 
    
    The right-hand-side of the constraint can be written as 
    \begin{align}\label{eq: infinity lower bound constraint}
        \beta_jb_2+\gamma_jb_3+\beta_{j+1}r_2\frac{c-j}{\lambda d} & \stackrel{(a)}{=} \lambda b_3 - (b_3-b_2)\beta_j + \beta_{j+1} r_2 \frac{c-j}{\lambda d}  \leq \lambda b_3 - (b_3-b_2)\beta_j + \beta_{j+1} r_2 \frac{c}{\lambda d} \notag  \\
        & \stackrel{(b)}{=} \lambda b_3 - (b_3-b_2)\beta_j + \beta_{j+1} r_2 x^* \stackrel{(c)}{=} 
        \lambda b_3 - \beta_{j+1} r_2 x^* \notag \\
        & \stackrel{(d)}{\leq} \lambda b_3 - r_2 x^* \beta_c  \stackrel{(e)}{\leq} \lambda b_3 - (b_3-b_2) \beta_c = \zeta.
    \end{align}
    where $(a)$ holds as $\beta_j + \alpha_j = \lambda$ for $j \in [c]$; $(b)$ holds as $x^* = c/(\lambda d)$ by definition; $(c)$ holds as $2 r_2 x^* \beta_{j+1} = (b_3-b_2)\beta_j$ by definition; $(d)$ holds as $g(x^*) = b_3$, we must have $(b_3-b_2)/ r_2 < x^*$, thus $\beta_{j+1} < \beta_j$; $(e)$ holds as $(b_3-b_2)/r_2 < x^*$.
    \item \textbf{Constraint} $\zeta \ge \beta_c b_2 + \gamma_c b_3$:

     The right-hand-side of the constraint can be written as
    \[\beta_c b_2 + \gamma_c b_3 = \lambda b_3 - (b_3 - b_2)\beta_c = \zeta.\]

    \item \textbf{Constraint} $\zeta \ge \beta_1r_2x^*$:
    The right-hand-side of the constraint can be written as 
    \[ \beta_1 r_2 x^* \stackrel{(a)}{=} \lambda b_3 - \beta_1 r_2 x^* \stackrel{(a)}{<} \lambda b_3 - \beta_2 r_2 x^* \leq \zeta, \]
    where $(a)$ holds as $\beta_1 = {\lambda b_3}/({2 r_2 x^*})$; $(b)$ holds as $\beta_1 > \beta_2$; $(c)$ holds as $\lambda b_3 - \beta_2 r_2 x^* $ is the right-hand-side of the constraint for $j = 2$, which is bounded by $\zeta$ as shown in \eqref{eq: infinity lower bound constraint}.
\end{itemize}
Therefore, the solution stated in \eqref{eq: infinity flat case feasible solution} is a feasible solution. Therefore the performance loss of the optimal stock-dependent policy is lower bounded by 
\[\lambda (b_3-b_2)\beta_c = \lambda (b_3 - b_2) \frac{\lambda b_3}{2 r_2 x^*} \left( \frac{b_3-b_2}{2 r_2 x^*} \right)^{c-1} = \Omega(c e^{-c}).\]
where the last inequality holds as $({b_3-b_2})/{2 r_2 x^*}<1$.\Halmos \endproof

\section{Discussion on Assumptions \ref{assmp: condition for upper bound} and \ref{assmp: condition for lower bound}}\label{apx: two assumptions}
\begin{lemma}\label{lem: sufficient condition of assmp ub}
A sufficient condition for \Cref{assmp: condition for upper bound} is $|\partial g(x) - \partial g(x^*)| \leq k_1|x-x^*|^{\alpha-1}$ for all $x\in[x^*-\varepsilon,x^*+\varepsilon]$.
\end{lemma}
\proof{Proof.}
If $x > x^*$,
\begin{align*}
    g(x) & \stackrel{(a)}{\geq} g(x^*)+\partial g(x)(x-x^*) \\
    & \stackrel{(b)}{\geq} g(x^*) + (\partial g(x^*)-k_1(x-x^*)^{\alpha-1})(x-x^*)\\
    &= g(x^*) + \partial g(x^*)(x-x^*) - k_1(x-x^*)^\alpha, 
\end{align*}
where $(a)$ follows as $g$ is concave, i.e., $g(x^*)\leq g(x) + \partial g(x) (x^*-x)$; $(b)$ follows as $ \partial g(x^*)-\partial g(x) \leq k_1(x-x^*)^{\alpha-1}$. 
The proof for $x < x^*$ follows similarly. \Halmos \endproof
\begin{remark}
For welfare maximization, $g(x) = \int_{1-x}^1 F^{-1}(x)\, dx$ and $\partial g(x) = F^{-1}(1-x)$. Then by \Cref{lem: sufficient condition of assmp ub}, a sufficient condition for Assumption \ref{assmp: condition for upper bound} is 
\[ |F^{-1}(1-q) - F^{-1}(1-x^*)| \leq k_1|x^* - q|^{\alpha-1}, \]
which is a local version of the cluster density requirement of $(\beta,\epsilon_0)$-clustered distribution (Definition 1) in \cite{besbes2022multi}, with $\alpha-1 = 1/(\beta+1)$ for $\alpha\in(1,2]$. One can refer to \citet{besbes2022multi} for an example of a valuation distribution that leads to a function $g$ for which the best-possible value of $\alpha$ lies strictly between 1 and 2.
\end{remark}

\begin{proposition}\label{pro: 2continuous differentiable g}
If $g$ is twice continuously differentiable, then the optimal performance loss of stock-dependent policies is $\tilde{\mathcal{O}}(c^{-2/3})$.
\end{proposition}
\proof{Proof.}
If $g$ is twice continuously differentiable, then the first-order derivative $g'(x)$ is continuously differentiable and thus Lipschitz continuous: there exists $L>0$ such that
\[|g'(x) - g'(x^*)| \leq L|x-x^*|,\]
which implies \Cref{assmp: condition for upper bound} with $\alpha = 2$ by \Cref{lem: sufficient condition of assmp ub}.
\Halmos\endproof

\begin{proposition}
For welfare maximization, if the WTP density function $f$ is continuous and
positive almost everywhere over bounded support (i.e., $F$ is strictly increasing), then the optimal performance loss of stock-dependent policies is $\tilde{\Theta}(c^{-2/3})$.
\end{proposition}
\proof{Proof.}
For welfare maximization, $g(x) = \int_{1-x}^1 F^{-1}(v)\, dv$, thus $g'(x) = F^{-1}(1-x)$. When $f$ is continuous and $F$ is strictly increasing, both $F'(x) = f(x)$ and $(F^{-1})'(x) = 1/f(x)$ are bounded. Therefore, there must exist $L, L' >0$ such that 
\begin{align}
    & |F(x) - F(x^*)| \leq L|x-x^*| \label{eq: F lipschitz} \\
    & |F^{(-1)}(x) - F^{(-1)}(x^*)| \leq L' |x-x^*|. \label{eq: F-1 lipschitz}
\end{align}

Note that $F$ is one-to-one---$F$ is $L$-Lipschitz continuous if and only if $F^{-1}$ is $1/L$-Lipschitz continuous. (To see this: let $F^{-1}(1-x_1) = y_1, F^{-1}(1-x_2) = y_2$, thus $1-x_1 = F(y_1)$ and $1-x_2 = F(y_2)$ uniquely determined.) Then inequality \eqref{eq: F lipschitz} implies
\begin{align}
    |1-x - (1-x^*)| \leq L |F^{-1}(1-x)-F^{-1}(1-x^*)|. \label{eq: F-1 inverse lipschitz}
\end{align}

Combing inequality \eqref{eq: F-1 lipschitz} and \eqref{eq: F-1 inverse lipschitz}, by \Cref{remark: balseiro2022mechanism}, this implies \Cref{assmp: condition for upper bound} and \Cref{assmp: condition for lower bound} with $\alpha = 2$, then the bound on performance loss follows. \Halmos\endproof

\begin{lemma}\label{lem: gradient condition}
A necessary condition for \Cref{assmp: condition for lower bound} is $|\partial g(x^*) - \partial g(x)| \geq k_2|x-x^*|^{\alpha-1}$ for all $x\in[x^*-\varepsilon,x^*+\varepsilon]$.
\end{lemma}
\proof{Proof.}
\[g(x^*) + \partial g(x)(x-x^*) \stackrel{(a)}{\leq} g(x) \stackrel{(b)}{\leq} g(x^*) + \partial g(x^*)(x-x^*)-k_2|x-x^*|^\alpha, \]
where $(a)$ is by concavity of $g$, i.e.,$g(x^*)\leq g(x)+\partial g(x) (x^*-x)$; $(b)$ is given. 
Thus $\partial g(x)(x-x^*) \leq  \partial g(x^*)(x-x^*)-k_2|x-x^*|^\alpha$. When $x\geq x^*$, we have $\partial g(x^*)-\partial g(x)\geq k_2(x-x^*)^{\alpha-1}$; when $x<x^*$, we have $\partial g(x) - \partial g(x^*) \geq k_2(x^*-x)^{\alpha-1}$. 

Therefore, we have $|\partial g(x^*) - \partial g(x)| \geq k_2|x-x^*|^{\alpha-1}$.
\Halmos \endproof

Finally, we use the following known result to show that a notion of Lipschitzness of gradients around $x^*$ is sufficient for Assumptions \ref{assmp: condition for upper bound} and \ref{assmp: condition for lower bound} to hold.

\begin{lemma} [Lemma 7 in \cite{balseiro2022mechanism}]\label{lem: balseiro2022mechanism}
Suppose that $R(v)$ is locally concave around $p^{\star}$ and let $\alpha \in(1,+\infty)$. Then there exists positive constants $\tilde{\kappa}_L, \tilde{\kappa}_U$ and a neighborhood of $p^{\star}, \mathcal{N}_{\tilde{\ell}}=\left(p^{\star}-\tilde{\ell}, p^{\star}+\tilde{\ell}\right) \subset(0, \bar{v})$, such that
$$
\tilde{\kappa}_L \cdot\left|v-p^{\star}\right|^\alpha \leq R\left(p^{\star}\right)-R(v) \leq \tilde{\kappa}_U \cdot\left|v-p^{\star}\right|^\alpha, \quad \forall v \in \mathcal{N}_{\tilde{\ell}},
$$
if and only if there exists positive constants $\kappa_L, \kappa_U$ and a neighborhood of $p^{\star}, \mathcal{N}_{\ell}=\left(p^{\star}-\ell, p^{\star}+\ell\right) \subset$ $(0, \bar{v})$, such that
$$
\kappa_L \alpha \cdot\left|v-p^{\star}\right|^\alpha \leq\left(p^{\star}-v\right) \cdot \dot{R}(v) \leq \kappa_U \alpha \cdot\left|v-p^{\star}\right|^\alpha, \quad \forall v \in \mathcal{N}_{\ell} .
$$
\end{lemma}
\begin{remark}\label{remark: balseiro2022mechanism}
The proof for \Cref{lem: balseiro2022mechanism} also holds for $\alpha = \infty$, and the result trivially holds for $\alpha = 1$.

Applying \cref{lem: balseiro2022mechanism} with $R(x) = g(x) - g'(x^*)(x-x^*)$ and $p^* = x^*$, we have that Assumptions \ref{assmp: condition for upper bound} and \ref{assmp: condition for lower bound} are equivalent to there existing positive constants $K_1, K_2 \geq 0$ such that 
\[ K_2 \cdot |x- x^*|^\alpha \leq (x^*-x)\cdot (g'(x) - g'(x^*)) \leq K_1 |x-x^*|^\alpha, \forall x\in [x^*-\varepsilon, x^*+\varepsilon]. \]
\end{remark}

\end{APPENDICES}






\end{document}